\input amstex
\documentstyle{amsppt}

\magnification1200
\hsize15.6truecm
\vsize22.8truecm
\def\ZeilBL{18}
\def\ZeilAV{17}
\def\ZeilAN{16}
\def\ZeilAM{15}
\def\WiZeAC{14}
\def\SlatAC{13}
\def\KratBD{12}
\def\KratBI{11}
\def\KratBH{10}
\def\KratBG{9}
\def\KratBF{8}
\def\KratAM{7}
\def\GrKPAA{6}
\def\GoJaAJ{5}
\def\GeViAA{4}
\def\GaRaAA{3}
\def\BoHPAA{2}
\def\BailAA{1}

\def\cl#1{\left\lceil#1\right\rceil}
\def\po#1#2{(#1)_#2}

\catcode`\@=11
\font@\twelverm=cmr10 scaled\magstep1
\font@\twelveit=cmti10 scaled\magstep1
\font@\twelvebf=cmbx10 scaled\magstep1
\font@\twelvei=cmmi10 scaled\magstep1
\font@\twelvesy=cmsy10 scaled\magstep1
\font@\twelveex=cmex10 scaled\magstep1

\newtoks\twelvepoint@
\def\twelvepoint{\normalbaselineskip15\p@
 \abovedisplayskip15\p@ plus3.6\p@ minus10.8\p@
 \belowdisplayskip\abovedisplayskip
 \abovedisplayshortskip\z@ plus3.6\p@
 \belowdisplayshortskip8.4\p@ plus3.6\p@ minus4.8\p@
 \textonlyfont@\rm\twelverm \textonlyfont@\it\twelveit
 \textonlyfont@\sl\twelvesl \textonlyfont@\bf\twelvebf
 \textonlyfont@\smc\twelvesmc \textonlyfont@\tt\twelvett
%
 \ifsyntax@ \def\big##1{{\hbox{$\left##1\right.$}}}%
  \let\Big\big \let\bigg\big \let\Bigg\big
 \else
  \textfont\z@=\twelverm  \scriptfont\z@=\tenrm  \scriptscriptfont\z@=\sevenrm
  \textfont\@ne=\twelvei  \scriptfont\@ne=\teni  \scriptscriptfont\@ne=\seveni
  \textfont\tw@=\twelvesy \scriptfont\tw@=\tensy \scriptscriptfont\tw@=\sevensy
  \textfont\thr@@=\twelveex \scriptfont\thr@@=\tenex
        \scriptscriptfont\thr@@=\tenex
  \textfont\itfam=\twelveit \scriptfont\itfam=\tenit
        \scriptscriptfont\itfam=\tenit
  \textfont\bffam=\twelvebf \scriptfont\bffam=\tenbf
        \scriptscriptfont\bffam=\sevenbf
  \setbox\strutbox\hbox{\vrule height10.2\p@ depth4.2\p@ width\z@}%
  \setbox\strutbox@\hbox{\lower.6\normallineskiplimit\vbox{%
        \kern-\normallineskiplimit\copy\strutbox}}%
 \setbox\z@\vbox{\hbox{$($}\kern\z@}\bigsize@=1.4\ht\z@
 \fi
 \normalbaselines\rm\ex@.2326ex\jot3.6\ex@\the\twelvepoint@}

\font@\fourteenrm=cmr10 scaled\magstep2
\font@\fourteenit=cmti10 scaled\magstep2
\font@\fourteensl=cmsl10 scaled\magstep2
\font@\fourteensmc=cmcsc10 scaled\magstep2
\font@\fourteentt=cmtt10 scaled\magstep2
\font@\fourteenbf=cmbx10 scaled\magstep2
\font@\fourteeni=cmmi10 scaled\magstep2
\font@\fourteensy=cmsy10 scaled\magstep2
\font@\fourteenex=cmex10 scaled\magstep2
\font@\fourteenmsa=msam10 scaled\magstep2
\font@\fourteeneufm=eufm10 scaled\magstep2
\font@\fourteenmsb=msbm10 scaled\magstep2
\newtoks\fourteenpoint@
\def\fourteenpoint{\normalbaselineskip15\p@
 \abovedisplayskip18\p@ plus4.3\p@ minus12.9\p@
 \belowdisplayskip\abovedisplayskip
 \abovedisplayshortskip\z@ plus4.3\p@
 \belowdisplayshortskip10.1\p@ plus4.3\p@ minus5.8\p@
 \textonlyfont@\rm\fourteenrm \textonlyfont@\it\fourteenit
 \textonlyfont@\sl\fourteensl \textonlyfont@\bf\fourteenbf
 \textonlyfont@\smc\fourteensmc \textonlyfont@\tt\fourteentt
%
 \ifsyntax@ \def\big##1{{\hbox{$\left##1\right.$}}}%
  \let\Big\big \let\bigg\big \let\Bigg\big
 \else
  \textfont\z@=\fourteenrm  \scriptfont\z@=\twelverm  \scriptscriptfont\z@=\tenrm
  \textfont\@ne=\fourteeni  \scriptfont\@ne=\twelvei  \scriptscriptfont\@ne=\teni
  \textfont\tw@=\fourteensy \scriptfont\tw@=\twelvesy \scriptscriptfont\tw@=\tensy
  \textfont\thr@@=\fourteenex \scriptfont\thr@@=\twelveex
        \scriptscriptfont\thr@@=\twelveex
  \textfont\itfam=\fourteenit \scriptfont\itfam=\twelveit
        \scriptscriptfont\itfam=\twelveit
  \textfont\bffam=\fourteenbf \scriptfont\bffam=\twelvebf
        \scriptscriptfont\bffam=\tenbf
  \setbox\strutbox\hbox{\vrule height12.2\p@ depth5\p@ width\z@}%
  \setbox\strutbox@\hbox{\lower.72\normallineskiplimit\vbox{%
        \kern-\normallineskiplimit\copy\strutbox}}%
 \setbox\z@\vbox{\hbox{$($}\kern\z@}\bigsize@=1.7\ht\z@
 \fi
 \normalbaselines\rm\ex@.2326ex\jot4.3\ex@\the\fourteenpoint@}

\font@\seventeenrm=cmr10 scaled\magstep3
\font@\seventeenit=cmti10 scaled\magstep3
\font@\seventeensl=cmsl10 scaled\magstep3
\font@\seventeensmc=cmcsc10 scaled\magstep3
\font@\seventeentt=cmtt10 scaled\magstep3
\font@\seventeenbf=cmbx10 scaled\magstep3
\font@\seventeeni=cmmi10 scaled\magstep3
\font@\seventeensy=cmsy10 scaled\magstep3
\font@\seventeenex=cmex10 scaled\magstep3
\font@\seventeenmsa=msam10 scaled\magstep3
\font@\seventeeneufm=eufm10 scaled\magstep3
\font@\seventeenmsb=msbm10 scaled\magstep3
\newtoks\seventeenpoint@
\def\seventeenpoint{\normalbaselineskip18\p@
 \abovedisplayskip21.6\p@ plus5.2\p@ minus15.4\p@
 \belowdisplayskip\abovedisplayskip
 \abovedisplayshortskip\z@ plus5.2\p@
 \belowdisplayshortskip12.1\p@ plus5.2\p@ minus7\p@
 \textonlyfont@\rm\seventeenrm \textonlyfont@\it\seventeenit
 \textonlyfont@\sl\seventeensl \textonlyfont@\bf\seventeenbf
 \textonlyfont@\smc\seventeensmc \textonlyfont@\tt\seventeentt
%
 \ifsyntax@ \def\big##1{{\hbox{$\left##1\right.$}}}%
  \let\Big\big \let\bigg\big \let\Bigg\big
 \else
  \textfont\z@=\seventeenrm  \scriptfont\z@=\fourteenrm  \scriptscriptfont\z@=\twelverm
  \textfont\@ne=\seventeeni  \scriptfont\@ne=\fourteeni  \scriptscriptfont\@ne=\twelvei
  \textfont\tw@=\seventeensy \scriptfont\tw@=\fourteensy \scriptscriptfont\tw@=\twelvesy
  \textfont\thr@@=\seventeenex \scriptfont\thr@@=\fourteenex
        \scriptscriptfont\thr@@=\fourteenex
  \textfont\itfam=\seventeenit \scriptfont\itfam=\fourteenit
        \scriptscriptfont\itfam=\fourteenit
  \textfont\bffam=\seventeenbf \scriptfont\bffam=\fourteenbf
        \scriptscriptfont\bffam=\twelvebf
  \setbox\strutbox\hbox{\vrule height14.6\p@ depth6\p@ width\z@}%
  \setbox\strutbox@\hbox{\lower.86\normallineskiplimit\vbox{%
        \kern-\normallineskiplimit\copy\strutbox}}%
 \setbox\z@\vbox{\hbox{$($}\kern\z@}\bigsize@=2\ht\z@
 \fi
 \normalbaselines\rm\ex@.2326ex\jot5.2\ex@\the\seventeenpoint@}

\catcode`\@=13

\def\LL{\leavevmode\setbox0=\hbox{L}\hbox to\wd0{\hss\char'40L}}
\def\al{\alpha}

\def\de{\delta}

\def\la{\lambda}

\def\Ga{\Gamma}
\def\De{\Delta}


\def\today{\ifcase\month\or
 January\or February\or March\or April\or May\or June\or
 July\or August\or September\or October\or November\or December\fi
 \space\number\day, \number\year}
\def\dfrac#1#2{{\displaystyle{#1\over#2}}}

\def\({\left(}
\def\){\right)}
\def\[{\left[}
\def\]{\right]}

\def\3{\ss}
\catcode`\@=11
\def\dddot#1{\vbox{\ialign{##\crcr
      .\hskip-.5pt.\hskip-.5pt.\crcr\noalign{\kern1.5\p@\nointerlineskip}
      $\hfil\displaystyle{#1}\hfil$\crcr}}}

\newif\iftab@\tab@false
\newif\ifvtab@\vtab@false
\def\tab{\bgroup\tab@true\vtab@false\vst@bfalse\Strich@false%
   \def\\{\global\hline@@false%
     \ifhline@\global\hline@false\global\hline@@true\fi\cr}
   \edef\l@{\the\leftskip}\ialign\bgroup\hskip\l@##\hfil&&##\hfil\cr}
\def\endtab{\cr\egroup\egroup}
\def\vtab{\vtop\bgroup\vst@bfalse\vtab@true\tab@true\Strich@false%
   \bgroup\def\\{\cr}\ialign\bgroup&##\hfil\cr}
\def\endvtab{\cr\egroup\egroup\egroup}
\def\stab{\D@cke0.5pt\null 
 \bgroup\tab@true\vtab@false\vst@bfalse\Strich@true\Let@@\vspace@
 \normalbaselines\offinterlineskip
  \openup\spreadmlines@
 \edef\l@{\the\leftskip}\ialign
 \bgroup\hskip\l@##\hfil&&##\hfil\crcr}
\def\endstab{\crcr\egroup
 \egroup}
\newif\ifvst@b\vst@bfalse
\def\vstab{\D@cke0.5pt\null
 \vtop\bgroup\tab@true\vtab@false\vst@btrue\Strich@true\bgroup\Let@@\vspace@
 \normalbaselines\offinterlineskip
  \openup\spreadmlines@\bgroup}
\def\endvstab{\crcr\egroup\egroup
 \egroup\tab@false\Strich@false}

\newdimen\htstrut@
\htstrut@8.5\p@
\newdimen\htStrut@
\htStrut@12\p@
\newdimen\dpstrut@
\dpstrut@3.5\p@
\newdimen\dpStrut@
\dpStrut@3.5\p@
\def\openup{\afterassignment\@penup\dimen@=}
\def\@penup{\advance\lineskip\dimen@
  \advance\baselineskip\dimen@
  \advance\lineskiplimit\dimen@
  \divide\dimen@ by2
  \advance\htstrut@\dimen@
  \advance\htStrut@\dimen@
  \advance\dpstrut@\dimen@
  \advance\dpStrut@\dimen@}
\def\Let@@{\relax%
    \def\\{\global\hline@@false%
     \ifhline@\global\hline@false\global\hline@@true\fi\cr}%
    \iffalse}\fi}
\def\matrix{\null\,\vcenter\bgroup
 \tab@false\vtab@false\vst@bfalse\Strich@false\Let@@\vspace@
 \normalbaselines\openup\spreadmlines@\ialign
 \bgroup\hfil$\m@th##$\hfil&&\quad\hfil$\m@th##$\hfil\crcr
 \Mathstrut@\crcr\noalign{\kern-\baselineskip}}
\def\endmatrix{\crcr\Mathstrut@\crcr\noalign{\kern-\baselineskip}\egroup
 \egroup\,}
\def\smatrix{\D@cke0.5pt\null\,
 \vcenter\bgroup\tab@false\vtab@false\vst@bfalse\Strich@true\Let@@\vspace@
 \normalbaselines\offinterlineskip
  \openup\spreadmlines@\ialign
 \bgroup\hfil$\m@th##$\hfil&&\quad\hfil$\m@th##$\hfil\crcr}
\def\endsmatrix{\crcr\egroup
 \egroup\,\Strich@false}
\newdimen\D@cke
\def\Dicke#1{\global\D@cke#1}
\newtoks\tabs@\tabs@{&}
\newif\ifStrich@\Strich@false
\newif\iff@rst

\def\Stricherr@{\iftab@\ifvtab@\errmessage{\noexpand\s not allowed
     here. Use \noexpand\vstab!}%
  \else\errmessage{\noexpand\s not allowed here. Use \noexpand\stab!}%
  \fi\else\errmessage{\noexpand\s not allowed
     here. Use \noexpand\smatrix!}\fi}
\def\format{\ifvst@b\else\crcr\fi\egroup\iffalse{\fi\ifnum`}=0 \fi\format@}
\def\format@#1\\{\def\preamble@{#1}%
 \def\Str@chfehlt##1{\ifx##1\s\Stricherr@\fi\ifx##1\\\let\Next\relax%
   \else\let\Next\Str@chfehlt\fi\Next}%
 \def\c{\hfil\noexpand\ifhline@@\hbox{\vrule height\htStrut@%
   depth\dpstrut@ width\z@}\noexpand\fi%
   \ifStrich@\hbox{\vrule height\htstrut@ depth\dpstrut@ width\z@}%
   \fi\iftab@\else$\m@th\fi\the\hashtoks@\iftab@\else$\fi\hfil}%
 \def\r{\hfil\noexpand\ifhline@@\hbox{\vrule height\htStrut@%
   depth\dpstrut@ width\z@}\noexpand\fi%
   \ifStrich@\hbox{\vrule height\htstrut@ depth\dpstrut@ width\z@}%
   \fi\iftab@\else$\m@th\fi\the\hashtoks@\iftab@\else$\fi}%
 \def\l{\noexpand\ifhline@@\hbox{\vrule height\htStrut@%
   depth\dpstrut@ width\z@}\noexpand\fi%
   \ifStrich@\hbox{\vrule height\htstrut@ depth\dpstrut@ width\z@}%
   \fi\iftab@\else$\m@th\fi\the\hashtoks@\iftab@\else$\fi\hfil}%
 \def\s{\ifStrich@\ \the\tabs@\vrule width\D@cke\the\hashtoks@%
          \fi\the\tabs@\ }%
 \def\sa{\ifStrich@\vrule width\D@cke\the\hashtoks@%
            \the\tabs@\ %
            \fi}%
 \def\se{\ifStrich@\ \the\tabs@\vrule width\D@cke\the\hashtoks@\fi}%
 \def\cd{\hfil\noexpand\ifhline@@\hbox{\vrule height\htStrut@%
   depth\dpstrut@ width\z@}\noexpand\fi%
   \ifStrich@\hbox{\vrule height\htstrut@ depth\dpstrut@ width\z@}%
   \fi$\dsize\m@th\the\hashtoks@$\hfil}%
 \def\rd{\hfil\noexpand\ifhline@@\hbox{\vrule height\htStrut@%
   depth\dpstrut@ width\z@}\noexpand\fi%
   \ifStrich@\hbox{\vrule height\htstrut@ depth\dpstrut@ width\z@}%
   \fi$\dsize\m@th\the\hashtoks@$}%
 \def\ld{\noexpand\ifhline@@\hbox{\vrule height\htStrut@%
   depth\dpstrut@ width\z@}\noexpand\fi%
   \ifStrich@\hbox{\vrule height\htstrut@ depth\dpstrut@ width\z@}%
   \fi$\dsize\m@th\the\hashtoks@$\hfil}%
 \ifStrich@\else\Str@chfehlt#1\\\fi%
 \setbox\z@\hbox{\xdef\Preamble@{\preamble@}}\ifnum`{=0 \fi\iffalse}\fi
 \ialign\bgroup\span\Preamble@\crcr}
\newif\ifhline@\hline@false
\newif\ifhline@@\hline@@false
\def\hlinefor#1{\multispan@{\strip@#1 }\leaders\hrule height\D@cke\hfill%
    \global\hline@true\ignorespaces}
\def\Item "#1"{\par\noindent\hangindent2\parindent%
  \hangafter1\setbox0\hbox{\rm#1\enspace}\ifdim\wd0>2\parindent%
  \box0\else\hbox to 2\parindent{\rm#1\hfil}\fi\ignorespaces}
\def\ITEM #1"#2"{\par\noindent\hangafter1\hangindent#1%
  \setbox0\hbox{\rm#2\enspace}\ifdim\wd0>#1%
  \box0\else\hbox to 0pt{\rm#2\hss}\hskip#1\fi\ignorespaces}
\def\item"#1"{\par\noindent\hang%
  \setbox0=\hbox{\rm#1\enspace}\ifdim\wd0>\the\parindent%
  \box0\else\hbox to \parindent{\rm#1\hfil}\enspace\fi\ignorespaces}
\let\plainitem@\item
\catcode`\@=13

\catcode`\@=11
\font\tenln    = line10
\font\tenlnw   = linew10

\newskip\Einheit \Einheit=0.5cm
\newcount\xcoord \newcount\ycoord
\newdimen\xdim \newdimen\ydim \newdimen\PfadD@cke \newdimen\Pfadd@cke

\newcount\@tempcnta
\newcount\@tempcntb

\newdimen\@tempdima
\newdimen\@tempdimb

\newdimen\@wholewidth
\newdimen\@halfwidth

\newcount\@xarg
\newcount\@yarg
\newcount\@yyarg
\newbox\@linechar
\newbox\@tempboxa
\newdimen\@linelen
\newdimen\@clnwd
\newdimen\@clnht

\newif\if@negarg

\def\@whilenoop#1{}
\def\@whiledim#1\do #2{\ifdim #1\relax#2\@iwhiledim{#1\relax#2}\fi}
\def\@iwhiledim#1{\ifdim #1\let\@nextwhile=\@iwhiledim
        \else\let\@nextwhile=\@whilenoop\fi\@nextwhile{#1}}

\def\@whileswnoop#1\fi{}
\def\@whilesw#1\fi#2{#1#2\@iwhilesw{#1#2}\fi\fi}
\def\@iwhilesw#1\fi{#1\let\@nextwhile=\@iwhilesw
         \else\let\@nextwhile=\@whileswnoop\fi\@nextwhile{#1}\fi}

\def\thinlines{\let\@linefnt\tenln \let\@circlefnt\tencirc
  \@wholewidth\fontdimen8\tenln \@halfwidth .5\@wholewidth}
\def\thicklines{\let\@linefnt\tenlnw \let\@circlefnt\tencircw
  \@wholewidth\fontdimen8\tenlnw \@halfwidth .5\@wholewidth}
\thinlines

\PfadD@cke1pt \Pfadd@cke0.5pt
\def\PfadDicke#1{\PfadD@cke#1 \divide\PfadD@cke by2 \Pfadd@cke\PfadD@cke \multiply\PfadD@cke by2}
\long\def\LOOP#1\REPEAT{\def\BODY{#1}\ITERATE}
\def\ITERATE{\BODY \let\next\ITERATE \else\let\next\relax\fi \next}
\let\REPEAT=\fi
\def\Punkt{\hbox{\raise-2pt\hbox to0pt{\hss$\ssize\bullet$\hss}}}
\def\DuennPunkt(#1,#2){\unskip
  \raise#2 \Einheit\hbox to0pt{\hskip#1 \Einheit
          \raise-2.5pt\hbox to0pt{\hss$\bullet$\hss}\hss}}
\def\NormalPunkt(#1,#2){\unskip
  \raise#2 \Einheit\hbox to0pt{\hskip#1 \Einheit
          \raise-3pt\hbox to0pt{\hss\twelvepoint$\bullet$\hss}\hss}}
\def\DickPunkt(#1,#2){\unskip
  \raise#2 \Einheit\hbox to0pt{\hskip#1 \Einheit
          \raise-4pt\hbox to0pt{\hss\fourteenpoint$\bullet$\hss}\hss}}
\def\Kreis(#1,#2){\unskip
  \raise#2 \Einheit\hbox to0pt{\hskip#1 \Einheit
          \raise-4pt\hbox to0pt{\hss\fourteenpoint$\circ$\hss}\hss}}

\def\Line@(#1,#2)#3{\@xarg #1\relax \@yarg #2\relax
\@linelen=#3\Einheit
\ifnum\@xarg =0 \@vline
  \else \ifnum\@yarg =0 \@hline \else \@sline\fi
\fi}

\def\@sline{\ifnum\@xarg< 0 \@negargtrue \@xarg -\@xarg \@yyarg -\@yarg
  \else \@negargfalse \@yyarg \@yarg \fi
\ifnum \@yyarg >0 \@tempcnta\@yyarg \else \@tempcnta -\@yyarg \fi
\ifnum\@tempcnta>6 \@badlinearg\@tempcnta0 \fi
\ifnum\@xarg>6 \@badlinearg\@xarg 1 \fi
\setbox\@linechar\hbox{\@linefnt\@getlinechar(\@xarg,\@yyarg)}%
\ifnum \@yarg >0 \let\@upordown\raise \@clnht\z@
   \else\let\@upordown\lower \@clnht \ht\@linechar\fi
\@clnwd=\wd\@linechar
\if@negarg \hskip -\wd\@linechar \def\@tempa{\hskip -2\wd\@linechar}\else
     \let\@tempa\relax \fi
\@whiledim \@clnwd <\@linelen \do
  {\@upordown\@clnht\copy\@linechar
   \@tempa
   \advance\@clnht \ht\@linechar
   \advance\@clnwd \wd\@linechar}%
\advance\@clnht -\ht\@linechar
\advance\@clnwd -\wd\@linechar
\@tempdima\@linelen\advance\@tempdima -\@clnwd
\@tempdimb\@tempdima\advance\@tempdimb -\wd\@linechar
\if@negarg \hskip -\@tempdimb \else \hskip \@tempdimb \fi
\multiply\@tempdima \@m
\@tempcnta \@tempdima \@tempdima \wd\@linechar \divide\@tempcnta \@tempdima
\@tempdima \ht\@linechar \multiply\@tempdima \@tempcnta
\divide\@tempdima \@m
\advance\@clnht \@tempdima
\ifdim \@linelen <\wd\@linechar
   \hskip \wd\@linechar
  \else\@upordown\@clnht\copy\@linechar\fi}

\def\@hline{\ifnum \@xarg <0 \hskip -\@linelen \fi
\vrule height\Pfadd@cke width \@linelen depth\Pfadd@cke
\ifnum \@xarg <0 \hskip -\@linelen \fi}

\def\@getlinechar(#1,#2){\@tempcnta#1\relax\multiply\@tempcnta 8
\advance\@tempcnta -9 \ifnum #2>0 \advance\@tempcnta #2\relax\else
\advance\@tempcnta -#2\relax\advance\@tempcnta 64 \fi
\char\@tempcnta}

\def\Vektor(#1,#2)#3(#4,#5){\unskip\leavevmode
  \xcoord#4\relax \ycoord#5\relax
      \raise\ycoord \Einheit\hbox to0pt{\hskip\xcoord \Einheit
         \Vector@(#1,#2){#3}\hss}}

\def\Vector@(#1,#2)#3{\@xarg #1\relax \@yarg #2\relax
\@tempcnta \ifnum\@xarg<0 -\@xarg\else\@xarg\fi
\ifnum\@tempcnta<5\relax
\@linelen=#3\Einheit
\ifnum\@xarg =0 \@vvector
  \else \ifnum\@yarg =0 \@hvector \else \@svector\fi
\fi
\else\@badlinearg\fi}

\def\@hvector{\@hline\hbox to 0pt{\@linefnt
\ifnum \@xarg <0 \@getlarrow(1,0)\hss\else
    \hss\@getrarrow(1,0)\fi}}

\def\@vvector{\ifnum \@yarg <0 \@downvector \else \@upvector \fi}

\def\@svector{\@sline
\@tempcnta\@yarg \ifnum\@tempcnta <0 \@tempcnta=-\@tempcnta\fi
\ifnum\@tempcnta <5
  \hskip -\wd\@linechar
  \@upordown\@clnht \hbox{\@linefnt  \if@negarg
  \@getlarrow(\@xarg,\@yyarg) \else \@getrarrow(\@xarg,\@yyarg) \fi}%
\else\@badlinearg\fi}

\def\@upline{\hbox to \z@{\hskip -.5\Pfadd@cke \vrule width \Pfadd@cke
   height \@linelen depth \z@\hss}}

\def\@downline{\hbox to \z@{\hskip -.5\Pfadd@cke \vrule width \Pfadd@cke
   height \z@ depth \@linelen \hss}}

\def\@upvector{\@upline\setbox\@tempboxa\hbox{\@linefnt\char'66}\raise
     \@linelen \hbox to\z@{\lower \ht\@tempboxa\box\@tempboxa\hss}}

\def\@downvector{\@downline\lower \@linelen
      \hbox to \z@{\@linefnt\char'77\hss}}

\def\@getlarrow(#1,#2){\ifnum #2 =\z@ \@tempcnta='33\else
\@tempcnta=#1\relax\multiply\@tempcnta \sixt@@n \advance\@tempcnta
-9 \@tempcntb=#2\relax\multiply\@tempcntb \tw@
\ifnum \@tempcntb >0 \advance\@tempcnta \@tempcntb\relax
\else\advance\@tempcnta -\@tempcntb\advance\@tempcnta 64
\fi\fi\char\@tempcnta}

\def\@getrarrow(#1,#2){\@tempcntb=#2\relax
\ifnum\@tempcntb < 0 \@tempcntb=-\@tempcntb\relax\fi
\ifcase \@tempcntb\relax \@tempcnta='55 \or
\ifnum #1<3 \@tempcnta=#1\relax\multiply\@tempcnta
24 \advance\@tempcnta -6 \else \ifnum #1=3 \@tempcnta=49
\else\@tempcnta=58 \fi\fi\or
\ifnum #1<3 \@tempcnta=#1\relax\multiply\@tempcnta
24 \advance\@tempcnta -3 \else \@tempcnta=51\fi\or
\@tempcnta=#1\relax\multiply\@tempcnta
\sixt@@n \advance\@tempcnta -\tw@ \else
\@tempcnta=#1\relax\multiply\@tempcnta
\sixt@@n \advance\@tempcnta 7 \fi\ifnum #2<0 \advance\@tempcnta 64 \fi
\char\@tempcnta}

\def\Diagonale(#1,#2)#3{\unskip\leavevmode
  \xcoord#1\relax \ycoord#2\relax
      \raise\ycoord \Einheit\hbox to0pt{\hskip\xcoord \Einheit
         \Line@(1,1){#3}\hss}}
\def\AntiDiagonale(#1,#2)#3{\unskip\leavevmode
  \xcoord#1\relax \ycoord#2\relax 
      \raise\ycoord \Einheit\hbox to0pt{\hskip\xcoord \Einheit
         \Line@(1,-1){#3}\hss}}
\def\Pfad(#1,#2),#3\endPfad{\unskip\leavevmode
  \xcoord#1 \ycoord#2 \thicklines\ZeichnePfad#3\endPfad\thinlines}
\def\ZeichnePfad#1{\ifx#1\endPfad\let\next\relax
  \else\let\next\ZeichnePfad
    \ifnum#1=1
      \raise\ycoord \Einheit\hbox to0pt{\hskip\xcoord \Einheit
         \vrule height\Pfadd@cke width1 \Einheit depth\Pfadd@cke\hss}%
      \advance\xcoord by 1
    \else\ifnum#1=2
      \raise\ycoord \Einheit\hbox to0pt{\hskip\xcoord \Einheit
        \hbox{\hskip-\PfadD@cke\vrule height1 \Einheit width\PfadD@cke depth0pt}\hss}%
      \advance\ycoord by 1
    \else\ifnum#1=3
      \raise\ycoord \Einheit\hbox to0pt{\hskip\xcoord \Einheit
         \Line@(1,1){1}\hss}
      \advance\xcoord by 1
      \advance\ycoord by 1
    \else\ifnum#1=4
      \raise\ycoord \Einheit\hbox to0pt{\hskip\xcoord \Einheit
         \Line@(1,-1){1}\hss}
      \advance\xcoord by 1
      \advance\ycoord by -1
    \fi\fi\fi\fi
  \fi\next}
\def\hSSchritt{\leavevmode\raise-.4pt\hbox to0pt{\hss.\hss}\hskip.2\Einheit
  \raise-.4pt\hbox to0pt{\hss.\hss}\hskip.2\Einheit
  \raise-.4pt\hbox to0pt{\hss.\hss}\hskip.2\Einheit
  \raise-.4pt\hbox to0pt{\hss.\hss}\hskip.2\Einheit
  \raise-.4pt\hbox to0pt{\hss.\hss}\hskip.2\Einheit}
\def\vSSchritt{\vbox{\baselineskip.2\Einheit\lineskiplimit0pt
\hbox{.}\hbox{.}\hbox{.}\hbox{.}\hbox{.}}}
\def\DSSchritt{\leavevmode\raise-.4pt\hbox to0pt{%
  \hbox to0pt{\hss.\hss}\hskip.2\Einheit
  \raise.2\Einheit\hbox to0pt{\hss.\hss}\hskip.2\Einheit
  \raise.4\Einheit\hbox to0pt{\hss.\hss}\hskip.2\Einheit
  \raise.6\Einheit\hbox to0pt{\hss.\hss}\hskip.2\Einheit
  \raise.8\Einheit\hbox to0pt{\hss.\hss}\hss}}
\def\dSSchritt{\leavevmode\raise-.4pt\hbox to0pt{%
  \hbox to0pt{\hss.\hss}\hskip.2\Einheit
  \raise-.2\Einheit\hbox to0pt{\hss.\hss}\hskip.2\Einheit
  \raise-.4\Einheit\hbox to0pt{\hss.\hss}\hskip.2\Einheit
  \raise-.6\Einheit\hbox to0pt{\hss.\hss}\hskip.2\Einheit
  \raise-.8\Einheit\hbox to0pt{\hss.\hss}\hss}}
\def\SPfad(#1,#2),#3\endSPfad{\unskip\leavevmode
  \xcoord#1 \ycoord#2 \ZeichneSPfad#3\endSPfad}
\def\ZeichneSPfad#1{\ifx#1\endSPfad\let\next\relax
  \else\let\next\ZeichneSPfad
    \ifnum#1=1
      \raise\ycoord \Einheit\hbox to0pt{\hskip\xcoord \Einheit
         \hSSchritt\hss}%
      \advance\xcoord by 1
    \else\ifnum#1=2
      \raise\ycoord \Einheit\hbox to0pt{\hskip\xcoord \Einheit
        \hbox{\hskip-2pt \vSSchritt}\hss}%
      \advance\ycoord by 1
    \else\ifnum#1=3
      \raise\ycoord \Einheit\hbox to0pt{\hskip\xcoord \Einheit
         \DSSchritt\hss}
      \advance\xcoord by 1
      \advance\ycoord by 1
    \else\ifnum#1=4
      \raise\ycoord \Einheit\hbox to0pt{\hskip\xcoord \Einheit
         \dSSchritt\hss}
      \advance\xcoord by 1
      \advance\ycoord by -1
    \fi\fi\fi\fi
  \fi\next}
\def\Koordinatenachsen(#1,#2){\unskip
 \hbox to0pt{\hskip-.5pt\vrule height#2 \Einheit width.5pt depth1 \Einheit}%
 \hbox to0pt{\hskip-1 \Einheit \xcoord#1 \advance\xcoord by1
    \vrule height0.25pt width\xcoord \Einheit depth0.25pt\hss}}
\def\Koordinatenachsen(#1,#2)(#3,#4){\unskip
 \hbox to0pt{\hskip-.5pt \ycoord-#4 \advance\ycoord by1
    \vrule height#2 \Einheit width.5pt depth\ycoord \Einheit}%
 \hbox to0pt{\hskip-1 \Einheit \hskip#3\Einheit 
    \xcoord#1 \advance\xcoord by1 \advance\xcoord by-#3 
    \vrule height0.25pt width\xcoord \Einheit depth0.25pt\hss}}
\def\Gitter(#1,#2){\unskip \xcoord0 \ycoord0 \leavevmode
  \LOOP\ifnum\ycoord<#2
    \loop\ifnum\xcoord<#1
      \raise\ycoord \Einheit\hbox to0pt{\hskip\xcoord \Einheit\Punkt\hss}%
      \advance\xcoord by1
    \repeat
    \xcoord0
    \advance\ycoord by1
  \REPEAT}
\def\Gitter(#1,#2)(#3,#4){\unskip \xcoord#3 \ycoord#4 \leavevmode
  \LOOP\ifnum\ycoord<#2
    \loop\ifnum\xcoord<#1
      \raise\ycoord \Einheit\hbox to0pt{\hskip\xcoord \Einheit\Punkt\hss}%
      \advance\xcoord by1
    \repeat
    \xcoord#3
    \advance\ycoord by1
  \REPEAT}
\def\Label#1#2(#3,#4){\unskip \xdim#3 \Einheit \ydim#4 \Einheit
  \def\lo{\advance\xdim by-.5 \Einheit \advance\ydim by.5 \Einheit}%
  \def\llo{\advance\xdim by-.25cm \advance\ydim by.5 \Einheit}%
  \def\loo{\advance\xdim by-.5 \Einheit \advance\ydim by.25cm}%
  \def\o{\advance\ydim by.25cm}%
  \def\ro{\advance\xdim by.5 \Einheit \advance\ydim by.5 \Einheit}%
  \def\rro{\advance\xdim by.25cm \advance\ydim by.5 \Einheit}%
  \def\roo{\advance\xdim by.5 \Einheit \advance\ydim by.25cm}%
  \def\l{\advance\xdim by-.30cm}%
  \def\r{\advance\xdim by.30cm}%
  \def\lu{\advance\xdim by-.5 \Einheit \advance\ydim by-.6 \Einheit}%
  \def\llu{\advance\xdim by-.25cm \advance\ydim by-.6 \Einheit}%
  \def\luu{\advance\xdim by-.5 \Einheit \advance\ydim by-.30cm}%
  \def\u{\advance\ydim by-.30cm}%
  \def\ru{\advance\xdim by.5 \Einheit \advance\ydim by-.6 \Einheit}%
  \def\rru{\advance\xdim by.25cm \advance\ydim by-.6 \Einheit}%
  \def\ruu{\advance\xdim by.5 \Einheit \advance\ydim by-.30cm}%
  #1\raise\ydim\hbox to0pt{\hskip\xdim
     \vbox to0pt{\vss\hbox to0pt{\hss$#2$\hss}\vss}\hss}%
}
\catcode`\@=13

\TagsOnRight

\topmatter 
\title Proof of a determinant evaluation conjectured by Bombieri,
Hunt and van der Poorten
\endtitle 
\author C.~Krattenthaler\footnote{\hbox{Research supported in part by the
MSRI, Berkeley}} and D.~Zeilberger
\endauthor 
\affil 
Institut f\"ur Mathematik der Universit\"at Wien,\\
Strudlhofgasse 4, A-1090 Wien, Austria.\\
e-mail: KRATT\@Pap.Univie.Ac.At\\
WWW: {\tt http://radon.mat.univie.ac.at/People/kratt}\\
\vskip6pt
Department of Mathematics, Temple University\\
Philadelphia, PA19122, USA\\
e-mail: zeilberg\@math.temple.edu\\
WWW: {\tt http://www.math.temple.edu/\~{}zeilb}
\endaffil
\address Institut f\"ur Mathematik der Universit\"at Wien,
Strudlhofgasse 4, A-1090 Wien, Austria.
\endaddress
\address
Department of Mathematics, Temple University,
Philadelphia, PA19122, USA
\endaddress
\abstract 
A determinant evaluation is proven, a special case of which establishes a
conjecture of Bombieri, Hunt, and van der Poorten (Experimental
Math\. {\bf 4} (1995), 87--96) that arose in the
study of Thue's method of approximating algebraic numbers.
\endabstract
\endtopmatter
\document

\leftheadtext{C. Krattenthaler and D. Zeilberger}
\rightheadtext{Proof of a determinant evaluation}

\subhead 1. Introduction\endsubhead
In their study \cite{\BoHPAA} of Thue's method of approximating an
algebraic number, Bombieri, Hunt, and van der Poorten conjectured two
determinant evaluations, one of which can be restated as follows.
\proclaim{Conjecture (Bombieri, Hunt, van der Poorten \cite{\BoHPAA,
next-to-last paragraph})}\linebreak
Let $b,c$ be nonnegative integers, $c\le b$, and
let $\De(b,c)$ be the determinant of the $(b+c)\times(b+c)$ matrix
(given in block form)
\vskip10pt
$$
\(
\PfadDicke{.3pt}
\SPfad(0,1),111111111111\endSPfad
\SPfad(0,-1),111111111111\endSPfad
\SPfad(4,-3),222222\endSPfad
\SPfad(8,-3),222222\endSPfad
\Label\o{\raise25pt\hbox{$\dsize\binom ji$}}(2,1)
\Label\o{\raise25pt\hbox{$\dsize\binom ji$}}(6,1)
\Label\o{\raise25pt\hbox{$\dsize\binom ji$}}(10,1)
\Label\o{\raise25pt\hbox{$\dsize\binom ji$}}(6,-1)
\Label\o{\raise25pt\hbox{$\dsize\binom ji$}}(10,-1)
\Label\o{\raise25pt\hbox{$\dsize 2\binom j{i-b}$}}(2,-3)
\Label\o{\raise25pt\hbox{$\dsize\binom j{i-b}$}}(6,-3)
\Label\o{\raise15pt\hbox{\fourteenpoint$0$}}(2,-1)
\Label\o{\raise15pt\hbox{\fourteenpoint$0$}}(10,-3)
\Label\o{\raise20pt\hbox{$0\le
i<c\hphantom{b+{}}$\hskip10pt}}(-3,1)
\Label\o{\raise20pt\hbox{$c\le
i<b\hphantom{c+{}}$\hskip10pt}}(-3,-1)
\Label\o{\raise20pt\hbox{$b\le i<b+c$\hskip10pt}}(-3,-3)
\hskip6cm
\)
\hbox{\hskip-6.3cm}
\Label\o{\raise15pt\hbox{$0\le j<c$}}(2,3)
\Label\o{\raise15pt\hbox{$c\le j<b$}}(6,3)
\Label\o{\raise15pt\hbox{$b\le j<b+c$}}(10,3)
\hskip6.5cm.\hskip-2.5cm
\tag1.1
$$
Then
\roster
\item"{\rm (i)}" $\De(b,c)=0$ if $b$ is even and $c$ is odd;
\item"{\rm (ii)}" if any of these conditions does not hold, and if $b\ge 2c$, then
$$\De(b,c)=\pm \frac {[b-c/2]^2\,[b-2c]\, [(b+c)/2]^2\, [(b-c)/2]^6\,
[c/2]^6} {[b-c]^3\, [b/2]^6\, [b/2-c]^2\, [c]^3},\tag1.2$$
where $[s]:=\prod _{k=0} ^{s-1}k!$ if $s$ is an integer, and
$[s]^2=\(\prod
_{k=0} ^{s-1/2}k!\)\(\prod _{k=0} ^{s-3/2}k!\)$ if $s$ is a
half-integer;
\item"{\rm (iii)}" whereas if $b<2c$ then
$$\De(b,c)=\pm 2^{2c-b}\De(b,b-c).\tag1.3$$
\endroster
\endproclaim

The purpose of this paper is to prove a determinant evaluation,
containing the parameter $x$, which for $x=0$ 
reduces to the above Conjecture.

\proclaim{Theorem 1}
Let $b,c$ be nonnegative integers, $c\le b$, and
let $\De(x;b,c)$ be the determinant of the $(b+c)\times(b+c)$ matrix
\vskip10pt
$$\(
\PfadDicke{.3pt}
\SPfad(0,1),111111111111\endSPfad
\SPfad(0,-1),111111111111\endSPfad
\SPfad(4,-3),222222\endSPfad
\SPfad(8,-3),222222\endSPfad
\Label\o{\raise25pt\hbox{$\dsize\binom {x+j}i$}}(2,1)
\Label\o{\raise25pt\hbox{$\dsize\binom {x+j}i$}}(6,1)
\Label\o{\raise25pt\hbox{$\dsize\binom {2x+j}i$}}(10,1)
\Label\o{\raise25pt\hbox{$\dsize\binom {x+j}i$}}(6,-1)
\Label\o{\raise25pt\hbox{$\dsize\binom {2x+j}i$}}(10,-1)
\Label\o{\raise25pt\hbox{$\dsize 2\binom {x+j}{i-b}$}}(2,-3)
\Label\o{\raise25pt\hbox{$\dsize\binom {x+j}{i-b}$}}(6,-3)
\Label\o{\raise15pt\hbox{\fourteenpoint$0$}}(2,-1)
\Label\o{\raise15pt\hbox{\fourteenpoint$0$}}(10,-3)
\Label\o{\raise20pt\hbox{$0\le
i<c\hphantom{b+{}}$\hskip10pt}}(-3,1)
\Label\o{\raise20pt\hbox{$c\le
i<b\hphantom{c+{}}$\hskip10pt}}(-3,-1)
\Label\o{\raise20pt\hbox{$b\le i<b+c$\hskip10pt}}(-3,-3)
\hskip6cm
\)
\hbox{\hskip-6.3cm}
\Label\o{\raise15pt\hbox{$0\le j<c$}}(2,3)
\Label\o{\raise15pt\hbox{$c\le j<b$}}(6,3)
\Label\o{\raise15pt\hbox{$b\le j<b+c$}}(10,3)
\hskip6.5cm.\hskip-2.5cm
\tag1.4
$$
Then
\roster
\item"{\rm (i)}" $\De(x;b,c)=0$ if $b$ is even and $c$ is odd;
\item"{\rm (ii)}" if any of these conditions does not hold, then
$$\multline
\De(x;b,c)=(-1)^{c}2^{c}\prod _{i=1} ^{b-c}\frac {\(i+\frac {1}
{2}-\cl{\frac {b} {2}}\)_c} {(i)_c}\\
\times\prod _{i=1} ^{c}\frac {\(x+\cl{\frac {c+i}
{2}}\)_{b-c+\cl{i/2}-\cl{(c+i)/2}} \, \(x+\cl{\frac {b-c+i}
{2}}\)_{\cl{(b+i)/2}-\cl{(b-c+i)/2}}} {\(\frac {1} {2}-\cl{\frac {b} {2}}
+\cl{\frac {c+i}
{2}}\)_{b-c+\cl{i/2}-\cl{(c+i)/2}} \, \(\frac {1} {2}-\cl{\frac {b} {2}}
+\cl{\frac {b-c+i} {2}}\)_{\cl{(b+i)/2}-\cl{(b-c+i)/2}}},\\
\endmultline\tag1.5$$
\endroster
where the shifted factorial
$(a)_k$ is defined by 
$$(a)_k:=\cases a(a+1)\cdots(a+k-1)&k\ge0\\
\dfrac {1} {(a-1)(a-2)\cdots (a+k)}&k<0.\endcases$$
\endproclaim
(A uniform way to define the shifted factorial is by
$(a)_k:=\Ga(a+k)/\Ga(a)$, respectively by an appropriate limit in
case $a$ or $a+k$ is a nonpositive integer, see \cite{\GrKPAA,
p211f??}.)

The Conjecture does indeed immediately follow, since a routine
calculation shows in particular that the expression (1.5) satisfies
the equation
$$\De(x;b,c)=(-1)^{b}2^{2c-b}\De(x;b,b-c),\tag1.6$$
which implies (1.3) on setting $x=0$.

We are going to prove this Theorem in the next section. For the sake
of clarity of exposition, we defer the proof of some auxiliary facts 
to Section~3.
The method of proof that we use is also applied 
successfully in \cite{\KratBD, \KratBG, \KratBH, \KratBI}
(see in particular the tutorial description in \cite{\KratBI, Sec.~2}).
In order to apply this method, it is actually important to have (at
least) one free
parameter. So, the main difficulty in proving the Conjecture 
was to find the appropriate
generalization of (1.1), such that the determinant still factors
nicely. The various hypergeometric calculations were done, with some
patience, using the first author's {\sl Mathematica\/} package HYP
\cite{\KratBF}. For curiosity, we mention that, although at present it
is quite hopeless to prove any of the identities in this paper by the
recent algorithmic tools \cite{\WiZeAC, \ZeilAM, \ZeilAN, \ZeilAV},
these did, implicitly, have their place in this work. For example,
the fact that the three seemingly very different sums in (3.6) can be
combined into one single sum was discovered by applying the
Gosper-Zeilberger algorithm \cite{\ZeilAM, \ZeilAV} to each of the
three sums in (3.6), being puzzled that one obtains always exactly the same
recurrence, until eventually realizing that, maybe, these sums are in
fact just parts of one and the same series.

\subhead 2. Proof of Theorem~1\endsubhead
We proceed by first reducing the determinant $\De(x;b,c)$, 
which by definition is
the determinant of the matrix (1.4), by elementary row operations to
a constant times a smaller determinant, $\De'(x;b,c)$, given in (2.3).
This smaller determinant is then evaluated in Theorem~2, using a
method which is described and illustrated in \cite{\KratBI, Sec.~2}.

Now we describe the row reductions. We
subtract $1/2$ times the $(i+b)$-th row from the $i$-th row,
$i=0,1,\dots,c$. The resulting matrix has block form, with all the
entries in the $b\times c$ upper-left block being equal 0.
Therefore, up to sign, the determinant decomposes into the product of
the determinant of the
$c\times c$ lower-left block times the determinant of the
$b\times b$ upper-right block:
$$\multline
\De(x;b,c)=(-1)^{bc}\det_{0\le i,j< c}\(2\binom{x+j}i\)\\
\times
\det_{0\le i<b,\, c\le j\le b+c}\!\!\(
\PfadDicke{.3pt}
\SPfad(0,0),11111111\endSPfad
\SPfad(4,-2),2222\endSPfad
\Label\o{\raise25pt\hbox{$\dsize\frac {1} {2}\binom {x+j}i$}}(2,0)
\Label\o{\raise25pt\hbox{$\dsize\binom {2x+j}i$}}(6,0)
\Label\o{\raise25pt\hbox{$\dsize\binom {x+j}i$}}(2,-2)
\Label\o{\raise25pt\hbox{$\dsize\binom {2x+j}i$}}(6,-2)
\Label\o{\raise20pt\hbox{$0\le i<c$\hskip10pt}}(12,0)
\Label\o{\raise20pt\hbox{$c\le i<b$\hskip10pt}}(12,-2)
\hskip4cm
\)
\hbox{\hskip-4.3cm}
\Label\o{\raise0pt\hbox{$c\le j<b$}}(2,-3)
\Label\o{\raise0pt\hbox{$b\le j<b+c$}}(6,-3)
\hskip4.2cm.\hskip2.5cm
\endmultline\tag2.1$$
\vskip10pt
The first determinant is easily evaluated (see \cite{\GeViAA,
Theorem~1 with $a_j=x+j-1$, $b_i=i-1$; there is just one family of
nonintersecting lattice paths in that case!} for an unusual proof),
$$\det_{0\le i,j< c}\(2\binom{x+j}i\)=2^c.\tag2.2$$
So, what we have to do is to evaluate the second determinant, or
equivalently,
$$
\det_{0\le i<b,\, c\le j\le b+c}\!\!\(
\PfadDicke{.3pt}
\SPfad(0,0),11111111\endSPfad
\SPfad(4,-2),2222\endSPfad
\Label\o{\raise25pt\hbox{$\dsize\binom {x+j}i$}}(2,0)
\Label\o{\raise25pt\hbox{$\dsize2\binom {2x+j}i$}}(6,0)
\Label\o{\raise25pt\hbox{$\dsize\binom {x+j}i$}}(2,-2)
\Label\o{\raise25pt\hbox{$\dsize\binom {2x+j}i$}}(6,-2)
\Label\o{\raise20pt\hbox{$0\le i<c$\hskip10pt}}(12,0)
\Label\o{\raise20pt\hbox{$c\le i<b$\hskip10pt}}(12,-2)
\hskip4cm
\)
\hbox{\hskip-4.3cm}
\Label\o{\raise0pt\hbox{$c\le j<b$}}(2,-3)
\Label\o{\raise0pt\hbox{$b\le j<b+c$}}(6,-3)
\hskip4.2cm.\hskip2.5cm
$$
\vskip10pt
This determinant can be further reduced. We subtract column $b-2$
from column $b-1$, column $b-3$ from column $b-2$, \dots, column $c$
from column $c+1$, in that order. Then we subtract column $b-2$
from column $b-1$, column $b-3$ from column $b-2$, \dots, column $c+1$
from column $c+2$ (but not column $c$ from column $c+1$!), 
etc. We do the same sort of operations with columns
$b$, $b+1$, \dots, $b+c-1$. The resulting determinant is
$$
\det_{0\le i<b,\, c\le j\le b+c}\!\!\(
\PfadDicke{.3pt}
\SPfad(0,0),11111111\endSPfad
\SPfad(4,-2),2222\endSPfad
\Label\o{\raise25pt\hbox{$\dsize\binom {x+c}{i-j+c}$}}(2,0)
\Label\o{\raise25pt\hbox{$\dsize2\binom {2x+b}{i-j+b}$}}(6,0)
\Label\o{\raise25pt\hbox{$\dsize\binom {x+c}{i-j+c}$}}(2,-2)
\Label\o{\raise25pt\hbox{$\dsize\binom {2x+b}{i-j+b}$}}(6,-2)
\Label\o{\raise20pt\hbox{$0\le i<c$\hskip10pt}}(12,0)
\Label\o{\raise20pt\hbox{$c\le i<b$\hskip10pt}}(12,-2)
\hskip4cm
\)
\hbox{\hskip-4.3cm}
\Label\o{\raise0pt\hbox{$c\le j<b$}}(2,-3)
\Label\o{\raise0pt\hbox{$b\le j<b+c$}}(6,-3)
\hskip4.2cm.\hskip2.5cm
\tag2.3$$
\vskip10pt
Let us denote this determinant by $\De'(x;b,c)$. Recall that by (2.1)
and (2.2) we have
$$\De(x;b,c)=(-1)^{bc}\De'(x;b,c).\tag2.4$$
The next theorem gives the evaluation of $\De'(x;b,c)$.
\proclaim{Theorem 2}
Let $b,c$ be nonnegative integers, $c\le b$, and let, as before, 
$\De'(x;b,c)$ denote the
determinant in {\rm(2.3)}. Then
\roster
\item"{\rm (i)}" $\De'(x;b,c)=0$ if $b$ is even and $c$ is odd;
\item"{\rm (ii)}" if any of these conditions does not hold, then
$$\multline
\De'(x;b,c)=(-1)^{c(b-c)}2^c\prod _{i=1} ^{b-c}\frac {\(i+\frac {1}
{2}-\cl{\frac {b} {2}}\)_c} {(i)_c}\\
\times\prod _{i=1} ^{c}\frac {\(x+\cl{\frac {c+i}
{2}}\)_{b-c+\cl{i/2}-\cl{(c+i)/2}} \, \(x+\cl{\frac {b-c+i}
{2}}\)_{\cl{(b+i)/2}-\cl{(b-c+i)/2}}} {\(\frac {1} {2}-\cl{\frac {b} {2}}
+\cl{\frac {c+i}
{2}}\)_{b-c+\cl{i/2}-\cl{(c+i)/2}} \, \(\frac {1} {2}-\cl{\frac {b} {2}}
+\cl{\frac {b-c+i} {2}}\)_{\cl{(b+i)/2}-\cl{(b-c+i)/2}}}.\\
\endmultline\tag2.5$$
\endroster
\endproclaim
Clearly, once we have proved Theorem~2, the relation (2.4)
establishes Theorem~1 immediately.

We now proceed with the proof of Theorem~2. It relies on several
Lemmas, which are stated and proved separately as Lemmas~1--4 in
Section~3.
\demo{Proof of Theorem~2} 
We treat both (i) and (ii) at once. That is, for now we just assume
that $b$ and $c$ are nonnegative integers with $c\le b$.

The method that we use to prove the Theorem consists of three steps 
(see \cite{\KratBI, Sec.~2}): 
In the first step we show that the right-hand side of (2.5) divides
$\De'(x;b,c)$ as a polynomial in $x$, regardless what the parity of
$b$ or $c$ is. The reader should observe that, although (2.5) is
going to hold only if $b$ is odd or if both $b$ and $c$ are even, the
right-hand side of (2.5) is nevertheless
well-defined in {\it all\/} cases, as
long as $b\ge c$. Then, in the second step we show
that the degree of $\De'(x;b,c)$, as a polynomial in $x$, is at most
$c(b-c)$. On the other hand, as is easily seen, the degree in $x$ 
of the right-hand side of (2.5) is exactly $c(b-c)$ if $b$ is odd or
if both $b$ and $c$ are even, and is exactly $c(b-c)+1$ if $b$ is
even and $c$ is odd. Therefore, if $b$ is odd or if both $b$ and $c$
are even, the determinant $\De'(x;b,c)$ must
equal the right-hand side of (2.5) times some constant independent of
$x$, and it must be 0 if $b$ is even and $c$ is odd.
The constant in the former case 
is finally determined to be 1 in the third step. This would prove
both (i) and (ii).

\medskip
{\it Step 1. The right-hand side of\/ {\rm (2.5)} divides $\De'(x;b,c)$.}
This is done in Lemmas~1 and 2 in Section~3.

\medskip
{\it Step 2. $\De'(x;b,c)$ is a polynomial in $x$ of degree at most
$c(b-c)$.}
Each term in the defining expansion
of the determinant $\De'(x;b,c)$ (which by definition is the
determinant in (2.3)) has degree $c(b-c)$. Therefore, $\De'(x;b,c)$,
being the sum of all these terms, has degree {\it at most\/}
$c(b-c)$. Therefore,  
since the degree in $x$ of the right-hand side of (2.5) is
exactly $c(b-c)$ if $b$ is odd or if both $b$ and $c$ are even,
$\De'(x;b,c)$ and the right-hand side of (2.5) differ
only by a multiplicative constant, whereas, since 
the degree in $x$ of the right-hand side of (2.5) is
exactly $c(b-c)+1$ if $b$ is even and $c$ is
odd, $\De'(x;b,c)$ can only be 0.

\medskip
{\it Step 3. Determining the multiplicative constant in the case that
$b$ is odd or that both $b$ and $c$ are even.} If we are able
to show that $\De'(x;b,c)$ and the right-hand side of (2.5) do not
vanish and equal each other for some particular value of $x$, then it
is established that the multiplicative constant connecting
$\De'(x;b,c)$ and the right-hand side of (2.5) must be 1. Thus,
equation (2.5) would be proved.

We distinguish between the cases $b$ even or odd.

Let first $b$ be odd. We compare the values of $\De'(x;b,c)$ and the
right-hand side of (2.5) at $x=-b/2$. We have to show that the two values
agree. Now, the right-hand side of (2.5) at $x=-b/2$ equals
$$(-1)^{c(b-c)}2^c\prod _{i=1} ^{b-c}\frac 
{\(i-{\frac {b} {2}}\)_c} {(i)_c}.\tag2.6$$

On the other hand, let us turn to the determinant $\De'(x;b,c)$, given by (2.3),
evaluated at $x=-b/2$. In that case, the upper-right block becomes 2
times the $c\times c$ identity matrix, and the lower-right block
becomes the $(b-c)\times c$ zero matrix. Hence, $\De'(-b/2;b,c)$ equals
$$(-1)^{c(b-c)}2^c\,\det_{c\le i,j<b}\(\binom {c-b/2}{i-j+c}\).$$
The evaluation of this determinant is given by Lemma~3 with $X=c-b/2$. 
Thus we obtain
for $\De'(-b/2;b,c)$ exactly the expression in (2.6).

\smallskip
Now let $b$ be even, and, hence, due to our assumption, also $c$ be
even. In this case, it is of no use to set $x=-b/2$ in (2.5), since
both sides vanish for $x=-b/2$. Instead, we compare $\De'(x;b,c)$ and
the right-hand side of (2.5) at $-b/2+1/2$. Clearly, the right-hand
side at $-b/2+1/2$ equals
$$(-1)^{c(b-c)}2^c\prod _{i=1} ^{b-c}\frac {\(i+\frac {1} {2}
-{\frac {b} {2}}\)_c} {(i)_c}.
\tag2.7$$

Next, we turn to the determinant $\De'(x;b,c)$ evaluated at $x=-b/2+1/2$.
For convenience, we first add
$$\sum _{s=0} ^{i-1}\binom {-1}{i-s} \cdot\big(\text {row $s$ of $\De'(x;b,c)$}\big)$$
to row $i$ of $\De'(x;b,c)$, $i=c-1,c-2,\dots,0$. Thus, making use
of the Chu--Vandermonde summation
(see e.g\. \cite{\GrKPAA, Sec.~5.1, (5.27)}), the determinant
is transformed into
$$
\det_{0\le i<b,\, c\le j\le b+c}\!\!\(
\PfadDicke{.3pt}
\SPfad(0,0),111111111\endSPfad
\SPfad(4,-2),2222\endSPfad
\Label\o{\raise25pt\hbox{$\dsize\binom {x+c-1}{i-j+c}$}}(2,0)
\Label\o{\raise25pt\hbox{\hskip-10pt$\dsize2\binom {2x+b-1}{i-j+b}$}}(7,0)
\Label\o{\raise25pt\hbox{$\dsize\binom {x+c}{i-j+c}$}}(2,-2)
\Label\o{\raise25pt\hbox{\hskip-10pt$\dsize\binom {2x+b}{i-j+b}$}}(7,-2)
\Label\o{\raise20pt\hbox{$0\le i<c$\hskip10pt}}(13,0)
\Label\o{\raise20pt\hbox{$c\le i<b$\hskip10pt}}(13,-2)
\hskip4.5cm
\)
\hbox{\hskip-4.3cm}
\Label\o{\raise0pt\hbox{$c\le j<b$}}(1,-3)
\Label\o{\raise0pt\hbox{\hskip10pt$b\le j<b+c$}}(5,-3)
\hskip4.2cm.\hskip2.5cm
\tag2.8$$
\vskip10pt
In this determinant we set $x=-b/2+1/2$. The effect is that the
upper-right block becomes 2 times the $c\times c$ identity matrix,
while the lower-right block consists of all zeros, except that the
$(c,b+c-1)$-entry equals 1. Accordingly, we expand the determinant
along column $b$, then along column $b+1$, \dots, finally, along
column $b+c-1$. All these columns contain just one entry 2 and 0's
else, with the exception of the last column, which contains two
non-zero entries if $b>c>0$, i.e., if there is a non-empty lower-right
block. By that way, we obtain for our determinant the
difference
$$\multline
(-1)^{c(b-c)}2^c\,\det_{c\le i,j<b}\(\binom {c-b/2+1/2}{i-j+c}\)\\
-\chi(b>c>0)\cdot (-1)^{c(b-c)}2^{c-1}\,
\det_{c\le i,j<b}\pmatrix \dsize\binom {c-b/2-1/2}{i-j+c}&i=c\\
\dsize\binom {c-b/2+1/2}{i-j+c}&i>c
\endpmatrix.
\endmultline\tag2.9$$
Here, $\chi(\Cal A)=1$ if $\Cal A$ is true and $\chi(\Cal A)=0$
otherwise.

The first determinant in (2.9) can be evaluated by means of Lemma~3 with
$X=c-b/2+1/2$, the second determinant is shown to equal 0 in Lemma~4.
Thus we obtain
for $\De'(-b/2+1/2;b,c)$ exactly the expression in (2.7).

\smallskip
This completes the proof of Theorem~2.\quad \quad \qed
\enddemo

\subhead 3. Auxiliary Lemmas\endsubhead
In this section we prove the auxiliary facts that are needed in the
proof of Theorem~2 in the previous section.
\proclaim{Lemma 1}Let $b$ and $c$ be nonnegative
integers such that $b\ge 2c$.
Then the product
$$\prod _{i=1} ^{c}{\(x+\cl{\frac {c+i}
{2}}\)_{b-c+\cl{i/2}-\cl{(c+i)/2}} \, \(x+\cl{\frac {b-c+i}
{2}}\)_{\cl{(b+i)/2}-\cl{(b-c+i)/2}}}\tag3.1$$
divides $\De'(x;b,c)$, the determinant given by {\rm (2.3)}, 
as a polynomial in $x$.
\endproclaim
\demo{Proof} Let us concentrate on some factor $(x+e)$ which appears
in (3.1), say with multiplicity $m(e)$. We have to prove that
$(x+e)^{m(e)}$ divides $\De'(x;b,c)$. We accomplish this by finding
$m(e)$ linear combinations of the rows of $\De'(x;b,c)$ (or of an
equivalent determinant) that vanish
for $x=-e$, and which are linearly independent. See the Lemma in
Section~2 of \cite{\KratBI} for a formal proof of the correctness of
this procedure. 

We have to distinguish between four cases, depending on the
magnitude of $e$. The first case is $c/2\le e\le c$, the second case is
$c\le e\le b/2$, the third case is $b/2\le e\le b-c$, and the fourth
case is $b-c\le e\le b-c/2$. 

\smallskip
{\it Case 1: $c/2\le e\le c$}. By inspection of the expression (3.1),
we see that we have to prove that $(x+e)^{m(e)}$
divides $\De'(x;b,c)$, where
$$m(e)=\cases (2e-c)&c/2\le e\le (b-c)/2\\
(2e-c)+(2e+c-b)&(b-c)/2< e\le c.\endcases
\tag3.2$$
Note that the second case in (3.2) could be empty, but not the first,
because of $b\ge 2c$.

The term $(2e-c)$ in (3.2) is easily explained: 
We take $(x+e)$ out of rows $b+c-2e,b+c-2e+1,\dots,b-1$ of the
determinant $\De'(x;b,c)$. 
It follows from the definition (2.3) of $\De'(x;b,c)$ 
that the remaining determinant is
$$\multline
\hskip-.5cm
\det_{0\le i<b,\, c\le j\le b+c}\!\!\(
\PfadDicke{.3pt}
\SPfad(0,2),1111111111111111\endSPfad
\SPfad(0,0),1111111111111111\endSPfad
\SPfad(8,-3),2222222\endSPfad
\Label\o{\raise25pt\hbox{$\dsize\binom {x+c}{i-j+c}$}}(4,2)
\Label\o{\raise25pt\hbox{$\dsize2\binom {2x+b}{i-j+b}$}}(12,2)
\Label\o{\raise25pt\hbox{$\dsize\binom {x+c}{i-j+c}$}}(4,0)
\Label\o{\raise25pt\hbox{$\dsize\binom {2x+b}{i-j+b}$}}(12,0)
\Label\o{\raise15pt\hbox{\hskip-10pt$\dsize 
\frac {(x+e+1)_{c-e}} {(i-j+c)!}$}}(3,-2)
\Label\o{\raise25pt\hbox{\hskip0pt$\dsize
\times(x-i+j+1)_{e+i-j-1}$}}(4,-4)
\Label\o{\raise15pt\hbox{\hskip3pt$\dsize 
2\frac {(2x+2e+1)_{b-2e}} {(i-j+b)!}$}}(11,-2)
\Label\o{\raise25pt\hbox{\hskip5pt$\dsize
\times(2x-i+j+1)_{2e+i-j-1}$}}(12,-4)
\Label\o{\raise15pt\hbox{$0\!\le\! i\!<\!c$\hskip-5pt}}(19,2)
\Label\o{\raise15pt\hbox{$c\!\le\! i\!<\!b+c-2e$\hskip-5pt}}(19,0)
\Label\o{\raise20pt\hbox{$b+c-2e\!\le\! i\!<\!b$\hskip-5pt}}(19,-3)
\Label\o{,}(17,-4)
\hskip8cm
\)
\hbox{\hskip-8.3cm}
\Label\o{\raise0pt\hbox{$c\le j<b$}}(4,-5)
\Label\o{\raise0pt\hbox{$b\le j<b+c$}}(12,-5)
\hskip10cm\\
\endmultline
\tag3.3$$
\vskip10pt
\noindent
which we denote by $\De_1(x;b,c,e)$. Obviously, we have taken out
$(x+e)^{2e-c}$. The determinant $\De_1(x;b,c,e)$ has still entries which are
polynomial in $x$. For, it is obvious that the entries in rows
$i=0,1,\dots,b+c-2e-1$ are polynomials in $x$, and for $i\ge b+c-2e$ we
have: $c-e\ge0$ by assumption, $e+i-j-1\ge b+c-e-j-1\ge c-e\ge 0$ if
$j<b$, $b-2e\ge b-2c\ge0$ by our assumptions, and $2e+i-j-1\ge
b+c-j-1\ge0$ if $j<b+c$. This explains the term $(2e-c)$ in (3.2).

Now let $e>(b-c)/2$.
In order to explain the term $(2e+c-b)$ in (3.2), we claim that for
$s=0,1,\dots,2e+c-b-1$ we have
$$\multline
 \sum_{i = 0}^{b - 2 e + s}
{{\left( -1 \right) }^{b + c + e + i + s+1}}
\frac {   \left(   b + c - 2 e - i -1\right) !\, 
         \left(   b - e - i -1\right) ! } 
{\left(   2 e - 2 s -2\right) !\, \left(   b - i - s -1\right) !}\\
\cdot     {{ 
      \left(   e - s -1\right) !\, 
         \left(  2e- c  - s -1\right) !}\over 
       {         \left( b - 2 e - i + s \right) !}}
\cdot \big( \text{row $i$ of $\De_1(-e;b,c,e)$}\big)
\\  + 
 \sum_{i = b - e}^{ b + c - 2 e-1}
2\,{{\left( -1 \right) }^{b + c + i}}
\frac {   \left(   b + c - 2 e - i -1\right) !\, \left(   e - s -1\right) !}
{\left( -b + e + i \right) !\, \left(   2 e - 2 s -2\right) ! }
\hskip3cm
\\ \cdot
     {{ 
               \left(  2e- c  - s -1\right) !\, 
         \left(  2e- b  + i - s -1\right) !}\over 
       {
         \left(   b - i - s -1\right) !}} 
\cdot \big( \text{row $i$ of $\De_1(-e;b,c,e)$}\big)
\\ +
   \sum_{i = b + c - 2 e}^{ b - s-1}
     {{({ \textstyle 1 - b - c + 2 e + i}) _{ b - i - s-1}\,  
         ({ \textstyle 1 - e + s}) _{ b - i - s-1}  }\over 
       {\left(   b - i - s -1\right) !\, 
         ({ \textstyle 2 - 2 e + 2 s}) _{ b - i - s-1} }} \\
\cdot \big( \text{row $i$ of $\De_1(-e;b,c,e)$}\big)\\
=0.\hskip8cm
\endmultline\tag3.4$$
Note that these are indeed $2e+c-b$ linear combinations of the rows,
which are linearly independent. The latter fact comes from the
observation that for fixed $s$ 
the last nonzero coefficient in the linear combination (3.4) 
is the one for row $b-s-1$. 

Because of the condition $s\le 2e+c-b-1$, we have $b-2e+s\le c-1$, and
therefore the rows which are involved in the first sum in (3.4) are
from rows $0,1,\dots,c-1$, which form the top block in (3.3). 
The assumptions $e\le c$ and $2c\le b$ imply $b-2e+s\ge 0$, and so the
bounds for the sum are proper bounds. 
Again using the assumptions $b\ge 2c$ and $c\ge e$, we infer $b-e\ge c$, and
therefore the rows which are involved in the second sum in (3.4) are
from rows $c,c+1,\dots,b+c-2e-1$, which form the middle block in (3.3).
The bounds for the sum are proper, since $e\le c$ (including the
possibility that $c=e$, in which case the second sum in (3.4) is the
empty sum).
Finally, because of the condition $s\ge 0$, we have $b-s-1\le b-1$, and
therefore the rows which are involved in the third sum in (3.4) are
from rows $b+c-2e,b+c-2e+1,\dots,b-1$, which form the bottom block in (3.3).
Clearly, the bounds for the sum are proper because of $s\le
2e+c-b-1\le 2e-c-1$.
It is also useful to observe that we need the restriction $(b-c)/2<e$
in order that there is at least one $s$ with $0\le s\le 2e+c-b-1$.

Hence, in order to verify (3.4), we have to check
$$\multline
 \sum_{i = 0}^{b - 2 e + s}
{{\left( -1 \right) }^{b + c + e + i + s+1}}
\frac {         \left(   b + c - 2 e - i -1\right) !\, 
         \left(   b - e - i -1\right) ! } 
{\left(   2 e - 2 s -2\right) !\, \left(   b - i - s -1\right) !}\\
\cdot     {{ \left(   e - s -1\right) !\, 
         \left(  2e- c  - s -1\right) !}\over 
       {  \left( b - 2 e - i + s \right) !}} 
{{c - e}\choose {i - j+c}}
\\ + 
   \sum_{i = b + c - 2 e}^{ b - s-1}
{{\left( -1 \right) }^{e + i + j+1}}
     {{         ({ \textstyle 1 - b - c + 2 e + i}) _{ b - i - s-1}\,  
         ({ \textstyle 1 - e + s}) _{ b - i - s-1} }\over 
       {\left(   b - i - s -1\right) !\, 
         ({ \textstyle 2 - 2 e + 2 s}) _{ b - i - s-1} }} \\
\cdot\frac { \left( c - e \right) !\, \left(   e + i - j -1\right) !} 
{\left( i - j +c\right) ! }\\
=0,\hskip8cm
\endmultline\tag3.5$$
which is (3.4) restricted to the $j$-th column, $j=c,c+1,\dots,b-1$ (note
that all the entries in rows $b-e,b-e+1,\dots,b+c-2e-1$ of
$\De_1(-e;b,c,e)$ vanish
in such a column), 
and
$$\multline
   \sum_{i = 0}^{b - 2 e + s}
{{\left( -1 \right) }^{b + c + e + i + s+1}}     
\frac {      \left(   b + c - 2 e - i -1\right) !\, 
         \left(   b - e - i -1\right) ! } 
{\left(   2 e - 2 s -2\right) !\, \left(   b - i - s -1\right) !}\\
\cdot
{{
   \left(   e - s -1\right) !\, 
         \left(  2e- c  - s -1\right) !}\over 
       {         \left( b - 2 e - i + s \right) !}} 
2{{b - 2 e}\choose {i - j+b}}
\\  + 
 \sum_{i = b - e}^{ b + c - 2 e-1}
2\,{{\left( -1 \right) }^{b + c + i}}
\frac {      \left(   b + c - 2 e - i -1\right) !\, \left(   e - s -1\right) !} 
{\left( -b + e + i \right) !\, \left(  2 e - 2 s -2\right) !}
\hskip3cm
\\
\cdot     {{  
            \left(  2e- c  - s -1\right) !\, 
         \left(  2e- b  + i - s -1\right) !}\over 
       {         \left(   b - i - s -1\right) !}} 
{{b - 2 e}\choose {i - j+b}}
\\ + 
  \sum_{i = b + c - 2 e}^{ b - s-1}
2\, {{\left( -1 \right) }^{i + j+1}}
     {{          ({ \textstyle 1 - b - c + 2 e + i}) _{ b - i - s-1}\,  
         ({ \textstyle 1 - e + s}) _{ b - i - s-1} }\over 
       {\left(   b - i - s -1\right) !\, 
         ({ \textstyle 2 - 2 e + 2 s}) _{ b - i - s-1} }} \\
\cdot \frac {\left( b - 2 e \right) !\, 
         \left(   2 e + i - j -1\right) !} 
{\left( i - j+b \right) !}\\
=0,\hskip8cm
\endmultline\tag3.6$$
which is (3.4) restricted to the $j$-th column, $j=b,b+1,\dots,b+c-1$.

We start by proving (3.5). We remind the reader that here $j$ is
restricted to $c\le j<b$.
The two sums in (3.5) can be combined into a single sum. To be
precise, the left-hand side in (3.5) can be written as
$$\multline
\lim_{\de\to0}\Bigg(
\sum_{i = j-c}^{b-s-1}{{\left( -1 \right) }^{e + i + j+1}}
\frac {       \left( c - e \right) !\, ({ \textstyle 1 + {\de}}) _{  e + i -
j-1}} 
{\left( i - j +c\right) !\, \left(   b - i - s -1\right) !}\\
\cdot
{{          ({ \textstyle 1 - b - c + 2 e + {\de} + i}) _{  b - i -
s-1}\,  
       ({ \textstyle 1 - e + {\de} + s}) _{  b - i - s-1} }\over 
     {       ({ \textstyle 2 - 2 e + {\de} + 2 s}) _{  b - i - s-1}
}}\Bigg).
\endmultline\tag3.7$$
In terms of the usual hypergeometric notation
$${}_r F_s\!\left[\matrix a_1,\dots,a_r\\ b_1,\dots,b_s\endmatrix; 
z\right]=\sum _{k=0} ^{\infty}\frac {\po{a_1}{k}\cdots\po{a_r}{k}}
{k!\,\po{b_1}{k}\cdots\po{b_s}{k}} z^k\ ,$$
where the shifted factorial
$(a)_k$ is given by $(a)_k:=a(a+1)\cdots(a+k-1)$,
$k\ge1$, $(a)_0:=1$, as before, this sum can be rewritten in the form
$$\multline
\kern-5pt 
\lim_{\de\to0}\Bigg(
{{\left( -1 \right) }^{c + e +1}}
{{ 
  ({ \textstyle 1}) _{c - e}\,  
      ({ \textstyle 1 - b - 2 c + 2 e + {\de} + j}) _{  b + c - j -
s-1}\,
        ({ \textstyle 1 - e + {\de} + s}) _{  b + c - j - s-1} }\over 
    {({ \textstyle 1}) _{  b + c - j - s-1}\,  
      ({ \textstyle -c + e + {\de}}) _{c - e+1}\,  
      ({ \textstyle 2 - 2 e + {\de} + 2 s}) _{  b + c - j - s-1} }}\\
\times
  {} _{3} F _{2} \!\left [ \matrix { -c + e + {\de}, -b - c + 2 e -
       {\de} + j - s, 1 - b - c + j + s}\\ { 1 - b - 2 c + 2 e + {\de}
       + j, 1 - b - c + e - {\de} + j}\endmatrix ; {\displaystyle 1}\right
       ]\Bigg).
\endmultline$$
The $_3F_2$-series can be evaluated by means of the
Pfaff--Saalsch\"utz summation (cf\. \cite{\SlatAC, (2.3.1.3); Appendix
(III.2)}),
$$
{} _{3} F _{2} \!\left [ \matrix { A, B, -n}\\ { C, 1 + A + B - C -
   n}\endmatrix ; {\displaystyle 1}\right ]  =
  {{({ \textstyle C-A}) _{n}  ({ \textstyle C-B}) _{n} }\over 
    {({ \textstyle C}) _{n}  ({ \textstyle C-A-B}) _{n} }},
\tag3.8$$
where $n$ is a nonnegative integer. We have to apply the case where
$n=b+c-j-s-1$. Note that this is indeed a nonnegative integer since
$j\le b-1$ and $s\le c-1$. The latter inequality comes from the
assumption $e\le c$ and the inequality chain
$$s\le 2e+c-b-1\le 2e-c-1\le e-1.\tag3.9$$
Thus we obtain, after some
simplification, the expression
$$
\lim_{\de\to0}\Bigg(
{{\left( -1 \right) }^{c + e + 1}}
{{ ({ \textstyle 1}) _{c - e}\,  
      ({ \textstyle 1 - b - c + e + j}) _{  b + c - j - s-1}\,  
      ({ \textstyle 1 - c + 2 {\de} + s}) _{  b + c - j - s-1} }\over 
    {({ \textstyle 1}) _{  b + c - j - s-1}\,  
      ({ \textstyle -c + e + {\de}}) _{c - e+1}\,  
      ({ \textstyle 2 - 2 e + {\de} + 2 s}) _{  b + c - j - s-1} }}
\Bigg)
$$
for the left-hand side in (3.5). This expression vanishes because of
the occurence of the term 
$$    ({ \textstyle 1 - b - c + e + j}) _{  b + c - j - s-1}=
(1 - b - c + e + j)(2 - b - c + e + j)\cdots (e-s-1)$$
in the numerator. For, we have $1-b-c+e+j\le 0$, since $e\le c$ and
$j\le b-1$, and we have $e-s-1\ge 0$, 
thanks to (3.9). This establishes (3.5).

Now we turn to (3.6). We remind the reader that here $j$ is
restricted to $b\le j<b+c$. To begin with, 
we make the similar observation as before that the
three sums on the left-hand side of (3.6) can be combined into a single sum. 
Here, the left-hand side in (3.6) can be written as
$$\multline
\lim_{\de\to0}\Bigg(
\sum_{i = j-b}^{b-s-1}
2\, {{\left( -1 \right) }^{i + j+1}}
\frac {       \left( b - 2 e \right) !\, 
       ({ \textstyle 1 + {\de}}) _{  2 e + i - j-1}} 
{\left(   b - i - s -1\right) !\, ({i - j+b})!}\\
\cdot
{{ 
       ({ \textstyle 1 - b - c + 2 e + {\de} + i}) _{  b - i -
s-1}\,  
       ({ \textstyle 1 - e + {\de} + s}) _{  b - i - s-1} }\over 
     {       ({ \textstyle 2 - 2 e + {\de} + 2 s}) _{  b - i - s-1}
}}\Bigg).
\endmultline\tag3.10$$
In hypergeometric notation, the sum can be rewritten as
$$\multline
\lim_{\de\to0}\Bigg(
2 {{\left( -1 \right) }^{b+1}}
{{ 
     ({ \textstyle 1}) _{b - 2 e}\,  
     ({ \textstyle 1 - 2 b - c + 2 e + {\de} + j}) _{  2 b - j -
s-1}\, 
     ({ \textstyle 1 - e + {\de} + s}) _{  2 b - j - s-1} }\over 
   {({ \textstyle 1}) _{  2 b - j - s-1}\,  
     ({ \textstyle -b + 2 e + {\de}}) _{b - 2 e+1}\,  
     ({ \textstyle 2 - 2 e + {\de} + 2 s}) _{  2 b - j - s-1} }}\\
\times
     {} _{3} F _{2} \!\left [ \matrix { -b + 2 e + {\de}, -2 b + 2 e -
      {\de} + j - s, 1 - 2 b + j + s}\\ { 1 - 2 b + e - {\de} + j, 1 -
      2 b - c + 2 e + {\de} + j}\endmatrix ; {\displaystyle 1}\right
]\Bigg).
\endmultline$$
To this $_3F_2$-series we apply one of Thomae's $_3F_2$-transformation formula 
(cf\. \cite{\BailAA, Ex.~7, p.~98})
$$
{} _{3} F _{2} \!\left [ \matrix { A, B, C}\\ { D, E}\endmatrix ;
   {\displaystyle 1}\right ]  =
\frac {\Gamma(E)\, \Gamma(-A - B - C + D + E)} {\Gamma(-A + E)\, 
\Gamma(-B - C + D + E)}
  {} _{3} F _{2} \!\left [ \matrix { A, -B + D, -C + D}\\ { D, -B - C + D +
    E}\endmatrix ; {\displaystyle 1}\right ]  .
\tag3.11$$
Thus we obtain
$$\multline
\lim_{\de\to0}\Bigg(
2 {{\left( -1 \right) }^{b+1}}
{{ 
     ({ \textstyle 1}) _{b - 2 e}\,  
     ({ \textstyle 1 - 2 b - c + 2 e + {\de} + j}) _{  2 b - j -
s-1}\, 
     ({ \textstyle 1 - e + {\de} + s}) _{  2 b - j - s-1} }\over 
   {     ({ \textstyle 1}) _{  2 b - j - s-1}\,  
     ({ \textstyle -b + 2 e + {\de}}) _{b - 2 e+1}\,  
     ({ \textstyle 2 - 2 e + {\de} + 2 s}) _{  2 b - j - s-1} }}\\
\times
\frac {     \Gamma({ \textstyle 1 - 2 b - c + 2 e + {\de} + j})  \,
   \Gamma({ \textstyle 1 + b - c - e})  } 
{    \Gamma({ \textstyle 1 - b - c + j}) \,
 \Gamma({ \textstyle 1 - c + e + {\de}})   }\\
\times
     {} _{3} F _{2} \!\left [ \matrix { -b + 2 e + {\de}, 1 - e + s, e -
      {\de} - s}\\ { 1 - 2 b + e - {\de} + j, 1 - c + e + 
   {\de}}\endmatrix ; {\displaystyle 1}\right ]\Bigg)
\endmultline$$
for the left-hand side in (3.6). The $_3F_2$-series in this
expression terminates because of the upper parameter $1-e+s$, which
is a nonpositive integer thanks to (3.9). Hence it is well-defined.
The complete expression vanishes because of the occurence of the
term $\Ga({ \textstyle 1 - b - c + j})$ in the denominator. For,
we have $1-b-c+j\le 0$ and so the gamma function equals $\infty$.
This establishes (3.6), and thus completes the proof that $(x+e)$
divides $\De'(x;b,c)$ with multiplicity $m(e)$ as given in (3.2).

\smallskip
{\it Case 2: $c\le e\le b/2$}. By inspection of the expression (3.1),
we see that we have to prove that $(x+e)^{m(e)}$
divides $\De'(x;b,c)$, where
$$m(e)= \cases c&c\le e\le (b-c)/2\\
c+(2e+c-b)&(b-c)/2<e\le b/2.\endcases\tag3.12$$
Note that the first case in (3.12) could be empty, but not the second.

We proceed in a similar manner as before. However, there is a slight
deviation at the beginning. Before we are able to extract the
appropriate number of factors $(x+e)$ out of the determinant
$\De'(x;b,c)$, we have to perform a few row manipulations. 
We add row $b-2$
to row $b-1$, row $b-3$ to row $b-2$, \dots, row $c$
to row $c+1$, in that order. Then we add row $b-2$
to row $b-1$, row $b-3$ to row $b-2$, \dots, row $c+1$
to row $c+2$ (but not row $c$ to row $c+1$!), 
etc. Finally we stop by adding $b-2$
to row $b-1$, row $b-3$ to row $b-2$, \dots, row $e-1$
to row $e$. The resulting determinant is
$$\multline
\hskip-.5cm
\det_{0\le i<b,\, c\le j\le b+c}\!\!\(
\PfadDicke{.3pt}
\SPfad(0,1),1111111111111111\endSPfad
\SPfad(0,-1),1111111111111111\endSPfad
\hbox{\hskip-5pt}
\SPfad(8,-3),222222\endSPfad
\hbox{\hskip5pt}
\Label\o{\raise25pt\hbox{$\dsize\binom {x+c}{i-j+c}$}}(4,1)
\Label\o{\raise25pt\hbox{$\dsize2\binom {2x+b}{i-j+b}$}}(12,1)
\Label\o{\raise25pt\hbox{$\dsize\binom {x+i}{i-j+c}$}}(4,-1)
\Label\o{\raise25pt\hbox{$\dsize\binom {2x+b-c+i}{i-j+b}$}}(12,-1)
\Label\o{\raise25pt\hbox{$\dsize\binom {x+e}{i-j+c}$}}(4,-3)
\Label\o{\raise25pt\hbox{$\dsize\binom {2x+b-c+e}{i-j+b}$}}(12,-3)
\Label\o{\raise15pt\hbox{$0\!\le\! i\!<\!c$\hskip-5pt}}(19,1)
\Label\o{\raise15pt\hbox{$c\!\le\! i\!<\!e$\hskip-5pt}}(19,-1)
\Label\o{\raise20pt\hbox{$e\!\le\! i\!<\!b$\hskip-5pt}}(19,-3)
\Label\o{.}(17,-4)
\hskip8cm
\)
\hbox{\hskip-8.3cm}
\Label\o{\raise0pt\hbox{$c\le j<b$}}(4,-5)
\Label\o{\raise0pt\hbox{$b\le j<b+c$}}(12,-5)
\hskip10cm\\
\endmultline\tag3.13
$$
\vskip10pt
Now we take $(x+e)$ out of rows $b-c,b-c+1,\dots,b-1$, and obtain the
determinant
$$\multline
\hskip-.5cm
\det_{0\le i<b,\, c\le j\le b+c}\!\!\(
\PfadDicke{.3pt}
\SPfad(0,3),1111111111111111\endSPfad
\SPfad(0,1),1111111111111111\endSPfad
\SPfad(0,-1),1111111111111111\endSPfad
\hbox{\hskip-5pt}
\SPfad(8,-5),2222222222\endSPfad
\hbox{\hskip5pt}
\Label\o{\raise25pt\hbox{$\dsize\binom {x+c}{i-j+c}$}}(4,3)
\Label\o{\raise25pt\hbox{$\dsize2\binom {2x+b}{i-j+b}$}}(12,3)
\Label\o{\raise25pt\hbox{$\dsize\binom {x+i}{i-j+c}$}}(4,1)
\Label\o{\raise25pt\hbox{$\dsize\binom {2x+b-c+i}{i-j+b}$}}(12,1)
\Label\o{\raise25pt\hbox{$\dsize\binom {x+e}{i-j+c}$}}(4,-1)
\Label\o{\raise25pt\hbox{$\dsize\binom {2x+b-c+e}{i-j+b}$}}(12,-1)
\Label\o{\raise15pt\hbox{\hskip-10pt$\dsize 
\frac {(x\!+\!e\!-\!c\!-\!i\!+\!j\!+\!1)_{c\!+\!i\!-\!j\!-\!1}}
{(i-j+c)!}$}}(4,-3)
\Label\o{\raise15pt\hbox{\hskip3pt$\dsize 
2\frac {(2x+2e+1)_{b-c-e}} {(i-j+b)!}$}}(11,-3)
\Label\o{\raise25pt\hbox{\hskip5pt$\dsize
\times(2x\!\kern-1pt +\!\kern-1pt e\!\kern-1pt -\!\kern-1pt c\!\kern-1pt 
-\!\kern-1pt i\!\kern-1pt +\!\kern-1pt j\!\kern-1pt +\!\kern-1pt 1)
_{c\!+\!e\!+\!i\!-\!j\!-\!1}$}}(12,-5)
\Label\o{\raise15pt\hbox{$0\!\le\! i\!<\!c$\hskip-5pt}}(19,3)
\Label\o{\raise15pt\hbox{$c\!\le\! i\!<\!e$\hskip-5pt}}(19,1)
\Label\o{\raise20pt\hbox{$e\!\le\! i\!<\!b-c$\hskip-5pt}}(19,-1)
\Label\o{\raise20pt\hbox{$b-c\!\le\! i\!<\!b$\hskip-5pt}}(19,-4)
\Label\o{,}(17,-5)
\hskip8cm
\)
\hbox{\hskip-8.3cm}
\Label\o{\raise0pt\hbox{$c\le j<b$}}(4,-6)
\Label\o{\raise0pt\hbox{$b\le j<b+c$}}(12,-6)
\hskip10cm\\
\endmultline
\tag3.14$$
\vskip10pt
\noindent
which we denote by $\De_2(x;b,c,e)$. Obviously, we have taken out
$(x+e)^{c}$. The remaining determinant has still entries which are
polynomial in $x$. For, it is obvious that the entries in rows
$i=0,1,\dots,b-c-1$ are polynomials in $x$, and for $i\ge b-c$ we
have: $c+i-j-1\ge b-j-1\ge0$ if $j<b$,
$b-c-e\ge b-2e\ge0$ by our assumptions, and $c+e+i-j-1\ge
b+e-j-1\ge b+c-j-1\ge0$ if $j<b+c$. 
 This explains the term $c$ in (3.12).

Now let $e>(b-c)/2$.
In order to explain the term $(2e+c-b)$ in (3.12), we claim that for
$s=0,1,\dots,2e+c-b-1$ we have
$$\align
&\sum _{i=0} ^{b-2e+s}
{{\left( -1 \right) }^{b + c + e + i + s+1}}{{ 
     \left( b + c - 2 e - i -1 \right) !\, \left( b - e - i -1 \right) !\, 
     }\over 
   {\left(  2 e - 2 s -2 \right) !\, \left( b - i - s -1 \right) !}}\\
&\hskip3cm
\cdot \frac {\left(  e - s -1 \right) !\, \left( 2e- c - s -1 \right) !} 
{\left( b - 2 e - i + s \right) !}
\cdot \big(\text {row $i$ of $\De_2(-e;b,c,e)$}\big)
\\
&+
\sum _{i=c} ^{e-1} \Bigg(
   \sum_{k = 0}^{i-c}{{\left( -1 \right) }^{b + e + s}}
{{2e- b + i - s-1 }\choose   {i -c - k}}
{{ \left( e- c - k -1 \right) !}
        \over {k!}}  
\hskip3.5cm
\\
&\hskip2cm
\cdot \frac {\left( c + k - s -1 \right) !\, \left( e + k - s -1 \right) !} 
{\left( 2 e - 2 s -2 \right) !}\Bigg)
\cdot \big(\text {row $i$ of $\De_2(-e;b,c,e)$}\big)
\\
&+
\sum _{i=e} ^{b-e-1} \Bigg(
\sum_{k = 0}^{e- c-1 }{{\left( -1 \right) }^{b + i + s+1}}
{{ 
       \left( 2e- b + i - s -1 \right) !\, 
       \left( c + k - s -1 \right) !\, \left( e + k - s -1 \right) !}\over 
     {k!\, \left( 2 e - 2 s -2 \right) !\, 
       \left( c + e -b + i + k - s \right) !}}\Bigg)\\
&\hskip5cm
\cdot \big(\text {row $i$ of $\De_2(-e;b,c,e)$}\big)
\\
&+
\sum _{i=b-e} ^{b-c-1} \Bigg(   \sum_{k = 0}^{e- c-1 }
{{\left( -1 \right) }^{b + i + s+1}}
     {{ 
         \left( 2e- b + i - s -1 \right) !\, 
         \left( c + k - s -1 \right) !\, \left( e + k - s -1 \right) !}
        \over {k!\, \left( 2 e - 2 s -2 \right) !\, 
         \left(  c + e - b + i + k - s \right) !}} \\
& + 
{{\left( -1 \right) }^{b + c + i}}
{{ \left( b - c - i -1 \right) !\, 
      \left( c - s -1 \right) !\, \left( 2e- b + i - s -1 \right) !\, 
      ({ \textstyle b - i - s}) _{i-b + e} }\over 
    {\left( i-b + e \right) !\, \left( 2 e - 2 s -2 \right) !}} \Bigg) \\
&\hskip5cm
\cdot \big(\text {row $i$ of $\De_2(-e;b,c,e)$}\big)
\\
&+
\sum _{i=b-c} ^{b-s-1}
{{({ \textstyle 1 - b + c + i}) _{ b - i - s-1 }  \,
     ({ \textstyle 1 - b + e + i}) _{b - i - s-1 } }\over 
   {\left( b - i - s -1 \right) !\, 
     ({ \textstyle 2e-b + i - s}) _{b - i - s-1 } }}
\cdot \big(\text {row $i$ of $\De_2(-e;b,c,e)$}\big)
\\
&\hskip4cm=0.\tag3.15
\endalign$$
Again, note that these are indeed $2e+c-b$ linear combinations of the rows,
which are linearly independent. 

Because of the condition $s\le 2e+c-b-1$, we have $b-2e+s\le c-1$, and
therefore the rows which are involved in the first sum in (3.15) are
from rows $0,1,\dots,c-1$, which form the top block in (3.14). The
assumption $e\le b/2$ guarantees that the bounds for the sum 
are proper bounds. Clearly, the rows which are involved in the second sum in
(3.15) are the rows $c,c+1,\dots,e-1$, which form the second block
from top in (3.14). The assumption $c\le e$
guarantees that the bounds for the sum are proper bounds,
(including the
possibility that $c=e$, in which case the sum is the
empty sum). Because of the same assumptions, the rows which are
involved in the third and fourth sum in (3.15) are from rows
$e,e+1,\dots,b-c-1$, which form the third block from top in (3.14). The
assumption $e\le b/2$ guarantees that the third sum in (3.15) has proper
bounds (including the
possibility that $e=b/2$, in which case the sum is the empty sum).
Finally, because of the condition $s\ge 0$, we have $b-s-1\le b-1$, and
therefore the rows which are involved in the fifth sum in (3.15) are
from rows $b-c,b-c+1,\dots,b-1$, which form the bottom block in (3.14).
The bounds for this fifth sum are proper because $s\le c-1$.
This inequality follows from the inequality chain
$$s\le 2e+c-b-1\le b+c-b-1=c-1.\tag3.16$$
Again, it is also useful to observe that we need the restriction $(b-c)/2<e$
in order that there is at least one $s$ with $0\le s\le 2e+c-b-1$.

Hence, in order to verify (3.15), we have to check
$$\align
&\sum _{i=0} ^{b-2e+s}
{{\left( -1 \right) }^{b + c + e + i + s+1}}{{ 
     \left( b + c - 2 e - i -1 \right) !\, \left( b - e - i -1 \right) !\, 
     }\over 
   {\left(  2 e - 2 s -2 \right) !\, \left( b - i - s -1 \right) !}}\\
&\hskip3cm
\cdot \frac {\left(  e - s -1 \right) !\, \left( 2e- c - s -1 \right) !} 
{\left( b - 2 e - i + s \right) !}
\binom {c-e}{i-j+c}
\tag3.17a\\
&+
\sum _{i=c} ^{e-1} \Bigg(
   \sum_{k = 0}^{i-c}{{\left( -1 \right) }^{b + c+e+i+k + s}}
{{b -c-e }\choose   {i -c - k}} {{e-s-1} \choose {k}}
\hskip3.5cm
\\
&\hskip2cm
\cdot \frac {\left( e- c - k -1 \right) !\,
\left( c - s -1 \right) !\, \left( e + k - s -1 \right) !} 
{\left( 2 e - 2 s -2 \right) !}\Bigg)
\binom {i-e}{i-j+c}
\tag3.17b\\
&+
\sum _{i=e} ^{b-c-1} \Bigg(
\sum_{k = 0}^{e- c-1 }{{\left( -1 \right) }^{b + i + s+1}}
{{ 
       \left( 2e- b + i - s -1 \right) !\, 
       \left( c + k - s -1 \right) !\, \left( e + k - s -1 \right) !}\over 
     {k!\, \left( 2 e - 2 s -2 \right) !\, 
       \left( c + e -b + i + k - s \right) !}}\Bigg)\\
&\hskip5cm
\cdot\binom {0}{i-j+c}
\tag3.17c\\
& + 
\sum _{i=b-e} ^{b-c-1}
{{\left( -1 \right) }^{b + c + i}}
{{ \left( b - c - i -1 \right) !\, 
      \left( c - s -1 \right) !\, \left( 2e- b + i - s -1 \right) !\, 
      ({ \textstyle b - i - s}) _{i-b + e} }\over 
    {\left( i-b + e \right) !\, \left( 2 e - 2 s -2 \right) !}} \\
&\hskip5cm
\cdot\binom {0}{i-j+c}
\tag3.17d\\
&+
\sum _{i=b-c} ^{b-s-1}
{{({ \textstyle 1 - b + c + i}) _{ b - i - s-1 }  \,
     ({ \textstyle 1 - b + e + i}) _{b - i - s-1 } }\over 
   {\left( b - i - s -1 \right) !\, 
     ({ \textstyle 2e-b + i - s}) _{b - i - s-1 } }}
(-1)^{c+i+j+1}\frac {1} {(i-j+c)}
\tag3.17e\\
&\hskip4cm=0,\tag3.17f
\endalign$$
which is (3.15) restricted to the $j$-th column, $j=c,c+1,\dots,b-1$, 
and
$$\align
&\sum _{i=0} ^{b-2e+s}
{{\left( -1 \right) }^{b + c + e + i + s+1}}{{ 
     \left( b + c - 2 e - i -1 \right) !\, \left( b - e - i -1 \right) !\, 
     }\over 
   {\left(  2 e - 2 s -2 \right) !\, \left( b - i - s -1 \right) !}}\\
&\hskip3cm
\cdot \frac {\left(  e - s -1 \right) !\, \left( 2e- c - s -1 \right) !} 
{\left( b - 2 e - i + s \right) !}
2\binom {b-2e}{i-j+b}
\tag3.18a\\
&+
\sum _{i=c} ^{e-1} \Bigg(
   \sum_{k = 0}^{i-c}{{\left( -1 \right) }^{b + c+e+i+k + s}}
{{b -c-e }\choose   {i -c - k}} {{e-s-1} \choose {k}}
\hskip3.5cm
\\
&\hskip2cm
\cdot \frac {\left( e- c - k -1 \right) !\,
\left( c - s -1 \right) !\, \left( e + k - s -1 \right) !} 
{\left( 2 e - 2 s -2 \right) !}\Bigg)
\binom {i+b-c-2e}{i-j+b}
\tag3.18b\\
&+
\sum _{i=e} ^{b-c-1} \Bigg(
\sum_{k = 0}^{e- c-1 }{{\left( -1 \right) }^{b + i + s+1}}
{{ 
       \left( 2e- b + i - s -1 \right) !\, 
       \left( c + k - s -1 \right) !\, \left( e + k - s -1 \right) !}\over 
     {k!\, \left( 2 e - 2 s -2 \right) !\, 
       \left( c + e -b + i + k - s \right) !}}\Bigg)\\
&\hskip5cm
\cdot\binom {b-c-e}{i-j+b}
\tag3.18c\\
& + 
\sum _{i=b-e} ^{b-c-1} 
{{\left( -1 \right) }^{b + c + i}}
{{ \left( b - c - i -1 \right) !\, 
      \left( c - s -1 \right) !\, \left( 2e- b + i - s -1 \right) !\, 
      ({ \textstyle b - i - s}) _{i-b + e} }\over 
    {\left( i-b + e \right) !\, \left( 2 e - 2 s -2 \right) !}} \\
&\hskip5cm
\cdot\binom {b-c-e}{i-j+b}
\tag3.18d\\
&+
\sum _{i=b-c} ^{b-s-1}
{{({ \textstyle 1 - b + c + i}) _{ b - i - s-1 }  \,
     ({ \textstyle 1 - b + e + i}) _{b - i - s-1 } }\over 
   {\left( b - i - s -1 \right) !\, 
     ({ \textstyle 2e-b + i - s}) _{b - i - s-1 } }}\\
&\hskip2cm \cdot2(-1)^{e+c+i+j+1}\frac {(b-c-e)!\,(i-j+c+e-1)!} {(i-j+b)!}
\tag3.18e\\
&\hskip4cm=0,\tag3.18f
\endalign$$
which is (3.15) restricted to the $j$-th column, $j=b,b+1,\dots,b+c-1$.

We start by proving (3.17). We remind the reader that here $j$ is
restricted to $c\le j<b$. Apparently, (3.17) is more
complex than (3.5) or (3.6), so it is not surprising that the
arguments here are more complex than the arguments for proving (3.5)
and (3.6). It turns out that the five terms in (3.17) cannot be
combined into one term, as was the case for (3.5) and (3.6). Rather we will
combine (3.17a), (3.17d), and (3.17e) into one term, (3.19), 
then we will
combine (3.17b) and (3.17c) into another term, (3.24), and then show how to
transform one of the two into the negative of the other.

So, we claim that the sum of (3.17a), (3.17d), and (3.17e) equals
$$
\lim_{\de\to 0}\Bigg(
\sum_{i = j-c}^{b-s-1} 
{{       ({ \textstyle 1 - b + c - { \de} + i}) _{  b - i - s-1}\,  
       ({ \textstyle 1 - b + e + { \de} + i}) _{  b - i - s-1} }\over 
     {{ \de} \left(   b - i - s -1\right) !\, 
       ({ \textstyle 2 e -b+ { \de} + i - s}) _{  b - i - s-1} }}
{{{ \de}}\choose {i - j+c}}\Bigg).
\tag3.19$$
It is straight-forward to check that (3.17e) agrees with the according
part $ i=b-c,b-c+1,\dots, b-s-1$ of (3.19), and that (3.17d) agrees with the
according part $ i=b-e,b-e+1,\dots, b-c-1$ of (3.19), and that the terms for
$i=b-2e+s+1,b-2e+s+2,\dots,e-1$ in (3.19) vanish. It remains to be
seen that (3.17a) 
agrees with the according part $i=j-c,j-c+1,\dots,b-2e+s$ of (3.19), which
is not directly evident.

In order to verify the last assertion, we replace $\binom {c-e}{i-j+c}$ in
(3.17a) by the expansion $\sum _{\ell=j-c} ^{i}\binom
0{\ell-j+c}\binom {c-e}{i-\ell}$. That the binomial equals the
expansion is due to the Chu--Vandermonde summation. Then we interchange sums
over $i$ and $\ell$, and write the now inner sum over $i$ in
hypergeometric notation. This gives for (3.17a) the expression
$$\multline
\hskip-3pt\sum _{\ell=j-c} ^{b-2e+s}
{{\left( -1 \right) }^{b + c + e + \ell  + s+1}} 
\binom 0{\ell-j+c}
  {{    (  b + c - 2 e - \ell -1)!\,  
      (  b - e - \ell -1)!  
(  e - s-1)!}\over 
    {(  2 e - 2 s-2)! }}
\\
\cdot
\frac {      (  c + 2 e - s-1)! } 
{      (  b - \ell  - s-1)!\,  
      (b - 2 e - \ell  + s)!}
      {} _{3} F _{2} \!\left [ \matrix { 1 - b + \ell  + s, -c + e, -b + 2 e + \ell 
       - s}\\ { 1 - b - c + 2 e + \ell , 1 - b + e + \ell }\endmatrix ;
       {\displaystyle 1}\right ]  .
\endmultline$$
The $_3F_2$-series can be evaluated by means of the
Pfaff--Saalsch\"utz summation (3.8). Thus, the expression for (3.17a)
simplifies to
$$\multline
\sum _{\ell=j-c} ^{b-2e+s}
{{\left( -1 \right) }^{b + c + e + \ell  + s+1}} 
\binom 0{\ell-j+c}\\
\cdot
  {{    \left(   b - c - \ell  -1\right) !\, 
      \left(   b + c - 2 e - \ell  -1\right) !\, 
      \left(   b - e - \ell  -1\right) !\, \left(   e - s -1\right) !}\over 
    {\left(   2 e - 2 s -2\right) !\, \left(   b - \ell  - s -1\right) !\, 
      \left( b - 2 e - \ell  + s \right) !\, 
      ({ \textstyle c - s}) _{b - 2 e - \ell  + s} }}.
\endmultline$$
Now it is straight-forward to check that this agrees with the
part $i=j-c,j-c+1,\dots,b-2e+s$ of (3.19).

\smallskip
Next we consider (3.17b) and (3.17c). We begin by replacing $\binom
{i-e}{i-j+c}$ in (3.17b) by the expansion $\sum _{\ell=j-c} ^{i}\binom
{0}{\ell-j+c} \binom {i-e}{i-\ell}$, the equality of binomial and
expansion being again due to Chu--Vandermonde summation. Then we
interchange the summations over $i$, $k$, $\ell$ so that the sum over
$\ell$ becomes the outer sum and the sum over $i$ becomes the inner
sum, and write the sum
over $i$ in hypergeometric notation. This gives
$$\multline
\hskip-13pt
\sum _{\ell=j-c} ^{e-1}
{{\left( -1 \right) }^{b + c + e + k + \ell  + s}} 
\binom 0{\ell-j+c}
\sum _{k=0} ^{e-c-1}
  {{   \left(  e- c  - k -1\right) !\, \left(   c - s -1\right) !\, 
     \left(   e + k - s -1\right) !}\over 
   {k!\, \left( \ell-c - k   \right) !}}
\\
\cdot
\frac {     ({ \textstyle 1 + b - e + k - \ell }) _{\ell-c - k  }\,  
     ({ \textstyle e - k - s}) _{k} }
 {\left(   2 e - 2 s -2\right) !}
     {} _{2} F _{1} \!\left [ \matrix { -b + e - k + \ell , 1 - e + \ell }\\ { 1 - c -
      k + \ell }\endmatrix ; {\displaystyle 1}\right ]  
\endmultline$$
as an equivalent expression for (3.17b). The $_2F_1$-series can be
evaluated by the hypergeometric form of the Chu--Vandermonde summation 
(see \cite{\SlatAC, (1.7.7); Appendix (III.4)}),
$$
{} _{2} F _{1} \!\left [ \matrix { A, -n}\\ { C}\endmatrix ; {\displaystyle
   1}\right ]  = {{({ \textstyle C-A}) _{n} }\over 
    {({ \textstyle C}) _{n} }},
\tag3.20$$
where $n$ is a nonnegative integer. In the resulting 
inner sum over $k$ we reverse the order of summation, i.e., we replace
$k$ by $e-c-1-k$, and then write the (new) sum over $k$ in
hypergeometric notation. We obtain the expression
$$\multline
\sum _{\ell=j-c} ^{e-1}
{{\left( -1 \right) }^{  b   + \ell  + s+1}} 
\binom 0{\ell-j+c}
  {{   (  c - s-1)!\,  
     ( 2e- c  - s-2)!\,  
     ({ \textstyle 1 + c - s}) _{ e- c -1} }\over 
   {( e- c -1)!\,  (  2 e - 2 s-2)! }}
\\
\cdot
     {} _{3} F _{2} \!\left [ \matrix { 1 - b + c + \ell , 1, 1 + c - e}\\ { 2 + c
      - 2 e + s, 1 + c - s}\endmatrix ; {\displaystyle 1}\right ]  .
\endmultline$$
To the $_3F_2$-series we apply another of Thomae's transformation
formulas (see \cite{\GaRaAA, (3.1.1)}),
$$
{} _{3} F _{2} \!\left [ \matrix { A, B, -n}\\ { D, E}\endmatrix ;
   {\displaystyle 1}\right ]  =
{{ ({ \textstyle -B + E}) _{n} }
     \over {({ \textstyle E}) _{n} }}  
{} _{3} F _{2} \!\left [ \matrix { -n, B, -A + D}\\ { D, 1 + B - E -
       n}\endmatrix ; {\displaystyle 1}\right ] ,
\tag3.21$$
where $n$ is a nonnegative integer. We have to apply the case where
$n=e-c-1$. Because of our assumption $c\le e$ this is indeed a
nonnegative integer, except if $e=c$. So, let us for the moment
exclude the case $e=c$.
After little manipulation, application of (3.21)
yields the following expression for (3.17b):
$$\multline 
\sum _{\ell=j-c} ^{e-1} 
{{\left( -1 \right) }^{b + \ell + s+1}}
\binom {0}{\ell-j+c}\\
\cdot
\sum_{k = 0}^{e- c-1 }
{{ 
       \left(  e-k - s -2 \right) !\, \left( 2e-c-k - s -2 \right) !\,
     (2e-b+\ell-k-s)_{k}}\over 
     {(e-c-k-1)!\, \left( 2 e - 2 s -2 \right) !}}.
\endmultline\tag3.22$$
This expression is equal to (3.17b) for $e=c$ as well since in that
case both expressions are zero due to empty summations over $k$. So,
in all possible cases (3.22) is equal to (3.17b).

The inner sum over $k$ in (3.22) 
is exactly the same as the inner sum over $k$
in (3.17c) when the order of summation is reversed, i.e., when $k$ is
replaced by $e-c-1-k$. This shows that (3.17b) and (3.17c) can be
combined into the single expression
$$\multline 
\sum _{i=j-c} ^{b-c-1} \Bigg(
\sum_{k = 0}^{e- c-1 }{{\left( -1 \right) }^{b + i + s+1}}
{{ 
       \left( c + k - s -1 \right) !\, \left( e + k - s -1 \right) !}\over 
     {k!}}\\
\hskip3cm
\cdot
\frac {     (c+e-b+i+k-s+1)_{e-c-k-1}} 
{\left( 2 e - 2 s -2 \right) !}\Bigg)
\binom {0}{i-j+c},
\endmultline\tag3.23$$
which of course equals
$$
\sum_{k = 0}^{e- c-1 }{{\left( -1 \right) }^{b + j+c + s+1}}
{{ 
       \left( c + k - s -1 \right) !\, \left( e + k - s -1 \right) !\,
   (e-b+j+k-s+1)_{e-c-k-1}}\over 
     {k!\, \left( 2 e - 2 s -2 \right) !}},
\tag3.24$$
since the binomial in (3.23) is 1 only for $i=j-c$ and 0 otherwise. In this
regard, it is
important that the range of summation over $i$ in (3.23)
is in fact not empty (so that the term for $i=j-c$ does indeed occur
in the sum (3.23); otherwise, the previous conclusion would have
been wrong) because for (3.17) we are considering a $j$ with $j\le b-1$.

\smallskip
Our computations so far allow the conclusion that, in order to
establish (3.17), we have to show that (3.19) and (3.24) add up to
zero.

\smallskip
In order to see this, we start with the expression (3.19). In the sum
over $i$, we reverse the order of summation, i.e., we replace $i$ by
$b-s-1-i$, and then write the new sum in hypergeometric notation, to
obtain
$$\multline
\lim_{\de\to0}\Bigg(
{{\left( -1 \right) }^{  b + c + j + s}} 
  {{   {{({ \textstyle 1-{ \de}}) _{  b + c - j - s-2}} \over  
   {(  b + c - j - s-1)!}} }}
\\
\times
     {} _{3} F _{2} \!\left [ \matrix { 1 - b - c + j + s, 1 - c + { \de} +
      s, 1 - e - { \de} + s}\\ { 2 - b - c + { \de} + j + s, 2 - 2 e -
      { \de} + 2 s}\endmatrix ; {\displaystyle 1}\right ]  \Bigg).
\endmultline$$
To the $_3F_2$-series we apply yet another of Thomae's transformation
formulas (see \cite{\SlatAC, (2.3.3.7)}),
$$\multline
{} _{3} F _{2} \!\left [ \matrix { A, B, C}\\ { D, E}\endmatrix ;
   {\displaystyle 1}\right ]  =
   \frac {\Gamma  (D)\,\Gamma( E)\,\Gamma( -A - B - C + D + E)} 
{\Gamma( B)\,\Gamma( -A - B + D + E)\,\Gamma( -B - C
    + D + E)}\\
\times
  {} _{3} F _{2} \!\left [ \matrix { -B + D, -B + E, -A - B - C + D + E}\\ {
    -A - B + D + E, -B - C + D + E}\endmatrix ; {\displaystyle 1}\right ]  .
\endmultline\tag3.25$$
Thus we obtain
$$\multline
\lim_{\de\to0}\Bigg(
{{\left( -1 \right) }^{ b + c + j + s}} 
  {{  ({ \textstyle e-c }) _{b + c - j - s-1} \,
      ({ \textstyle 1 - {\de}}) _{ b + c - j - s-2} }\over 
    {({ \textstyle 1}) _{ b + c - j - s-1} \,
      ({ \textstyle c - {\de} - s}) _{ b - j-1} \,
      ({ \textstyle -1 - c + 2 e + {\de} - s}) _{c - s} }}
\\
\times
   {} _{3} F _{2} \!\left [ \matrix { 1 + c - 2 e - 2 {\de} + s, 1 + c
       - e, 1 - b + j}\\ { 2 + c - 2 e - {\de} + s, 2 - b - e + j +
       s}\endmatrix ; {\displaystyle 1}\right ]  \Bigg)
\endmultline\tag3.26$$
as an equivalent expression for (3.19), after some simplification. The
$_3F_2$-series in this expression is terminating because of the upper
parameter $1-b+j$, which is a nonpositive integer due to $j\le b-1$,
so the $_3F_2$-series is well-defined. The complete expression
vanishes for $e=c$ because of the occurence of the term
$(e-c)_{b+c-j-s-1}$ in the numerator, for, we have $b+c-j-s-1>0$ since
$j\le b-1$ and $s<c$ (cf\. (3.16)). As we did already once, let us for
the moment exclude the case $e=c$.

Next, to the $_3F_2$-series in (3.26), 
we apply the transformation (3.21) (with $n=b-j-1$, which is indeed a
nonnegative integer as we noted just before), obtaining
$$\multline
\lim_{\de\to0}\Bigg(
{{\left( -1 \right) }^{b + c + j + s}} 
  {{  ({ \textstyle e-c }) _{ b + c - j - s-1} \,
      ({ \textstyle 1 - {\de}}) _{ b + c - j - s-2} }\over 
    {({ \textstyle 1}) _{ b + c - j - s-1} \,
      ({ \textstyle c - {\de} - s}) _{ b - j-1} \,
      ({ \textstyle  2e- c + {\de} - s-1}) _{c - s} }}
\\
\times
\frac {      ({ \textstyle 1 - b - c + j + s}) _{ b - j-1} } 
{      ({ \textstyle 2 - b - e + j + s}) _{ b - j-1} }
    {} _{3} F _{2} \!\left [ \matrix { 1 - b + j, 1 + {\de}, 1 + c - e}\\
       { 2 + c - 2 e - {\de} + s, 1 + c - s}\endmatrix ; {\displaystyle
       1}\right ]\Bigg),
\endmultline$$
and apply (3.21) once more (here we need that $e-c-1$ is nonegative,
which is only the case if $e>c$), obtaining
$$\multline
\lim_{\de\to0}\Bigg(
{{\left( -1 \right) }^{b + c + j + s}} 
\frac {({ \textstyle e-c }) _{ b + c - j - s-1} \,
      ({ \textstyle 1 - {\de}}) _{ b + c - j - s-2} }
 {({ \textstyle 1}) _{ b + c - j - s-1} \,
      ({ \textstyle 1 + c - s}) _{e-c-1 } }
\\
\times
  {{  ({ \textstyle c - {\de} - s}) _{e- c-1 } \,
      ({ \textstyle 1 - b - c + j + s}) _{ b - j-1} }\over 
    { ({ \textstyle c - {\de} - s}) _{ b - j-1} \,
      ({ \textstyle 2e - c + {\de} - s-1}) _{c - s} \,
      ({ \textstyle 2 - b - e + j + s}) _{ b - j-1} }}
\\
\times
    {} _{3} F _{2} \!\left [ \matrix { 1 + c - e, 1 + {\de}, 1 + b + c -
       2 e - {\de} - j + s}\\ { 2 + c - 2 e - {\de} + s, 2 - e + {
       \de} + s}\endmatrix ; {\displaystyle 1}\right ]  \Bigg).
\endmultline$$
Now we write the $_3F_2$-series explicitly as a sum over $k$ and perform
the termwise limit $\de\to0$. This gives, after some simplification,
$$
\sum_{k = 0}^{e- c-1 }{{\left( -1 \right) }^{b + j+c + s}}
{{ 
       \left( e-k - s -2 \right) !\, \left( 2e-c-k - s -2 \right) !\,
   (2e-b-c+j-k-s)_{k}}\over 
     {(e-c-k-1)!\, \left( 2 e - 2 s -2 \right) !}},
$$
which is exactly the negative of 
the sum (3.24) in reverse order, i.e., with $k$
replaced by $e-c-1-k$. Hence, the expressions (3.19) and (3.24) add up
to zero. This is also true for $e=c$, since, via
(3.26), we saw that in that case (3.19) vanishes; and so does (3.24)
because of the empty summation over $k$. This establishes the equation
(3.17). 

\smallskip
Now we turn to (3.18). We remind the reader that here $j$ is
restricted to $b\le j<b+c$. We pursue a similar strategy. 
We combine (3.18a), (3.18d), and (3.18e) into one term, and we combine
(3.18b) and (3.18c) into another term. 
Here, it turns out that each combination itself is already zero.

In the same way as before, it is seen that the sum of 
(3.18a), (3.18d), and (3.18e) equals 
$$
\lim_{\de\to 0}\Bigg(
\sum_{i = j-b}^{b-s-1} 
2\,{{       ({ \textstyle 1 - b + c - { \de} + i}) _{  b - i - s-1}\,  
       ({ \textstyle 1 - b + e + { \de} + i}) _{  b - i - s-1} }\over 
     {{ \de} \left(   b - i - s -1\right) !\, 
       ({ \textstyle 2 e -b+ { \de} + i - s}) _{  b - i - s-1} }}
{{{ \de+b-c-e}}\choose {i - j+b}}\Bigg).
\tag3.27$$
This follows by using the same arguments as those that showed that 
the sum of (3.17a), (3.17d), and (3.17e) equals (3.19), the only
deviation is that $\binom 0{i-j+c}$ has to be replaced by $\binom
{b-c-e}{i-j+b}$ everywhere. Actually, the term (3.18d) need not be
considered since it vanishes because $b-c-e<i-j+b$, and therefore the
binomial in (3.18d) vanishes. The inequality is a consequence of $i\ge
b-e$ in the sum (3.18d) and $j<b+c$.

Now we write the sum in (3.27) in hypergeometric notation,
$$\multline 
\lim_{\de\to0}\Bigg(2\,
{{   ({ \textstyle 1 - 2 b + c - {\de} + j}) _{2b-j-s-1} \,
     ({ \textstyle 1 - 2 b + e + {\de} + j}) _{2b-j-s-1} }\over 
   {{\de}\,({ \textstyle 1}) _{2b-j-s-1} \,
     ({ \textstyle -2 b + 2 e + {\de} + j - s}) _{2b-j-s-1} }}\\
\times
{} _{3} F _{2} \!\left [ \matrix { -2 b + 2 e + {\de} + j - s,
  1 - 2 b + j + s, -b + c + e - {\de}
   }\\ { 1 - 2 b + c - {\de} + j, 1 -
      2 b + e + {\de} + j}\endmatrix ; {\displaystyle 1}\right ]  \Bigg).
\endmultline$$
and then apply the transformation formula (3.11), to get
$$\multline 
\lim_{\de\to0}\Bigg(2\,
{{   ({ \textstyle 1 - 2 b + c - {\de} + j}) _{ 2 b - j - s-1} \,
     ({ \textstyle 1 - 2 b + e + {\de} + j}) _{ 2 b - j - s-1} }\over 
   {{\de} \,    ({ \textstyle 1}) _{ 2 b - j - s-1} \,
     ({ \textstyle -2 b + 2 e + {\de} + j - s}) _{2b-j-s-1} }}\\
\times\frac {  \Gamma({ \textstyle 1 + b - 2 e}) \,
     \Gamma({ \textstyle 1 - 2 b + e + {\de} + j}) }
 {\Gamma({ \textstyle 1 - b + {\de} + j - s}) \,
     \Gamma({ \textstyle 1 - e + s}) }\\
\times
{} _{3} F _{2} \!\left [ \matrix { -2 b + 2 e + {\de} + j - s, c - {
      \de} - s, 1 - b - e + j}\\ { 1 - 2 b + c - {\de} + j, 1 - b + {
      \de} + j - s}\endmatrix ; {\displaystyle 1}\right ] \Bigg) ,
\endmultline$$
The $_3F_2$-series in this expression terminates because of the upper
parameter $1-b-e+j$, which is a nonpositive integer since $j<b+c\le
b+e$. Hence it is well-defined.
The complete expression vanishes, because of the occurence of the term
$\Ga(1-e+s)$ in the denominator. For, we have $1-e+s\le0$, thanks to
(3.16) and the assumption $c\le e$, and so the
gamma function equals $\infty$. 

Hence, the sum of (3.18a), (3.18d), and (3.18e) vanishes.

\smallskip
Second, in the same way as before, it is seen that the sum of (3.18b)
and (3.18c) equals
$$\multline 
\sum _{i=j-b} ^{b-c-1
} \Bigg(
\sum_{k = 0}^{e- c-1 }{{\left( -1 \right) }^{b + i + s+1}}
{{ 
       \left( c + k - s -1 \right) !\, \left( e + k - s -1 \right) !}\over 
     {k!}}\\
\hskip3cm
\cdot
\frac {   (c+e-b+i+k-s+1)_{e-c-k-1}} 
{\left( 2 e - 2 s -2 \right) !}\Bigg)
\binom {b-c-e}{i-j+b}.
\endmultline\tag3.28$$
Similarly to before, 
the only change to be made in the arguments that showed that
the sum of (3.17b) and (3.17c) equals (3.23) is to start by replacing
the binomial $\binom {i+b-c-2e} {i-j+b}$ in (3.18b) by the 
expansion $\sum _{\ell=j-b}
^{i}\binom {b-c-e}{\ell-j+b} \binom {i-e}{i-\ell}$ (which is the
substitute of replacing the binomial 
$\binom {i-e}{i-j+c}$ in (3.17b) by some expansion), and in the
subsequent calculation replace the
binomial $\binom 0{\ell-j+c}$ by $\binom {b-c-e}{\ell-j+b}$ everywhere. 

Now,
of course, we cannot argue that the sum over $i$ in (3.28) consists
of just a single term, as opposed to (3.23), where we were could
derive the expression (3.24) accordingly. However, we may rewrite
(3.28) in the slightly fancier fashion
$$\multline 
\lim_{\de\to0}\Bigg(
\sum _{i=j-b} ^{b-c-1
} \Bigg(
\sum_{k = 0}^{e- c-1 }{{\left( -1 \right) }^{b + i + s+1}}
{{ 
       \left( c + k - s -1 \right) !\, \left( e + k - s -1 \right) !}\over 
     {k!}}\\
\hskip3cm
\cdot
\frac {   (\de+c+e-b+i+k-s+1)_{e-c-k-1}} 
{\left( 2 e - 2 s -2 \right) !}\Bigg)
\binom {b-c-e}{i-j+b}\Bigg),
\endmultline$$
interchange
the sums over $i$ and $k$, write the now inner sum over $i$
in hypergeometric notation,
$$\multline 
\lim_{\de\to0}\Bigg(
\sum _{k=0} ^{e-c-1}
{{\left( -1 \right) }^{j + s+1}} 
  {{  \left(  c + k - s-1 \right) !\,\left(  e + k - s-1 \right) !}\over 
    {k!}}\\
\cdot
\frac {    ( \de+c+e-2 b  + j + k - s+1)_{e-c-k-1}}
 {\left(  2 e - 2 s-2 \right) !}
   {} _{2} F _{1} \!\left [ \matrix { \de-2 b + 2 e + j - s, -b + c + e}\\ {
       \de - 2 b + c + e + j + k - s+1}\endmatrix ; {\displaystyle
1}\right ]\Bigg) ,
\endmultline$$
and sum the $_2F_1$-series, using the Chu--Vandermonde summation
(3.20) again, to get
$$\multline
\lim_{\de\to0}\Bigg(
\sum _{k=0} ^{e-c-1}
{{\left( -1 \right) }^{j + s+1}}{{ 
      \left(  c + k - s-1 \right) !\,\left(  e + k - s-1 \right) !}\over 
    {k!\,\left(  2 e - 2 s-2 \right) !}}\\
\cdot
\frac {   (1+c-e+k)_{b-c-e}\, 
      ( \de+c+e-2 b  + j + k - s+1)_{e-c-k-1}} 
{      ( \de+c+e-2 b  + j + k - s+1)_{b-c-e}}\Bigg).
\endmultline\tag3.29$$
Each summand in the sum over $k$ vanishes due to the occurence of the term
$$(1+c-e+k)_{b-c-e}=(1+c-e+k)(2+c-e+k)\cdots (b-2e+k)$$
in the numerator. For, excluding for the moment the case $e=c$, 
we have $1+c-e+k\le0$, because the summation
index $k$ is restricted above by $e-c-1$, and $b-2e+k\ge0$, 
because in the current case we are assuming $e\le b/2$.
The above argument does not apply when $e=c$, but in that case the sum
in (3.29) is empty, and so is zero anyway.

Hence, the sum of (3.18b) and (3.18c) vanishes, which, together with
our previous finding that the sum of (3.18a), (3.18d), and (3.18e)
vanishes, establishes the equation (3.18). 

Thus, the proof that $(x+e)$
divides $\De'(x;b,c)$ with multiplicity $m(e)$ as given in (3.12) is
complete.

\smallskip
{\it Case 3: $b/2\le e\le b-c$}.
By inspection of the expression (3.1),
we see that we have to prove that $(x+e)^{m(e)}$
divides $\De'(x;b,c)$, where
$$m(e)=\cases c+(b+c-2e)&b/2\le e<(b+c)/2\\
c&(b+c)/2\le e\le b-c.\endcases\tag3.30$$
Note that the second case in (3.30) could be empty, but not the first.

To extract the appropriate number of factors $(x+e)$ out of the
determinant $\De'(x;b,c)$, we start again with the modified determinant
(3.13). Again, 
we take $(x+e)$ out of rows $b-c,b-c+1,\dots,b-1$, and obtain the
determinant in (3.14), which we denoted by $\De_2(x;b,c,e)$. 
Obviously, we have taken out
$(x+e)^{c}$. 
Again, the remaining determinant has still entries which are
polynomial in $x$. This has to be argued here slightly differently
than it was for Case~2. Sure enough, the entries in rows
$i=0,1,\dots,b-c-1$ are polynomials in $x$. For $i\ge b-c$ we
have:  $c+i-j-1\ge b-j-1\ge0$ if $j<b$,
$b-c-e\ge0$ by assumption, and $c+e+i-j-1\ge
b+e-j-1\ge e-c\ge b/2-c\ge0$ if $j<b+c$.
This explains the term $c$ in (3.30).

Now let $e<(b+c)/2$.
In order to explain the term $(b+c-2e)$ in (3.30), we claim that for
$s=0,1,\dots,b+c-2e-1$ we have
$$\multline
\sum _{i=0} ^{2e-b+s} \Bigg(
\sum_{k = 0}^{2e-b-i+s }{{\left( -1 \right) }^{c + e + i + k + s+1}}
{{c - i-1 }\choose {b + c - 2 e + k - s-1 }}
\frac {\left( 2 b - c - 2 e - s -1 \right) !} 
{k!}
\\
\cdot{{  \left( b - e - s -1 \right) !\, \left( c + k - s -1 \right) !\, 
       \left( b - e + k - s -1 \right) !}\over 
     { \left( 2 b - 2 e - 2 s -2 \right) !\,  
   \left( 2 b - 2 e + k - 2 s -1 \right) !}}\Bigg)
\cdot \big(\text {row $i$ of $\De_2(-e;b,c,e)$}\big)\\
+
\sum _{i=2e+c-b} ^{e-1}
{{\left( -1 \right) }^{b + e + s}}
{{ \left( e - i -1 \right) !\, 
     \left( b - e - s -1 \right) !}\over 
   {\left( b - c - 2 e + i \right) !}}
\hskip5cm
\\
\cdot \frac {\left( b - 2 e + i - s -1 \right) !\, 
     ({ \textstyle b - i - s}) _{b - c - 2 e + i} } 
{  \left( 2 b - 2 e - 2 s -2 \right) !}
\cdot \big(\text {row $i$ of $\De_2(-e;b,c,e)$}\big)\\
+
\sum _{i=e} ^{b-c-1} \Bigg( {{\left( -1 \right) }^{b + c + i}}
{{ \left( b - c - i -1 \right) !\, 
      \left( c - s -1 \right) !\, \left( b - 2 e + i - s -1 \right) !\, 
      ({ \textstyle b - i - s}) _{i-e} }\over 
    {\left(i -e \right) !\, \left( 2 b - 2 e - 2 s -2 \right) !}} \\
 + 
\sum_{k = 0}^{b - c - e-1 }
{{\left( -1 \right) }^{b + i + s+1}}
     {{ 
         \left( b - 2 e + i - s -1 \right) !\, 
         \left( c + k - s -1 \right) !\, \left( b - e + k - s -1 \right) !}
        \over {k!\, \left( 2 b - 2 e - 2 s -2 \right) !\, 
         \left( c - e + i + k - s \right) !}}  \Bigg)\\
\cdot \big(\text {row $i$ of $\De_2(-e;b,c,e)$}\big)\\
+
\sum _{i=b-c} ^{b-s-1}
{{({ \textstyle 1 - b + c + i}) _{b - i - s-1 }  \,
     ({ \textstyle 1 - e + i}) _{b - i - s-1 } }\over 
   {\left( b - i - s -1 \right) !\, 
     ({ \textstyle b - 2 e + i - s}) _{b - i - s-1 } }}
\cdot \big(\text {row $i$ of $\De_2(-e;b,c,e)$}\big)\\
=0.\hskip8cm
\endmultline\tag3.31$$
Once more, note that these are indeed $b+c-2e$ linear combinations of the rows,
which are linearly independent.

Because of $s\le b+c-2e-1$, 
the rows which are involved in the first sum in (3.31) are
from rows $0,1,\dots,c-1$, which form the top block in (3.14). 
The assumption $e\ge b/2$ guarantees that the bounds for the sum are
proper bounds.
Because of the same assumption,
the rows which are involved in the second sum in
(3.31) are from rows $c,c+1,\dots,e-1$, which form the second block
from top in (3.14). The assumption $e\le b-c$
guarantees that the bounds for the sum are proper bounds,
(including the
possibility that $e=b-c$, in which case the sum is the
empty sum). The rows which are
involved in the third sum in (3.31) are clearly the rows
$e,e+1,\dots,b-c-1$, which form the third block from top in (3.14). The
bounds for the third sum are proper because of the
assumption $e\le b-c$ (including the
possibility that $e=b-c$, in which case the sum is the empty sum).
Finally, because of the condition $s\ge 0$, we have $b-s-1\le b-1$, and
therefore the rows which are involved in the fourth sum in (3.31) are
from rows $b-c,b-c+1,\dots,b-1$, which form the bottom block in (3.14).
The bounds for this fourth sum are proper because $s\le c-1$.
This inequality follows from the inequality chain
$$s\le b+c-2e-1\le b+c-b-1= c-1.\tag3.32$$
Again, it is also useful to observe that we need the restriction $e<(b+c)/2$
in order that there is at least one $s$ with $0\le s\le b+c-2e-1$.

Hence, in order to verify (3.31), we have to check
$$\align
&\sum _{i=0} ^{2e-b+s} \Bigg(
\sum_{k = 0}^{2e-b-i+s}{{\left( -1 \right) }^{c + e + i + k + s+1}}
{{c - i-1 }\choose {b + c - 2 e + k - s-1 }}
\frac {\left( 2 b - c - 2 e - s -1 \right) !} 
{k!}
\\
&\hskip1cm
\cdot{{  \left( b - e - s -1 \right) !\, \left( c + k - s -1 \right) !\, 
       \left( b - e + k - s -1 \right) !}\over 
     { \left( 2 b - 2 e - 2 s -2 \right) !\,  
   \left( 2 b - 2 e + k - 2 s -1 \right) !}}\Bigg)
\binom {c-e} {i-j+c}
\tag3.33a\\
&+
\sum _{i=2e+c-b} ^{e-1}
{{\left( -1 \right) }^{b + e + s}}
{{ \left( e - i -1 \right) !\, 
     \left( b - e - s -1 \right) !}\over 
   {\left( b - c - 2 e + i \right) !}}
\hskip5cm
\\
&\hskip2cm
\cdot \frac {\left( b - 2 e + i - s -1 \right) !\, 
     ({ \textstyle b - i - s}) _{b - c - 2 e + i} } 
{  \left( 2 b - 2 e - 2 s -2 \right) !}
\binom {i-e} {i-j+c}
\tag3.33b\\
& +\sum _{i=e} ^{b-c-1}\Bigg( 
\sum_{k = 0}^{b - c - e-1 }
{{\left( -1 \right) }^{b + i + s+1}}
     {{ 
         \left( b - 2 e + i - s -1 \right) !\, 
         \left( c + k - s -1 \right) !}
        \over {k!\, \left( 2 b - 2 e - 2 s -2 \right) !}}  \\
&\hskip4cm 
\cdot\frac {\left( b - e + k - s -1 \right) !} 
{         \left( c - e + i + k - s \right) !}\Bigg)
\binom {0} {i-j+c}
\tag3.33c\\
&+
\sum _{i=e} ^{b-c-1}  {{\left( -1 \right) }^{b + c + i}}
{{ \left( b - c - i -1 \right) !\, 
      \left( c - s -1 \right) !\, \left( b - 2 e + i - s -1 \right) !\, 
      ({ \textstyle b - i - s}) _{i-e} }\over 
    {\left(i -e \right) !\, \left( 2 b - 2 e - 2 s -2 \right) !}} \\
&\hskip8cm\cdot\binom {0} {i-j+c}
\tag3.33d\\
&+
\sum _{i=b-c} ^{b-s-1}
{{({ \textstyle 1 - b + c + i}) _{b - i - s-1 }  \,
     ({ \textstyle 1 - e + i}) _{b - i - s-1 } }\over 
   {\left( b - i - s -1 \right) !\, 
     ({ \textstyle b - 2 e + i - s}) _{b - i - s-1 } }}
(-1)^{c+i+j+1}\frac {1} {(i-j+c)}
\tag3.33e\\
&\hskip4cm=0,\tag3.33f
\endalign$$
which is (3.31) restricted to the $j$-th column, $j=c,c+1,\dots,b-1$, 
and
$$\align
&\sum _{i=0} ^{2e-b+s} \Bigg(
\sum_{k = 0}^{2e-b-i+s }{{\left( -1 \right) }^{c + e + i + k + s+1}}
{{c - i-1 }\choose {b + c - 2 e + k - s-1 }}
\frac {\left( 2 b - c - 2 e - s -1 \right) !} 
{k!}
\\
&\hskip1cm
\cdot{{  \left( b - e - s -1 \right) !\, \left( c + k - s -1 \right) !\, 
       \left( b - e + k - s -1 \right) !}\over 
     { \left( 2 b - 2 e - 2 s -2 \right) !\,  
   \left( 2 b - 2 e + k - 2 s -1 \right) !}}\Bigg)
2\binom {b-2e} {i-j+b}
\tag3.34a\\
&+
\sum _{i=b-c} ^{b-s-1}
{{({ \textstyle 1 - b + c + i}) _{b - i - s-1 }  \,
     ({ \textstyle 1 - e + i}) _{b - i - s-1 } }\over 
   {\left( b - i - s -1 \right) !\, 
     ({ \textstyle b - 2 e + i - s}) _{b - i - s-1 } }}\\
&\hskip2cm \cdot2(-1)^{e+c+i+j+1}\frac {(b-c-e)!\,(i-j+c+e-1)!} {(i-j+b)!}
\tag3.34b\\
&\hskip4cm=0,\tag3.34c
\endalign$$
which is (3.31) restricted to the $j$-th column,
$j=b,b+1,\dots,b+c-1$. Equation (3.34) is indeed the restriction of (3.31)
to the $j$-th column, $j=b,b+1,\dots,b+c-1$, because all the entries in
rows $2e+c-b,2e+c-b+1,\dots,b-c-1$ of $\De_2(-e;b,c,e)$ vanish
in such a column. For, due to the assumption $b/2\le e$, we have 
$$i+b-c-2e\le i-c<i-j+b,$$
therefore the 
entries $\binom {i+b-c-2e}
{i-j+b}$ in rows $2e+c-b,2e+c-b+1,\dots,e-1$ vanish in such a column, 
and for $i\ge e$ we have
$$b-c-e\le e-c\le i-c< i-j+b,$$
therefore the 
entries $\binom {b-c-e}
{i-j+b}$ in rows $e,e+1,\dots,b-c-1$ vanish in such a column.

We start by proving (3.33). We remind the reader that here $j$ is
restricted to $c\le j<b$.
Our strategy consists of exhibiting that (3.33) is equivalent to
(3.17) with $e$ replaced by $b-e$. Once this is done, the validity of
(3.33) follows immediately. (It should be observed that the ranges of
parameters in (3.33) and in (3.17) with $e$ replaced by $b-e$
correspond to each other perfectly: While in (3.33) the parameter
$e$ is restricted to $b/2\le e<(b+c)/2$, in (3.17) it is restricted to
$(b-c)/2<e\le b/2$, which matches nicely under the replacement $e\to
b-e$. A similar match occurs for the range of $s$.)

It is obvious that (3.33e) equals (3.17e) with $e$ replaced by $b-e$,
and that (3.33d) equals (3.17d) with $e$ replaced by $b-e$.
On other hand, it is not immediate that the sum of (3.33b) and (3.33c)
matches with the sum of (3.17b) and (3.17c), and
that (3.33a) matches (3.17a), under the same replacement.

First, we show how to convert (3.33a) into (3.17a) with $e$ replaced
by $b-e$. We replace $\binom {c-e} {i-j+c}$ in (3.33a) by the expansion $\sum
_{\ell=j-c} ^{i}\binom {c-b+e} {\ell-j+c} \binom {b-2e} {i-\ell}$, the
equality of binomial and expansion being again due to Chu--Vandermonde
summation. Then we interchange summations over $i$, $k$, $\ell$ 
so that the sum over
$\ell$ becomes the outer sum and the sum over $i$ becomes the inner
sum, and write the sum
over $i$ in hypergeometric notation. This gives for (3.33a) the
expression 
$$\multline
\sum _{\ell=j-c} ^{2e-b+s}\binom {c-b+e}{\ell-j+c}
\sum _{k=0} ^{2e-b-\ell+s}
{{\left( -1 \right) }^{b + e + \ell}} 
  {{  \left(   2 b - c - 2 e - s -1\right) !\,
      \left(   b - e - s -1\right) ! }\over 
    {k!\,\left(   2 b - 2 e - 2 s -2\right) !}}
\\
\cdot
\frac {    \left(   c + k - s -1\right) !\,  \left(   b - e + k - s -1\right) !\,
      ({ \textstyle 1 - c + \ell}) _{  b + c - 2 e + k - s-1}}
 {        \left(   2 b - 2 e + k - 2 s -1\right) !\,
      \left(   b + c - 2 e + k - s -1\right) !}\\
\cdot
   {} _{2} F _{1} \!\left [ \matrix {  2 e-b, b - 2 e + k + \ell - s}\\ {
       1 - c + \ell}\endmatrix ; {\displaystyle 1}\right ]  .
\endmultline$$
The $_2F_1$-series can be evaluated by the hypergeometric form (3.20)
of the Chu--Vandermonde summation. In the resulting expression we write
the sum over $k$ in hypergeometric notation, and obtain
$$\multline
\sum _{\ell=j-c} ^{2e-b+s}\binom {c-b+e}{\ell-j+c}
{{\left( -1 \right) }^{ c + e + \ell + s+1}} 
  {{  ( 2e+c- b - \ell-1)! \,
      (  2 b - c - 2 e - s-1)! }\over 
    {(  2 b - 2 e - 2 s-2)! \,
      (  2 b - 2 e - 2 s-1)! }}
\\
\cdot
\frac {      {(  b - e - s-1)! }^2}
     {  (2e-b  - \ell + s)! }
   {} _{2} F _{1} \!\left [ \matrix { b - e - s, b - 2 e + \ell - s}\\ { 2 b
       - 2 e - 2 s}\endmatrix ; {\displaystyle 1}\right ]  .
\endmultline$$
Again, Chu--Vandermonde summation (3.20) can be applied. After some
simplification, this gives the expression
$$\multline
\sum _{\ell=j-c} ^{2e-b+s}\binom {c-b+e}{\ell-j+c}
{{\left( -1 \right) }^{ c + e + \ell + s+1}} 
  {{    ( 2e+c- b  - \ell-1)! \,
      (  2 b - c - 2 e - s-1)! }\over 
    {(b-\ell-s-1)!\, (  2 b - 2 e - 2 s-2)! }}\\
\cdot
\frac {      {{(  b - e - s-1)! }}\,(e-\ell-1)!}
 {      (2e-b  - \ell + s)! },
\endmultline$$
which is exactly (3.17a) with $e$ replaced by $b-e$.

Now we turn to the relation of (3.33b)+(3.33c) and (3.17b)+(3.17c). 
In order to simplify matters, we make use of the fact that the sum of
(3.17b) and (3.17c) equals (3.23), as was shown in the
proof of (3.17).
It is apparent that (3.33c) matches with the part $i=e,e+1,\dots,b-c-1$ of
(3.23) with $e$ replaced by $b-e$. So, it remains to be seen that
(3.33b) matches with the remaining part $i=j-c,j-c+1,\dots,e-1$ of
(3.23), under the same replacement.

In order to verify the last assertion, we replace $\binom {i-e}{i-j+c}$ in
(3.33b) by the expansion $\sum _{\ell=j-c} ^{i}\binom
0{\ell-j+c}\binom {i-e}{i-\ell}$. That the binomial equals the
expansion is once again 
due to the Chu--Vandermonde summation. Subsequently, we interchange sums
over $i$ and $\ell$. In the now inner sum over $i$, we reverse the
order of summation, i.e., we replace $i$ by $e-1-i$, and then we write
the new sum over $i$ in
hypergeometric notation. This gives for (3.33b) the expression
$$\multline
\sum _{\ell=j-c} ^{e-1}\binom {0}{\ell-j+c}
{{\left( -1 \right) }^{b + \ell + s+1}} 
  {{  \left(   b - e - s -2\right) !\,\left(   b - e - s -1\right) !}\over 
   {\left(   b - c - e -1\right) !}
   }\\
\cdot
\frac {     ({ \textstyle 1 + b - e - s}) _{  b - c - e-1} }
 {\left(   2 b - 2 e - 2 s -2\right) !}
  {} _{3} F _{2} \!\left [ \matrix { 1 - e + \ell, 1, 1 - b + c + e}\\ { 2 - b
      + e + s, 1 + b - e - s}\endmatrix ; {\displaystyle 1}\right ]  .
\endmultline\tag3.35$$
To the $_3F_2$-series we apply, once again, the transformation formula
(3.21). We need to apply the case where $n=b-c-e-1$. Due to our assumption
$e\le b-c$, this is indeed a
nonnegative integer, except if $e=b-c$. So, let us for the moment
exclude the case $e=b-c$. After little manipulation, application of
(3.21) to (3.35) yields the following expression for (3.33b):
$$\multline 
\sum _{i=j-c} ^{e-1} \Bigg(
\sum_{k = 0}^{b-e- c-1 }{{\left( -1 \right) }^{b + i + s+1}}
{{ 
       \left( c + k - s -1 \right) !\, \left( b-e + k - s -1 \right) !}\over 
     {k!}}\\
\hskip3cm
\cdot
\frac {     (c-e+i+k-s+1)_{b-c-e-k-1}} 
{\left( 2 b-2e - 2 s -2 \right) !}\Bigg)
\binom {0}{i-j+c}.
\endmultline$$
This is exactly the part $i=j-c,j-c+1,\dots,e-1$ of (3.23) with $e$ replaced
by $b-e$. In the excluded case $e=b-c$, 
(3.33b) and the part $i=j-c,j-c+1,\dots,e-1$ of (3.23) with $e$ replaced
by $b-e$ 
are also in agreement, since both expressions
vanish in that case, due to an empty summation over $i$ in (3.33b) and
an empty summation over $k$ in (3.23) with $e$ replaced
by $b-e$. Hence, in all cases, we have
established the equality of (3.33b) and the part $i=j-c,j-c+1,\dots,e-1$ of
(3.23) with $e$ replaced
by $b-e$. Therefore, in all cases, the sum of (3.33b) and (3.33c) is
equal to the sum of (3.17b) and (3.17c) with $e$ replaced
by $b-e$. 

This completes the argument that (3.33) is equivalent to (3.17) 
with $e$ replaced by $b-e$, and so establishes the equation (3.33).

\smallskip
Now we turn to (3.34). We claim that the two sums in (3.34) can be
combined into one term,
$$\multline
\lim_{\de\to0}\Bigg(
\sum_{i = j-b}^{b-s-1}2\,{{\left( -1 \right) }^{c + e + i + j+1}} 
{{       \left( b - c - e \right) !\,
       ({ \textstyle 1 + {\de}}) _{  c + e + i - j-1} }\over 
     {\left( b + i - j \right) !\,\left(   b - i - s -1\right) !}}
\\
\cdot
\frac {       ({ \textstyle 1 - b + c + {\de} + i}) _{  b - i - s-1} \,
       ({ \textstyle 1 - e + {\de} + i}) _{  b - i - s-1} }
 {       ({ \textstyle b - 2 e + {\de} + i - s}) _{  b - i - s-1} }
\Bigg).
\endmultline\tag3.36$$
It is obvious that (3.34b) agrees with the according part
$i=b-c,b-c+1,\dots,b-s-1$ of (3.36). It is also straight-forward to check
that the terms for $i=2e-b+s+1,2e-b+s+2,\dots,b-c-1$ in (3.36) vanish.
It remains to be seen that (3.34a)
agrees with the according part $i=j-b,j-b+1,\dots,2e-b+s$ of (3.36).

In order to verify the last assertion, we replace $\binom {b-2e}
{i-j+b}$ in (3.34a) by the expansion $\sum _{\ell=j-b} ^{i}\binom
{b-c-e}{\ell-j+c}\binom {c-e}{i-\ell}$, again making use of the
Chu--Vandermonde summation. Then we interchange summations over $i$, $k$, $\ell$ 
so that the sum over
$\ell$ becomes the outer sum and the sum over $i$ becomes the inner
sum, and write the sum
over $i$ in hypergeometric notation. This gives for (3.34a) the
expression 
$$\multline
\sum _{\ell=j-b} ^{2e-b+s}\binom {b-c-e}{\ell-j+b}
\sum _{k=0} ^{2e-b-\ell+s}
2\,{{\left( -1 \right) }^{b + e + \ell}} 
     {{       (  2 b - c - 2 e - s-1)! \,
      (  b - e - s-1)!  }\over 
    {k! \,(  2 b - 2 e - 2 s-2)!  }}
\\
\cdot
\frac {         (  c + k - s-1)!\,   (  b - e + k - s-1)! \,
      ({ \textstyle 1 - c + \ell}) _{  b + c - 2 e + k - s-1} }
 {    (  2 b - 2 e + k - 2 s-1)!\,   (  b + c - 2 e + k - s-1)! }\\
\cdot
{} _{2} F _{1} \!\left [ \matrix { e-c, b - 2 e + k + \ell - s}\\ { 1 -
       c + \ell}\endmatrix ; {\displaystyle 1}\right ]  .
\endmultline$$
We sum the $_2F_1$-series by means of the hypergeometric form (3.20)
of the Chu--Vandermonde summation, and in the resulting expression write
the inner sum over $k$ in hypergeometric notation, to obtain the expression
$$\multline
\sum _{\ell=j-b} ^{2e-b+s}\binom {b-c-e}{\ell-j+b}
2\,{{\left( -1 \right) }^{ c + e + \ell + s+1}} 
     {{       (  e - \ell-1)! \,(  c - s-1)! }\over 
    {(  2 b - 2 e - 2 s-2)! }}
\\
\cdot
\frac {      (  2 b - c - 2 e - s-1)! \,
      (  b - e - s-1)! }
 {      (  2 b - 2 e - 2 s-1)! \,
      (2e-b - \ell + s)! }
{} _{2} F _{1} \!\left [ \matrix { c - s, b - 2 e + \ell - s}\\ { 2 b -
       2 e - 2 s}\endmatrix ; {\displaystyle 1}\right ]  .
\endmultline$$
Another application of the Chu--Vandermonde summation (3.20) yields the
expression
$$\multline
\sum _{\ell=j-b} ^{2e-b+s}\binom {b-c-e}{\ell-j+b}
2\,{{\left( -1 \right) }^{  c + e + \ell + s+1}} 
     {{ (  b - c - \ell-1)! \,(  e - \ell-1)!  }
     \over {(  2 b - 2 e - 2 s-2)! }}\\
\cdot
\frac {      (  c - s-1)! \,(  b - e - s-1)!}
 {      (  b - \ell - s-1)! \,
      (2e-b  - \ell + s)! }
\endmultline$$
for (3.34a). It is now straight-forward to check that this sum
is equal the according part
$i=j-b,j-b+1,\dots,2e-b+s$ of (3.36).

So, in order to prove (3.34), we need to show that (3.36) vanishes. To
accomplish this, we write the sum in (3.36) in hypergeometric
notation,
$$\multline
\lim_{\de\to0}\Bigg(
2\,{{\left( -1 \right) }^{b + c + e +1}} 
  {{      ({ \textstyle 1}) _{b - c - e} \,
     ({ \textstyle 1 + {\de}}) _{  b + c + e-1} }\over 
   {({ \textstyle 1}) _{  2 b - j - s-1}  }}
\\
\times
\frac {     ({ \textstyle 1 - 2 b + c + {\de} + j}) _{  2 b - j - s-1} \,
     ({ \textstyle 1 - b - e + {\de} + j}) _{  2 b - j - s-1} }
 {     ({ \textstyle -2 e + {\de} + j - s}) _{  2 b - j - s-1}}
\\
\times
  {} _{3} F _{2} \!\left [ \matrix { -b + c + e +
{\de}, -2 e + {\de} + j - s,  1 - 2 b + j + s}\\ 
{ 1 - 2 b + c + {\de} + j, 1 - b -
      e + {\de} + j}\endmatrix ; {\displaystyle 1}\right ]  
\Bigg),
\endmultline$$
and to the $_3F_2$-series apply the Saalsch\"utz summation (3.8). We
have to apply the case where $n=2b-j-s-1$, which is indeed a
nonnegative integer because of the inequality chain
$$2b-j-s-1\ge b-c-s\ge c-s\ge 1,$$
the last inequality being due to (3.32).
Thus, we obtain for (3.36) the expression
$$\multline
\lim_{\de\to0}\Bigg(
2\,{{\left( -1 \right) }^{b + c + e + 1}} 
  {{   ({ \textstyle 1}) _{b - c - e} \,
     ({ \textstyle 1 + {\de}}) _{ c+e- b -1} \,
     ({ \textstyle 1 - b - e + j}) _{  2 b - j - s-1} }\over 
   {({ \textstyle 1}) _{  2 b - j - s-1} \,
     ({ \textstyle -2 e + {\de} + j - s}) _{  2 b - j - s-1} }}\\
\cdot
\frac {     ({ \textstyle 1 - b - e + {\de} + j}) _{  2 b - j - s-1} \,
     ({ \textstyle 1 - 2 b + c + 2 e + s}) _{  2 b - j - s-1} }
 {     ({ \textstyle 1 - b + e - {\de} + s}) _{  2 b - j - s-1} }
\Bigg).
\endmultline$$
This expression is indeed zero, because of the occurence of the term 
$$(1-b-e+j)_{2b-j-s-1}=(1-b-e+j)(2-b-e+j)\cdots (b-e-s-1)$$
in the numerator.
For, we have $1-b-e+j\le -e+c\le 0$, and we have $b-e-s-1\ge -c+e\ge
0$.
This establishes equation (3.34).

Thus, the proof that $(x+e)$
divides $\De'(x;b,c)$ with multiplicity $m(e)$ as given in (3.30) is
complete.

\smallskip
{\it Case 4: $b-c\le e\le b-c/2$}.
By inspection of the expression (3.1),
we see that we have to prove that $(x+e)^{m(e)}$
divides $\De'(x;b,c)$, where
$$m(e)=\cases (2b-2e-c)+(b+c-2e)&b-c\le e< (b+c)/2\\
(2b-2e-c)&(b+c)/2\le e\le b-c/2.\endcases
\tag3.37$$
Note that the first case in (3.37) could be empty, but not the second,
because of $b\ge 2c$.

In order to explain the term $2b-2e-c$ in (3.37), we start with the
determinant (3.13) with $e=b$. (This determinant equals 
$\De'(x;b,c)$ as we showed by a few row manipulations at the beginning
of Case~2.) The choice of $e=b$ has the effect that the bottom block
in (3.13) is empty. Now we take $(x+e)$ out of rows
$2e+c-b,2e+c-b+1,\dots,b-1$, and obtain the determinant
$$\multline
\hskip-.5cm
\det_{0\le i<b,\, c\le j\le b+c}\!\!\(
\PfadDicke{.3pt}
\SPfad(0,2),1111111111111111\endSPfad
\SPfad(0,0),1111111111111111\endSPfad
\SPfad(8,-3),2222222\endSPfad
\Label\o{\raise25pt\hbox{$\dsize\binom {x+c}{i-j+c}$}}(4,2)
\Label\o{\raise25pt\hbox{$\dsize2\binom {2x+b}{i-j+b}$}}(12,2)
\Label\o{\raise25pt\hbox{$\dsize\binom {x+i}{i-j+c}$}}(4,0)
\Label\o{\raise25pt\hbox{$\dsize\binom {2x+b-c+i}{i-j+b}$}}(12,0)
\Label\o{\raise15pt\hbox{\hskip-10pt$\dsize 
\frac {(x+e+1)_{i-e}} {(i-j+c)!}$}}(3,-2)
\Label\o{\raise25pt\hbox{\hskip0pt$\dsize
\times(x-c+j+1)_{e+c-j-1}$}}(4,-4)
\Label\o{\raise15pt\hbox{\hskip23pt$\dsize 
2\frac {(2x+2e+1)_{i+b-c-2e}} {(i-j+b)!}$}}(11,-2)
\Label\o{\raise25pt\hbox{\hskip5pt$\dsize
\times(2x-c+j+1)_{2e+c-j-1}$}}(12,-4)
\Label\o{\raise15pt\hbox{$0\!\le\! i\!<\!c$\hskip-5pt}}(19,2)
\Label\o{\raise15pt\hbox{$c\!\le\! i\!<\!2e+c-b$\hskip-5pt}}(19,0)
\Label\o{\raise20pt\hbox{$2e+c-b\!\le\! i\!<\!b$\hskip-5pt}}(19,-3)
\Label\o{,}(17,-4)
\hskip8cm
\)
\hbox{\hskip-8.3cm}
\Label\o{\raise0pt\hbox{$c\le j<b$}}(4,-5)
\Label\o{\raise0pt\hbox{$b\le j<b+c$}}(12,-5)
\hskip10cm\\
\endmultline
\tag3.38$$
\vskip10pt
\noindent
which we denote by $\De_3(x;b,c,e)$. Obviously, we have taken out
$(x+e)^{2b-2e-c}$. The remaining determinant has still entries which are
polynomial in $x$. For, it is obvious that the entries in rows
$i=0,1,\dots,2e+c-b-1$ are polynomials in $x$, and for $i\ge 2e+c-b$ we
have: $i-e\ge e+c-b\ge0$ by assumption, 
 $e+c-j-1\ge e+c-b\ge 0$ if
$j<b$, $i+b-c-2e\ge0$, and $2e+c-j-1\ge
2e-b\ge b-2c\ge0$ if $j<b+c$. This explains the term $(2b-2e-c)$ in (3.37).

Now let $e<(b+c)/2$. In order to explain the term $(b+c-2e)$ in
(3.37), we claim that for $s=0,1,\dots,b+c-2e-1$ we have
$$\multline
\sum _{i=0} ^{2e-b+s} \Bigg(
\sum_{k = 0}^{2e-b-i+s }{{\left( -1 \right) }^{c + e + i + k + s+1}}
       {{c - i-1 }\choose {b + c - 2 e + k - s-1 }} 
\frac {\left( 2 b - c - 2 e - s -1 \right) !} 
{k!}
\\
\cdot{{       \left( b - e - s -1 \right) !\, \left( c + k - s -1 \right) !\, 
       \left( b - e + k - s -1 \right) !}\over 
     { \left( 2 b - 2 e - 2 s -2 \right) !\,
      \left( 2 b - 2 e + k - 2 s -1 \right) !}}\Bigg)
\cdot \big(\text {row $i$ of $\De_3(-e;b,c,e)$}\big)\\
+
\sum _{i=e} ^{2e+c-b-1}
2 {{\left( -1 \right) }^{c + s}}
{{ \left( 2e + c- b - i -1 \right) !\, 
     \left( 2 b - c - 2 e - s -1 \right) !\, 
     ({ \textstyle 1 - e + i}) _{b - i - s-1 } }\over 
   {\left( b - i - s -1 \right) !\, 
     ({ \textstyle b - 2 e + i - s}) _{b - i - s-1 } }}\\
\cdot \big(\text {row $i$ of $\De_3(-e;b,c,e)$}\big)\\
+
\sum _{i=2e+c-b} ^{b-s-1}
{{\left( -1 \right) }^{b - i - s-1 }}
{{ 
     ({ \textstyle 1 + b - c - 2 e + i}) _{b - i - s-1 }  \,
     ({ \textstyle 1 - e + i}) _{b - i - s-1 } }\over 
   {\left( b - i - s -1 \right) !\, 
     ({ \textstyle b - 2 e + i - s}) _{b - i - s-1 } }}\\
\cdot \big(\text {row $i$ of $\De_3(-e;b,c,e)$}\big)\\
=0.\hskip8cm
\endmultline\tag3.39$$
Once again, note that these are indeed $b+c-2e$ linear combinations of the rows,
which are linearly independent. 

Because of $s\le b+c-2e-1$, 
the rows which are involved in the first sum in (3.39) are
from rows $0,1,\dots,c-1$, which form the top block in (3.38). 
The inequality chain $2e-b+s\ge b-2c+s\ge 0$ 
guarantees that the bounds for the sum are proper bounds.
Because of $e\ge b-c\ge c$,
the rows which are involved in the second sum in
(3.39) are from rows $c,c+1,\dots,2e+c-b-1$, which form the middle block
in (3.38). The assumption $b-c\le e$
guarantees that the bounds for the sum are proper bounds
(including the
possibility that $e=b-c$, in which case the sum is the
empty sum). 
Finally, because of the condition $s\ge 0$, we have $b-s-1\le b-1$, and
therefore the rows which are involved in the third sum in (3.39) are
from rows $2e+c-b,2e+c-b+1,\dots,b-1$, which form the bottom block in (3.38).
The bounds for this third sum are proper because of $2e+c-b\le 2e-c\le
b-s-1$. 
Again, it is also useful to observe that we need the restriction $e<(b+c)/2$
in order that there is at least one $s$ with $0\le s\le b+c-2e-1$.

Hence, in order to verify (3.39), we have to check
$$\align
&\sum _{i=0} ^{2e-b+s} \Bigg(
\sum_{k = 0}^{2e-b-i+s }{{\left( -1 \right) }^{c + e + i + k + s+1}}
       {{c - i-1 }\choose {b + c - 2 e + k - s-1 }} 
\frac {\left( 2 b - c - 2 e - s -1 \right) !} 
{k!}
\\
&\hskip1cm
\cdot{{       \left( b - e - s -1 \right) !\, \left( c + k - s -1 \right) !\, 
       \left( b - e + k - s -1 \right) !}\over 
     { \left( 2 b - 2 e - 2 s -2 \right) !\,
      \left( 2 b - 2 e + k - 2 s -1 \right) !}}\Bigg)
\binom {c-e}{i-j+c}
\tag3.40a\\
&+
\sum _{i=2e+c-b} ^{b-s-1}
{{\left( -1 \right) }^{b +c+e+ i+j + s }}
{{ 
     ({ \textstyle 1 + b - c - 2 e + i}) _{b - i - s-1 }  \,
     ({ \textstyle 1 - e + i}) _{b - i - s-1 } }\over 
   {\left( b - i - s -1 \right) !\, 
     ({ \textstyle b - 2 e + i - s}) _{b - i - s-1 } }}\\
&\hskip2cm
\cdot \frac {(i-e)!\,(e+c-j-1)!} {(i-j+c)!}
\tag3.40b\\
&\hskip4cm=0,\tag3.40c
\endalign$$
which is (3.39) restricted to the $j$-th column, $j=c,c+1,\dots,b-1$,
(note that this is indeed the restriction of (3.39) to the $j$-th
column, $c\le j<b$, since, due to $0\le i-e\le
i-b+c<i-j+c$, the entries $\binom {i-e} {i-j+c}$ in rows
$e,e+1,\dots,2e+c-b-1$ of $\De_3(-e;b,c,e)$ vanish in such a column),
and
$$\align
&\sum _{i=0} ^{2e-b+s} \Bigg(
\sum_{k = 0}^{2e-b-i+s }{{\left( -1 \right) }^{c + e + i + k + s+1}}
       {{c - i-1 }\choose {b + c - 2 e + k - s-1 }} 
\frac {\left( 2 b - c - 2 e - s -1 \right) !} 
{k!}
\\
&\hskip1cm
\cdot{{       \left( b - e - s -1 \right) !\, \left( c + k - s -1 \right) !\, 
       \left( b - e + k - s -1 \right) !}\over 
     { \left( 2 b - 2 e - 2 s -2 \right) !\,
      \left( 2 b - 2 e + k - 2 s -1 \right) !}}\Bigg)
2\binom {b-2e}{i-j+b}
\tag3.41a\\
&+
\sum _{i=e} ^{2e+c-b-1}
2 {{\left( -1 \right) }^{c + s}}
{{ \left( 2e + c- b - i -1 \right) !\, 
     \left( 2 b - c - 2 e - s -1 \right) !\, 
     ({ \textstyle 1 - e + i}) _{b - i - s-1 } }\over 
   {\left( b - i - s -1 \right) !\, 
     ({ \textstyle b - 2 e + i - s}) _{b - i - s-1 } }}\\
&\hskip5cm\cdot \binom {i+b-c-2e}{i-j+b}
\tag3.41b\\
&+
\sum _{i=2e+c-b} ^{b-s-1}
{{\left( -1 \right) }^{b+c + i+j + s }}
{{ 
     ({ \textstyle 1 + b - c - 2 e + i}) _{b - i - s-1 }  \,
     ({ \textstyle 1 - e + i}) _{b - i - s-1 } }\over 
   {\left( b - i - s -1 \right) !\, 
     ({ \textstyle b - 2 e + i - s}) _{b - i - s-1 } }}\\
&\hskip2cm
\cdot 2\,\frac {(i+b-c-2e)!\,(2e+c-j-1)!} {(i-j+b)!}
\tag3.41c\\
&\hskip4cm=0,\tag3.41d
\endalign$$
which is (3.39) restricted to the $j$-th column, $j=b,b+1,\dots,b+c-1$.

We start by proving (3.40). We remind the reader that here $j$ is
restricted to $c\le j<b$.
We claim that the left-hand side of (3.40) can be written as
$$\multline
\lim_{\de\to0}\Bigg(
\sum_{i = j-c}^{  b - s-1}{{\left( -1 \right) }^
         {b + c + e + i + j + s}}
{{ \left(   b - e - s -1\right) !\, 
       ({ \textstyle 1 + { \de}}) _{  c + e - j-1} }
      \over {\left( i - j +c\right) !\, \left(   b - i - s -1\right) ! }}
\\
\cdot
\frac {       ({ \textstyle 1 + b - c - 2 e + { \de} + i}) _{  b - i - s-1}}
 {       ({ \textstyle b - 2 e + { \de} + i - s}) _{  b - i - s-1}}
\Bigg).
\endmultline\tag3.42$$
It is apparent that (3.40b) equals the according part
$i=2e+c-b,2e+c-b+1,\dots,b-s-1$ of (3.42). It is also easy to check
that the terms for $i=2e-b+s+1,2e-b+s+2,\dots,2e+c-b-1$ in (3.42)
vanish. It remains to be seen that
(3.40a) equals the remaining part $i=j-c,j-c+1,\dots,2e-b+s$ of
(3.42), which is not directly evident.

In order to verify this last assertion, we replace $\binom {c-e}
{i-j+c}$ in (3.40a) by the expansion $\sum _{\ell=j-c} ^{i}\binom
{\ell-e} {\ell -j+c} \binom {c-\ell-1} {i-\ell}$, making again use of
the Chu--Vandermonde summation. Then we interchange summations over
$i$, $k$, $\ell$ so that the sum over
$\ell$ becomes the outer sum and the sum over $i$ becomes the inner
sum, and write the sum
over $i$ in hypergeometric notation. This gives for (3.40a) the
expression 
$$\multline
\sum _{\ell=j-c} ^{2e-b+s}
{{  \ell -e}\choose {\ell - j + c }} 
\sum _{k=0} ^{2e-b-\ell+s}
\kern-2pt {{\left( -1 \right) }^{c + e + k + \ell  + s+1}}
\frac { (  2 b - c - 2 e - s-1)!\,  
      (  b - e - s-1)!   } 
{k!}\\
\cdot
{{    (  c + k - s-1)!     (  b - e + k - s-1)!\,  
      ({ \textstyle 1 - b + 2 e - k - \ell  + s}) _{  b + c - 2 e + k - s-1} }
     \over { (  2 b - 2 e - 2 s-2)!  
      (  2 b - 2 e + k - 2 s-1)!\,  
      (  b + c - 2 e + k - s-1)! }}\\
\cdot
      {} _{1} F _{0} \!\left [ \matrix { b - 2 e + k + \ell  - s}\\ {
       -}\endmatrix ; {\displaystyle 1}\right ]  .
\endmultline\tag3.43$$
Clearly, because of the hypergeometric form of the 
binomial theorem (see \cite{\SlatAC, Appendix (III.1)}),
$$
{} _{1} F _{0} \!\left [ \matrix { a}\\ { -}\endmatrix ; {\displaystyle
   z}\right ]  = {{\left( 1 - z \right) }^{-a}},
$$
the $_1F_0$-series in (3.43) is nonzero only if $k=2e-b-\ell+s$, in
which case it is $1$. Hence, we obtain for (3.40a) the expression
$$\multline
\sum _{\ell=j-c} ^{2e-b+s}
{{\left( -1 \right) }^{b + c + e+1}} 
      {{  \ell -e}\choose {\ell - j + c }} \\
\times{{(  e - \ell -1)!  \,
      ( 2e+c- b  - \ell -1)!\,  
      (  2 b - c - 2 e - s-1)!  \,
      (  b - e - s-1)! }\over 
    {(  2 b - 2 e - 2 s-2)!  \,
      (  b - \ell  - s-1)!  \,
    (2e-b  - \ell  + s)! }}.
\endmultline$$
It is now readily checked that this agrees with the part
$i=j-c,j-c+1,\dots,2e-b+s$ of (3.42).

Hence, in order to prove (3.40), we have to show that (3.42) vanishes. 
We do this by writing the sum (3.42) in hypergeometric notation,
$$\multline
\lim_{\de\to0}\Bigg(
{{\left( -1 \right) }^{b + e + s}} 
  {{     (  b - e - s-1)!\,  
     ({ \textstyle 1 + { \de}}) _{  c + e - j-1}  \,
     ({ \textstyle 1 + b - 2 c - 2 e + { \de} + j}) _{  b + c - j - s-1} 
     }\over {(  b + c - j - s-1)!\,  
     ({ \textstyle b - c - 2 e + { \de} + j - s}) _{  b + c - j - s-1} }}
\\
\times 
{} _{2} F _{1} \!\left [ \matrix {  b - c - 2 e + {
      \de} + j - s, 1 - b - c + j + s}\\ 
 { 1 + b - 2 c - 2 e + { \de} + j}\endmatrix ;
      {\displaystyle 1}\right ]
\Bigg),
\endmultline$$
and summing the $_2F_1$-series, once again, by means of the
hypergeometric form (3.20) of Chu--Vandermonde summation. We have to
apply the case where $n=b+c-j-s-1$. This is indeed a nonnegative
integer, because of $j\le b-1$ and because of the inequality chain
$$s\le b+c-2e-1\le b+c-2(b-c)-1=3c-b-1\le c-1.\tag3.44$$
Thus, we obtain for (3.42) the expression
$$
\lim_{\de\to0}\Bigg(
{{\left( -1 \right) }^{b + e + s}} 
  {{   (  b - e - s-1)!\,  
     ({ \textstyle 1 + { \de}}) _{  c + e - j-1}\,  
     ({ \textstyle 1 - c + s}) _{  b + c - j - s-1} }\over 
   {(  b + c - j - s-1)!\,  
     ({ \textstyle b - c - 2 e + { \de} + j - s}) _{  b + c - j - s-1} }}
\Bigg),
$$
which does indeed vanish due to the occurence of the term
$$   ({ \textstyle 1 - c + s}) _{  b + c - j - s-1} =
(1-c+s)(2-c+s)\cdots (b-j-1)$$
in the numerator.
For, by (3.44) we have $1-c+s\le0$ and we have $b-j-1\ge0$.
This establishes the equation (3.40).

\smallskip
Finally, we turn to (3.41). We remind the reader that here $j$ is
restricted to $b\le j<b+c$.
We claim that the left-hand side of
(3.41) can be combined into the single term
$$\multline
\lim_{\de\to0}\Bigg(
\sum_{i = j-b}^{b-s-1}2 {{\left( -1 \right) }^{b + c + i + j + s}} 
\frac {  \left(   c + 2 e - j -1\right) !\, 
       ({ \textstyle 1 + { \de}}) _{b - c - 2 e + i}}
 {\left( i - j +b\right) !\, \left(   b - i - s -1\right) !}\\
\cdot
{{     ({ \textstyle 1 + b - c - 2 e + { \de} + i}) _{  b - i - s-1}\,  
       ({ \textstyle 1 - e + { \de} + i}) _{  b - i - s-1} }\over 
     { ({ \textstyle b - 2 e + { \de} + i - s}) _{  b - i - s-1} }}
\Bigg).
\endmultline\tag3.45$$
In fact, it is straight-forward to check that (3.41c) agrees with the
according part $i=2e+c-b,2e+c-b+1,\dots,b-s-1$ of (3.45), and that
(3.41b) agrees with the according part $i=e,e+1,\dots,2e+c-b-1$ of
(3.45). It is also easy to see that the terms for
$i=2e-b+s+1,2e-b+s+2,\dots,e-1$ in (3.45) vanish. 
That (3.41a) agrees with the part $i=j-b,j-b+1,\dots,2e-b+s$ of (3.45)
is proved in the same way as it was proved before that (3.40a) agrees
with the part $i=j-c,j-c+1,\dots,2e-b+s$ of (3.42). The only change to
be made is to start by replacing the binomial $\binom {b-2e} {i-j+b}$
in (3.41a) by the
expansion $\sum _{\ell=j-b} ^{i}\binom {\ell+b-c-2e} {\ell-j+b} \binom
{c-\ell-1} {i-\ell}$ 
(which is the substitute 
of replacing the binomial $\binom {c-e} {i-j+c}$ in (3.40a) by an
expansion), and in the subsequent calculation replace the binomial
$\binom {\ell-e}{\ell-j+c}$ by $\binom {\ell +b-c-2e}{\ell-j+b}$
everywhere. 

So, in order to prove equation (3.41), we need to show that (3.45)
vanishes. In hypergeometric notation, the expression (3.45) reads
$$\multline
\lim_{\de\to0}\Bigg(
2 {{\left( -1 \right) }^{c + s}} 
  {{     (  c + 2 e - j-1)!\,  
     ({ \textstyle 1 + { \de}}) _{-c - 2 e + j}\,  
     ({ \textstyle 1 - c - 2 e + { \de} + j}) _{  2 b - j - s-1}}\over 
   {(  2 b - j - s-1)!}}
\\
\times
\frac {     ({ \textstyle 1 - b - e + { \de} + j}) _{  2 b - j - s-1} } 
{({ \textstyle -2 e + { \de} + j - s}) _{  2 b - j - s-1} }
     {} _{2} F _{1} \!\left [ \matrix {  -2 e + { \de} + j
      - s, 1 - 2 b + j + s}\\ { 1 - b - e + { \de} + j}\endmatrix ;
{\displaystyle 1}\right ]
\Bigg).
\endmultline$$
Clearly, we want to apply the hypergeometric form (3.20) of
Chu--Vandermonde summation again, with $n=2b-j-s-1$. This is indeed a
nonnegative integer because of the inequality chain
$$2b-j-s-1\ge b-c-s\ge 2e-2c+1\ge 2b-4c+1\ge 1. $$
Thus, we obtain for (3.45) the expression
$$\multline
\lim_{\de\to0}\Bigg(
2 {{\left( -1 \right) }^{c + s}} 
\frac {(  c + 2 e - j-1)!\,  
     ({ \textstyle 1 + { \de}}) _{-c - 2 e + j}}
 {(  2 b - j - s-1)!}\\
\times
  {{   ({ \textstyle 1 - c - 2 e + { \de} + j}) _{  2 b - j - s-1}\,  
     ({ \textstyle 1 - b + e + s}) _{  2 b - j - s-1} }\over 
   { ({ \textstyle -2 e + { \de} + j - s}) _{  2 b - j - s-1} }}
\Bigg).
\endmultline$$
This expression does indeed vanish, because of the occurence of the
term
$$ ({ \textstyle 1 - b + e + s}) _{  2 b - j - s-1} =
(1-b+e+s)(2-b+e+s)\cdots (b+e-j-1)$$
in the numerator. For, we have $1-b+e+s\le c-e\le 2c-b\le0$ and
$b+e-j-1\ge e-c\ge0$. This establishes equation (3.41).

Thus, the proof that $(x+e)$
divides $\De'(x;b,c)$ with multiplicity $m(e)$ as given in (3.37) is
complete.

\medskip
This finishes the proof of Lemma~1.\quad \quad \qed

\enddemo

\proclaim{Lemma 2}Let $b$ and $c$ be nonnegative
integers such that $c\le b\le 2c$.
Then the product
$$\prod _{i=1} ^{c}{\(x+\cl{\frac {c+i}
{2}}\)_{b-c+\cl{i/2}-\cl{(c+i)/2}} \, \(x+\cl{\frac {b-c+i}
{2}}\)_{\cl{(b+i)/2}-\cl{(b-c+i)/2}}}\tag3.46$$
divides $\De'(x;b,c)$, the determinant given by {\rm (2.3)}, 
as a polynomial in $x$.
\endproclaim
\demo{Proof} As in the proof of Lemma~1,
let us concentrate on some factor $(x+e)$ which appears
in (3.46), say with multiplicity $m(e)$. We have to prove that
$(x+e)^{m(e)}$ divides $\De'(x;b,c)$. As before, 
we accomplish this by finding
$m(e)$ linear combinations of the rows of $\De'(x;b,c)$ (or of an
equivalent determinant) that vanish
for $x=-e$, and which are linearly independent. 

Also here, we have to distinguish between four cases, depending on the
magnitude of $e$. The first case is $(b-c)/2\le e\le b-c$, the second case is
$b-c\le e\le b/2$, the third case is $b/2\le e\le c$, and the fourth
case is $c\le e\le (b+c)/2$. 

\smallskip
{\it Case 1: $(b-c)/2\le e\le b-c$}. By inspection of the expression
(3.46),
we see that we have to prove that $(x+e)^{m(e)}$
divides $\De'(x;b,c)$, where
$$m(e)=\cases (2e+c-b)&(b-c)/2\le e\le c/2\\
(2e+c-b)+(2e-c)&c/2< e\le b-c.\endcases
\tag3.47$$
Note that the second case in (3.47) could be empty, but not the first,
because of $b\le 2c$.

The term $(2e-c)$ in the ($c/2<e$)-case of (3.47) is easily explained: 
As in Case~1 of the proof of Lemma~1, 
we take $(x+e)$ out of rows $b+c-2e,b+c-2e+1,\dots,b-1$ of the
determinant $\De'(x;b,c)$ (clearly, such rows exist only if $c/2<e$),
and thus obtain the determinant (3.3), which we denoted by $\De_1(x;b,c,e)$.
Also here, this determinant has still entries which are
polynomial in $x$. For, it is obvious that the entries in rows
$i=0,1,\dots,b+c-2e-1$ are polynomials in $x$, and for $i\ge b+c-2e$ we
have: $c-e\ge 2c-b\ge0$ by our assumptions, 
$e+i-j-1\ge b+c-e-j-1\ge c-e\ge 0$ if
$j<b$, $b-2e\ge 2c-b\ge0$ by our assumptions, and $2e+i-j-1\ge
b+c-j-1\ge0$ if $j<b+c$. This explains the term $(2e-c)$ in (3.47).

In order to explain the term $(2e+c-b)$ in (3.47), we claim that for
$s=0,1,\dots,2e+c-b-1$ we have
$$\align
&\sum _{i=0} ^{b-2e+s}
{{\left( -1 \right) }^{b + c + e + i + s+1}}
{{ 
     \left(  b + c - 2 e - i -1\right) !\, \left(  b - e - i -1\right) !}\over
     {\left(  2 e - 2 s -2\right) !\, \left(  b - i - s -1\right) !}}\\
&\hskip3cm\cdot
\frac {     \left(  e - s -1\right) !\, \left(  2e+c- b  - s -1\right) !} 
  {     \left( b - 2 e - i + s \right) !}
\cdot \big(\text {row $i$ of $\De_1(-e;b,c,e)$}\big)\\
&+
\sum _{i=b-e} ^{\min\{b+c-2e-1,b-s-1\}}
2\, {{\left( -1 \right) }^{b + c + i}} 
{{     \left(  b + c - 2 e - i -1\right) !\, \left(  e - s -1\right) !}\over 
   {\left( i+e-b  \right) !\, \left(  2 e - 2 s -2\right) !}}\\
&\hskip2cm\cdot
\frac {     \left(  2e+c- b - s -1\right) !\, 
     \left(  2e- b  + i - s -1\right) !}
 {     \left(  b - i - s -1\right) !}
\cdot \big(\text {row $i$ of $\De_1(-e;b,c,e)$}\big)\\
&+\chi(s\ge 2e-c)
\sum _{i=b+c-2e} ^{b-s-1}
{{({ \textstyle 1 - b - c + 2 e + i}) _{ 2 c - i - s-1}  \,
     ({ \textstyle 1 - e + s}) _{ b - i - s-1} }\over 
   {\left(  b - i - s -1\right) !\, 
     ({ \textstyle 2 - 2 e + 2 s}) _{ b - i - s-1} }}\\
&\hskip4cm
\cdot \big(\text {row $i$ of $\De_1(-e;b,c,e)$}\big)
\\
&\hskip4cm=0.\tag3.48
\endalign$$
The notation in this 
assertion needs some explanation. Whereas the meaning of
$\De_1(x;b,c,e)$ is clear if $c/2<e$, in the alternative case $e\le
c/2$ the symbol $\De_1(x;b,c,e)$ stands for the original determinant
$\De'(x;b,c)$, in abuse of notation. (An alternative way to see this
is to say that $\De_1(x;b,c,e)$, in that case, is also given by
(3.3), but because of $e\le c/2$ the bottom block is empty, and
therefore the middle block ranges over $i=c,c+1,\dots,b-1$.) 
As earlier, the truth symbol $\chi(.)$ is defined by 
$\chi(\Cal A)=1$ if $\Cal A$ is true and $\chi(\Cal A)=0$
otherwise. So, the third sum in (3.48) only appears if $s\ge
2e-c$.

Note that these are indeed $2e+c-b$ linear combinations of the rows,
which are linearly independent. 
The latter fact comes from the observation that for fixed $s$ 
the last nonzero coefficient in the linear combination (3.48) 
is the one for row $b-s-1$, regardless whether $s\ge 2e-c$ or not. 

Because of the condition $s\le 2e+c-b-1$, we have $b-2e+s\le c-1$, and
therefore the rows which are involved in the first sum in (3.48) are
from rows $0,1,\dots,c-1$, which form the top block in (3.3). 
The assumptions $e\le b-c$ and $b\le 2c$ imply $b-2e+s\ge 0$, and so the
bounds for the sum are proper bounds. 
Because of $b-e\ge c$, the rows which are involved in the second sum in 
(3.48) are
from rows $c,c+1,\dots,b+c-2e-1$, which form the middle block in (3.3).
The bounds for the sum are proper, since by our assumptions we have
$$s\le 2e+c-b-1\le e-1\le b-c-1\le c-1,\tag3.49$$
and therefore $b-e\le b-s-1$ and $b-e\le b+c-2e$ (including the
possibility that $c=e$, in which case the second sum in (3.48) is the
empty sum).
Finally, because of the condition $s\ge 0$, we have $b-s-1\le b-1$, and
therefore the rows which are involved in the third sum in (3.48) (if
existent) are
from rows $b+c-2e,b+c-2e+1,\dots,b-1$, which form the bottom block in (3.3).

Hence, in order to verify (3.48), we have to check
$$\align
&\sum _{i=0} ^{b-2e+s}
{{\left( -1 \right) }^{b + c + e + i + s+1}}
{{ 
     \left(  b + c - 2 e - i -1\right) !\, \left(  b - e - i -1\right) !}\over
     {\left(  2 e - 2 s -2\right) !\, \left(  b - i - s -1\right) !}}\\
&\hskip3cm\cdot
\frac {     \left(  e - s -1\right) !\, \left(  2e+c- b  - s -1\right) !} 
  {     \left( b - 2 e - i + s \right) !}
\binom {c-e}{i-j+c}\\
&+\chi(s\ge 2e-c)
\sum _{i=b+c-2e} ^{b-s-1}(-1)^{e+i+j+1}
{{({ \textstyle 1 - b - c + 2 e + i}) _{ 2 c - i - s-1}  \,
     ({ \textstyle 1 - e + s}) _{ b - i - s-1} }\over 
   {\left(  b - i - s -1\right) !\, 
     ({ \textstyle 2 - 2 e + 2 s}) _{ b - i - s-1} }}\\
&\hskip4cm
\cdot\frac { \left( c - e \right) !\, \left(   e + i - j -1\right) !} 
{\left( i - j +c\right) ! }
\\
&\hskip4cm=0,\tag3.50
\endalign$$
which is (3.48) restricted to the $j$-th column, $j=c,c+1,\dots,b-1$ (note
that all the entries in rows $b-e,b-e+1,\dots,b+c-2e-1$ of
$\De_1(-e;b,c,e)$ vanish
in such a column), 
and
$$\align
&\sum _{i=0} ^{b-2e+s}
{{\left( -1 \right) }^{b + c + e + i + s+1}}
{{ 
     \left(  b + c - 2 e - i -1\right) !\, \left(  b - e - i -1\right) !}\over
     {\left(  2 e - 2 s -2\right) !\, \left(  b - i - s -1\right) !}}\\
&\hskip3cm\cdot
\frac {     \left(  e - s -1\right) !\, \left(  2e+c- b  - s -1\right) !} 
  {     \left( b - 2 e - i + s \right) !}
\,2\binom {b-2e}{i-j+b}\\
&+
\sum _{i=b-e} ^{\min\{b+c-2e-1,b-s-1\}}
2\, {{\left( -1 \right) }^{b + c + i}} 
{{     \left(  b + c - 2 e - i -1\right) !\, \left(  e - s -1\right) !}\over 
   {\left( i+e-b  \right) !\, \left(  2 e - 2 s -2\right) !}}\\
&\hskip2cm\cdot
\frac {     \left(  2e+c- b - s -1\right) !\, 
     \left(  2e- b  + i - s -1\right) !}
 {     \left(  b - i - s -1\right) !}
\binom {b-2e}{i-j+b}\\
&+\chi(s\ge 2e-c)
\sum _{i=b+c-2e} ^{b-s-1}
2\,(-1)^{i+j+1}{{({ \textstyle 1 - b - c + 2 e + i}) _{ 2 c - i - s-1}  \,
     ({ \textstyle 1 - e + s}) _{ b - i - s-1} }\over 
   {\left(  b - i - s -1\right) !\, 
     ({ \textstyle 2 - 2 e + 2 s}) _{ b - i - s-1} }}\\
&\hskip4cm
\cdot \frac {\left( b - 2 e \right) !\, 
         \left(   2 e + i - j -1\right) !} 
{\left( i - j+b \right) !}
\\
&\hskip4cm=0,\tag3.51
\endalign$$
which is (3.48) restricted to the $j$-th column, $j=b,b+1,\dots,b+c-1$.

We start by proving (3.50). We remind the reader that here $j$ is
restricted to $c\le j<b$.
The two sums in (3.50) can be combined into a single sum. To be
precise, the left-hand side in (3.50) can be written as
$$\multline
\lim_{\de\to0}\Bigg(
\sum_{i = j-c}^{b-s-1}{{\left( -1 \right) }^{e + i + j+1}}
\frac {       \left( c - e \right) !\, ({ \textstyle 1 + {\de}}) _{  e + i -
j-1}} 
{\left( i - j +c\right) !\, \left(   b - i - s -1\right) !}\\
\cdot
{{          ({ \textstyle 1 - b - c + 2 e + {\de} + i}) _{  2c - i -
s-1}\,  
       ({ \textstyle 1 - e + {\de} + s}) _{  b - i - s-1} }\over 
     {       ({ \textstyle 2 - 2 e + {\de} + 2 s}) _{  b - i - s-1}
}}\Bigg).
\endmultline\tag3.52$$
This expression is in fact just a multiple of the expression (3.7).
So, in the same way as it was done for (3.7), it is shown that (3.52)
vanishes. All the previous arguments apply because the crucial
inequalities $e\le c$, $s\le c-1$, $s\le e-1$ are also valid here,
thanks to (3.49).
This establishes (3.50).

Similarly, for the proof of (3.51) 
(we remind the reader that here $j$ is restricted to $b\le j<b+c$), 
we observe that the three sums in
(3.51) can be combined into the single expression
$$\multline
\lim_{\de\to0}\Bigg(
\sum_{i = j-b}^{b-s-1}
2\, {{\left( -1 \right) }^{i + j+1}}
\frac {       \left( b - 2 e \right) !\, 
       ({ \textstyle 1 + {\de}}) _{  2 e + i - j-1}} 
{\left(   b - i - s -1\right) !\, ({i - j+b})!}\\
\cdot
{{ 
       ({ \textstyle 1 - b - c + 2 e + {\de} + i}) _{  2c - i -
s-1}\,  
       ({ \textstyle 1 - e + {\de} + s}) _{  b - i - s-1} }\over 
     {       ({ \textstyle 2 - 2 e + {\de} + 2 s}) _{  b - i - s-1}
}}\Bigg),
\endmultline$$
and note that this expression is a multiple of the expression (3.10).
That it vanishes is then seen in the same way as it was for (3.10).
Again, the inequalities (3.49) guarantee that all the previous
arguments go through.
This establishes (3.51), and thus completes the proof that $(x+e)$
divides $\De'(x;b,c)$ with multiplicity $m(e)$ as given in (3.47).

\smallskip
{\it Case 2: $b-c\le e\le b/2$}. By inspection of the expression
(3.46),
we see that we have to prove that $(x+e)^{m(e)}$
divides $\De'(x;b,c)$, where
$$m(e)= \cases (b-c)&b-c\le e\le c/2\\
(b-c)+(2e-c)&c/2<e\le b/2.\endcases\tag3.53$$
Note that the first case in (3.53) could be empty, but not the second
(except if $b=c$).

The term $(2e-c)$ 
in the ($c/2<e$)-case of (3.53) is basically explained in the same
way as in Case~1: 
We take $(x+e)$ out of rows $b+c-2e,b+c-2e+1,\dots,b-1$ of the
determinant $\De'(x;b,c)$ (clearly, such rows exist only if $c/2<e$),
and thus obtain the determinant (3.3), which we denoted by $\De_1(x;b,c,e)$.
As before, to see that 
this determinant has still entries which are
polynomial in $x$, it suffices to check that the entries in rows
$i=b+c-2e,b+c-2e+1,\dots,b-1$ are polynomials in $x$.
This follows in almost the same way as in Case~1: We have
$c-e\ge c-b/2\ge0$ by our assumptions,
$e+i-j-1\ge b+c-e-j-1\ge c-e\ge 0$ if
$j<b$, $b-2e\ge 0$ by assumption, and $2e+i-j-1\ge
b+c-j-1\ge0$ if $j<b+c$. This explains the term $(2e-c)$ in (3.53).

In order to explain the term $(b-c)$ in (3.53), we claim that for
$s=0,1,\dots,b-c-1$ we have
$$\align
&\sum _{i=0} ^{2c-b+s}
{{\left( -1 \right) }^{c + i}}
{{ \left(  b - c - s -1\right) !\, 
     ({ \textstyle 1 + c - 2 e + s}) _{ b - i - s-1}  }\over 
   {2 \left(  2 b - 2 c - 2 s -2\right) !}}\\
&\hskip3cm\cdot
\frac {     ({ \textstyle 1 - b + 2 c - i + s}) _{ b - c - s-1} }
 {     \left(  b - i - s -1\right) !}
\cdot \big(\text {row $i$ of $\De_1(-e;b,c,e)$}\big)\\
&+
\sum _{i=c} ^{b-s-1}
{{\left( -1 \right) }^{c + i}}
{{ \left(  b - c - s -1\right) !\, 
     ({ \textstyle 1 + c - 2 e + s}) _{ b - i - s-1}  }\over 
   {\left(  2 b - 2 c - 2 s -2\right) !}}\\
&\hskip3cm\cdot
\frac {     ({ \textstyle 1 - b + 2 c - i + s}) _{ b - c - s-1} }
 {     \left(  b - i - s -1\right) !}
\cdot \big(\text {row $i$ of $\De_1(-e;b,c,e)$}\big)
\\
&\hskip4cm=0\tag3.54
\endalign$$
if $s\ge 2e-c$, and
$$\align
&\sum _{i=0} ^{2c-b+s}
{{\left( -1 \right) }^{i + s+1}}
{{    \left(  b - c - s -1\right) !\,
 \left(  b + c - 2 e - i -1\right) !}\over 
   {\left(  2 b - 2 c - 2 s -2\right) !}}\\
&\hskip1cm\cdot
\frac {\left(  2e- c - s -1\right) !\, 
     ({ \textstyle 1 - b + 2 c - i + s}) _{ b - c - s-1} }
 {     \left(  b - i - s -1\right) !}
\cdot \big(\text {row $i$ of $\De_1(-e;b,c,e)$}\big)\\
&+
\sum _{i=c} ^{b+c-2e-1}
2\,{{\left( -1 \right) }^{ i + s+1}}
{{     \left(  b - c - s -1\right) !\,
 \left(  b + c - 2 e - i -1\right) !}\over 
   {\left(  2 b - 2 c - 2 s -2\right) !}}\\
&\hskip1cm\cdot
\frac {\left(  2e- c - s -1\right) !\, 
     ({ \textstyle 1 - b + 2 c - i + s}) _{ b - c - s-1} }
 {     \left(  b - i - s -1\right) !}
\cdot \big(\text {row $i$ of $\De_1(-e;b,c,e)$}\big)\\
&+
\sum _{i=b+c-2e} ^{b-s-1}
{{\left( -1 \right) }^{c + i}}
{{ \left(  b - c - s -1\right) !\, 
     ({ \textstyle 1 + c - 2 e + s}) _{ b - i - s-1}  }\over 
   {\left(  2 b - 2 c - 2 s -2\right) !}}\\
&\hskip3cm\cdot
\frac {     ({ \textstyle 1 - b + 2 c - i + s}) _{ b - c - s-1} }
 {     \left(  b - i - s -1\right) !}
\cdot \big(\text {row $i$ of $\De_1(-e;b,c,e)$}\big)
\\
&\hskip4cm=0\tag3.55
\endalign$$
if $s\le 2e-c$. In (3.54) and (3.55) we make the same convention as
in Case~1 of how to understand $\De_1(x;b,c,e)$ in the case that $e\le
c/2$.

It should be noted that in both cases these are indeed $b-c$ 
linear combinations of the rows, which are linearly independent. 

Let us first consider (3.54), i.e., in the following paragraphs we
assume $s\ge 2e-c$. 
Because of the condition $s\le b-c-1$, we have $2c-b+s\le c-1$, and
therefore the rows which are involved in the first sum in (3.54) are
from rows $0,1,\dots,c-1$, which form the top block in (3.3). 
Because of $2c-b+s\ge0$
the bounds for the sum are proper bounds. 
Since $s\ge 0$, we have $b-s-1\le b-1$, and
therefore the rows which are involved in the second sum in (3.54) are
from rows $c,c+1,\dots,b-1$, which form the ``middle" block in (3.3) if
$s\ge 2e-c$ (recall: the bottom block is empty in this case). 
Finally, the assumption $s\le b-c-1$ implies $c\le
b-s-1$, and so the bounds for the sum are proper.

Hence, in order to verify (3.54), we have to check
$$\multline
\sum _{i=0} ^{2c-b+s}
{{\left( -1 \right) }^{c + i}}
{{ \left(  b - c - s -1\right) !\, 
     ({ \textstyle 1 + c - 2 e + s}) _{ b - i - s-1}  }\over 
   {2 \left(  2 b - 2 c - 2 s -2\right) !}}\\
\cdot
\frac {     ({ \textstyle 1 - b + 2 c - i + s}) _{ b - c - s-1} }
 {     \left(  b - i - s -1\right) !}
\binom {c-e}{i-j+c}\,=\,0,
\endmultline\tag3.56$$
which is (3.54) restricted to the $j$-th column, $j=c,c+1,\dots,b-1$ 
(note that this is indeed the restriction of (3.54) to the $j$-th
column, $c\le j<b$, since, due to $0\le c-e\le 2c-b<2c-j\le i-j+c$, 
the entries $\binom {c-e} {i-j+c}$ in rows
$c,c+1,\dots,b-s-1$ of $\De_1(-e;b,c,e)$ vanish in such a column),
and
$$\align
&\sum _{i=0} ^{2c-b+s}
{{\left( -1 \right) }^{c + i}}
{{ \left(  b - c - s -1\right) !\, 
     ({ \textstyle 1 + c - 2 e + s}) _{ b - i - s-1}  }\over 
   {2 \left(  2 b - 2 c - 2 s -2\right) !}}\\
&\hskip3cm\cdot
\frac {     ({ \textstyle 1 - b + 2 c - i + s}) _{ b - c - s-1} }
 {     \left(  b - i - s -1\right) !}
\,2\binom {b-2e}{i-j+b}\\
&+
\sum _{i=c} ^{b-s-1}
{{\left( -1 \right) }^{c + j+1}}
{{ \left(  b - c - s -1\right) !\, 
     ({ \textstyle 1 + c - 2 e + s}) _{ b - i - s-1}  }\over 
   {\left(  2 b - 2 c - 2 s -2\right) !}}\\
&\hskip3cm\cdot
\frac {     ({ \textstyle 1 - b + 2 c - i + s}) _{ b - c - s-1} }
 {     \left(  b - i - s -1\right) !}
\binom {b-2e}{i-j+b}\\
\\
&\hskip4cm=0,\tag3.57
\endalign$$
which is (3.54) restricted to the $j$-th column,
$j=b,b+1,\dots,b+c-1$.

In order to verify (3.56) (we remind the reader that here $j$ is
restricted to $c\le j<b$), we rewrite the left-hand side in a fancier
way as
$$\multline
\lim_{\de\to0}\Bigg(
\sum_{i = j-c }^{2c-b+s}{{\left( -1 \right) }^{c + i}} 
{{
      \left(  b - c - s -1\right) !\, 
       ({ \textstyle 1 + c - 2 e + {\de} + s}) _{ b - i - s-1} }
      \over {2        \left(  2 b - 2 c - 2 s -2\right) !}}\\
\cdot 
\frac {       ({ \textstyle 1 - b + 2 c + {\de} - i + s}) _{ b - c - s-1}
\,(c-e)! }
 {          \left(  b - i - s -1\right) !\, \left( i - j+c \right) !\, 
   ({ \textstyle 1 + {\de}}) _{j-i-e} }
\Bigg),
\endmultline\tag3.58$$
and convert the series into hypergeometric notation,
$$\multline
\lim_{\de\to0}\Bigg(
{{\left( -1 \right) }^j}
{{\left( c - e \right) !\, 
     ({ \textstyle 1 + c - 2 e + {\de} + s}) _{ b + c - j - s-1} \,
     ({ \textstyle 1 - b + 3 c + {\de} - j + s}) _{ b - c - s-1} }\over 
   {2 \left(  2 b - 2 c - 2 s -2\right) !\, 
     \left(  b + c - j - s -1\right) !}}\\
\times
\frac {  \left(  b - c - s -1\right) ! }
 {     ({ \textstyle 1 + {\de}}) _{c - e} }
 {} _{3} F _{2} \!\left [ \matrix { -c + e - {
      \de}, b - 3 c - {\de} + j - s, 1 - b - c + j + s}\\ { 1 - 2 c - {
      \de} + j, 1 - b - 2 c + 2 e - {\de} + j}\endmatrix ; {\displaystyle
      1}\right ]
\Bigg).
\endmultline$$
To the $_3F_2$-series we apply, once again, the transformation
formula (3.11). Thus we obtain the expression
$$\multline
\lim_{\de\to0}\Bigg(
{{\left( -1 \right) }^j}
{{   \left( c - e \right) !\, 
     ({ \textstyle 1 + c - 2 e + {\de} + s}) _{ b + c - j - s-1} \,
     ({ \textstyle 1 - b + 3 c + {\de} - j + s}) _{ b - c - s-1} }\over 
   {2 \,\left(  2 b - 2 c - 2 s -2\right) !\, 
     \left(  b + c - j - s -1\right) !}}\\
\times
\frac {\left(  b - c - s -1\right) !}
 {     ({ \textstyle 1 + {\de}}) _{c - e} }
\frac {\Gamma(1 - b - 2 c + 2 e - {\de} + j)\, \Gamma(1 - b + c + e)}
 {\Gamma(1 - b - c + e + j)\, \Gamma(1 - b + 2 e - {\de})}\\
\times
{} _{3} F _{2} \!\left [ \matrix { -c + e - {
      \de}, 1 - b + c + s, b - c - {\de} - s}\\ { 1 - 2 c - {\de} + j, 1
      - b + 2 e - {\de}}\endmatrix ; {\displaystyle 1}\right ]
\Bigg)
\endmultline$$
for the left-hand side in (3.56). The $_3F_2$-series in this
expression terminates because of the upper parameter $1-b+c+s$, which
is a nonpositive integer because of an assumption.
Hence it is well-defined.
The complete expression vanishes because of the occurence of the term
$\Gamma(1 - b - c + e + j)$ in the denominator. For, by our assumptions, 
we have $1-b-c+e+j\le e-c\le0$, and so the gamma function equals
$\infty$. This establishes (3.56).

For proving (3.57) (we remind the reader that here $j$ is
restricted to $b\le j<b+c$), we observe that the two sums in (3.57) can be
combined into the single expression
$$\multline
\lim_{\de\to0}\Bigg(
\sum_{i = j-b}^{2c-b+s}{{\left( -1 \right) }^{c + i}}
{{ 
       \left(  b - c - s -1\right) !\, 
       ({ \textstyle 1 + c - 2 e + {\de} + s}) _{ b - i - s-1}}
      \over {       \left(  2 b - 2 c - 2 s -2\right) !}}\\
\cdot
\frac {        ({ \textstyle 1 - b + 2 c + {\de} - i + s}) _{ b - c - s-1}
\,      ({ \textstyle 1 - 2 e + {\de} - i + j}) _{i - j+b} }
 {       \left(  b - i - s -1\right) !\, \left( i - j+b \right) !}
\Bigg).
\endmultline\tag3.59$$
Using hypergeometric notation, this expression can be rewritten as
$$\multline
\lim_{\de\to0}\Bigg(
{{\left( -1 \right) }^{b + c + j}}
{{         ({ \textstyle 1 + c - 2 e + {\de} + s}) _{ 2 b - j - s-1} \,
     ({ \textstyle 1 + 2 c + {\de} - j + s}) _{ b - c - s-1} }\over 
   {\left(  2 b - 2 c - 2 s -2\right) !}}\\
\times
\frac {  \left(  b - c - s -1\right) !}
 {     \left(  2 b - j - s -1\right) !}
     {} _{3} F _{2} \!\left [ \matrix { -b + 2 e - {
      \de}, -2 c - {\de} + j - s},1 - 2 b + j + s\\
   { 1 - 2 b - c + 2 e - {\de} + j, 1
      - b - c - {\de} + j}\endmatrix ; {\displaystyle 1}\right ]
\Bigg).
\endmultline$$
The $_3F_2$-series can be summed by means of the Pfaff-Saalsch\"utz
summation (3.8). We have to apply the case where $n=2b-j-s-1$, which is
indeed a nonnegative integer because of $2b-j-s-1\ge b-c-s\ge1$. This gives
$$\multline
\lim_{\de\to0}\Bigg(
{{\left( -1 \right) }^{b + c + j}}
{{      ({ \textstyle 1 - b - c + j}) _{ 2 b - j - s-1}  \,
     ({ \textstyle 1 - 2 b + c + 2 e + s}) _{ 2 b - j - s-1}  }\over 
   {\left(  2 b - 2 c - 2 s -2\right) !}}\\
\times
\frac { \left(  b - c - s -1\right) !\, 
    ({ \textstyle 1 + c - 2 e + {\de} + s}) _{ 2 b - j - s-1} \,
     ({ \textstyle 1 + 2 c + {\de} - j + s}) _{ b - c - s-1} }
 {         \left(  2 b - j - s -1\right) !\,
 ({ \textstyle 1 - 2 b - c + 2 e - {\de} + j}) _{ 2 b - j - s-1}\,
     ({ \textstyle 1 - b + c + {\de} + s}) _{ 2 b - j - s-1} }
\Bigg)
\endmultline$$
as an equivalent expression for (3.57).
It vanishes because of the occurence of the term
$$(1-b-c+j)_{2b-j-s-1}=(1-b-c+j)(2-b-c+j)\cdots (b-c-s-1)$$
in the numerator.
For, by our assumptions, we have $1-b-c+j\le0$, and we have $b-c-s-1\ge0$.
This establishes (3.57).

Now let us consider (3.55), i.e., in the following paragraphs we
assume $s\le 2e-c$. 
In the same way as for (3.54), it is checked that the the rows which
are involved in the first sum in (3.55) are from rows
$0,1,\dots,c-1$, and that the bounds for the sum are proper bounds.
Clearly, the rows which are involved in the second sum in (3.55)
are from rows $c,c+1,\dots,b+c-2e-1$, which form the middle block in
(3.3). The assumption $e\le b/2$ guarantees that the bounds for the
sum are proper (including the
possibility that $e=b/2$, in which case the sum is the empty sum).
Finally, since $s\ge0$, the rows which are involved in the third sum in (3.55)
are from rows $b+c-2e,b+c-2e+1,\dots,b-1$, which form the bottom block in
(3.3). That the bounds for the sum are proper follows from the 
condition $s\le 2e-c$ (including the
possibility that $s=2e-c$, in which case the sum is the
empty sum).

Hence, in order to verify (3.55), we have to check
$$\align
&\sum _{i=0} ^{2c-b+s}
{{\left( -1 \right) }^{i + s+1}}
{{    \left(  b - c - s -1\right) !\,
 \left(  b + c - 2 e - i -1\right) !}\over 
   {\left(  2 b - 2 c - 2 s -2\right) !}}\\
&\hskip2cm\cdot
\frac {\left(  2e- c - s -1\right) !\, 
     ({ \textstyle 1 - b + 2 c - i + s}) _{ b - c - s-1} }
 {     \left(  b - i - s -1\right) !}
\binom {c-e}{i-j+c}\\
&+
\sum _{i=b+c-2e} ^{b-s-1}
{{\left( -1 \right) }^{c + e+j+1}}
{{ \left(  b - c - s -1\right) !\, 
     ({ \textstyle 1 + c - 2 e + s}) _{ b - i - s-1}  }\over 
   {\left(  2 b - 2 c - 2 s -2\right) !}}\\
&\hskip3cm\cdot
\frac {     ({ \textstyle 1 - b + 2 c - i + s}) _{ b - c - s-1} }
 {     \left(  b - i - s -1\right) !}
\frac { \left( c - e \right) !\, \left(   e + i - j -1\right) !} 
{\left( i - j +c\right) ! }
\\
&\hskip4cm=0,\tag3.60
\endalign$$
which is (3.55) restricted to the $j$-th column, $j=c,c+1,\dots,b-1$ 
(again note that all the entries in rows $c,c+1,\dots,b+c-2e-1$ of
$\De_1(-e;b,c,e)$ vanish in such a column), 
and
$$\align
&\sum _{i=0} ^{2c-b+s}
{{\left( -1 \right) }^{i + s+1}}
{{    \left(  b - c - s -1\right) !\,
 \left(  b + c - 2 e - i -1\right) !}\over 
   {\left(  2 b - 2 c - 2 s -2\right) !}}\\
&\hskip2cm\cdot
\frac {\left(  2e- c - s -1\right) !\, 
     ({ \textstyle 1 - b + 2 c - i + s}) _{ b - c - s-1} }
 {     \left(  b - i - s -1\right) !}
\,2\binom {b-2e}{i-j+b}\\
&+
\sum _{i=c} ^{b+c-2e-1}
2\,{{\left( -1 \right) }^{ i + s+1}}
{{     \left(  b - c - s -1\right) !\,
 \left(  b + c - 2 e - i -1\right) !}\over 
   {\left(  2 b - 2 c - 2 s -2\right) !}}\\
&\hskip2cm\cdot
\frac {\left(  2e- c - s -1\right) !\, 
     ({ \textstyle 1 - b + 2 c - i + s}) _{ b - c - s-1} }
 {     \left(  b - i - s -1\right) !}
\binom {b-2e}{i-j+b}\\
&+
\sum _{i=b+c-2e} ^{b-s-1}
2\,{{\left( -1 \right) }^{c + j+1}}
{{ \left(  b - c - s -1\right) !\, 
     ({ \textstyle 1 + c - 2 e + s}) _{ b - i - s-1}  }\over 
   {\left(  2 b - 2 c - 2 s -2\right) !}}\\
&\hskip3cm\cdot
\frac {     ({ \textstyle 1 - b + 2 c - i + s}) _{ b - c - s-1} }
 {     \left(  b - i - s -1\right) !}
\frac {\left( b - 2 e \right) !\, 
         \left(   2 e + i - j -1\right) !} 
{\left( i - j+b \right) !}
\\
&\hskip4cm=0,\tag3.61
\endalign$$
which is (3.55) restricted to the $j$-th column, $j=b,b+1,\dots,b+c-1$.

Both identities are now easily verified. In fact, the left-hand side
of (3.60) can be written as $\lim_{\de\to0}(2 E_1/\de)$,
where $E_1$ is the expression in big parentheses in (3.58). Likewise, 
the left-hand side
of (3.61) can be written as $\lim_{\de\to0}(2 E_2/\de)$,
where $E_2$ is the expression in big parentheses in (3.59). The same
arguments as in the proofs of (3.56) and (3.57) then show that (3.60)
and (3.61) vanish.

This completes the proof that $(x+e)$
divides $\De'(x;b,c)$ with multiplicity $m(e)$ as given in (3.53).

\smallskip
{\it Case 3: $b/2\le e\le c$}.
By inspection of the expression (3.46),
we see that we have to prove that $(x+e)^{m(e)}$
divides $\De'(x;b,c)$, where
$$m(e)=\cases (b-c)+(2b-2e-c)&b/2\le e<b-c/2\\
(b-c)&b-c/2\le e\le c.\endcases\tag3.62$$
Note that the second case in (3.62) could be empty, but not the first
(except if $b=c$).

As in Case~4 of the proof of Lemma~1,
in order to extract the appropriate number of factors $(x+e)$ out of the
determinant $\De'(x;b,c)$, 
we start with the modified determinant (3.13) with $e=b$. Recall,
that the choice of $e=b$ has the effect that the bottom block
in (3.13) is empty. If $e<b-c/2$,
we take $(x+e)$ out of rows $2e+c-b,2e+c-b+1,\dots,b-1$ (such
rows only exist under the assumption $e<b-c/2$), and obtain the
determinant in (3.38), which we denoted by $\De_3(x;b,c,e)$. 
Obviously, we have taken out
$(x+e)^{2b-2e-c}$. 
The remaining determinant has still entries which are
polynomial in $x$. 
For, it is obvious that the entries in rows
$i=0,1,\dots,2e+c-b-1$ are polynomials in $x$, and for $i\ge 2e+c-b$ we
have: $i-e\ge e+c-b\ge c-b/2\ge0$ by our assumptions, 
 $e+c-j-1\ge e+c-b\ge 0$ if
$j<b$, $i+b-c-2e\ge0$, and $2e+c-j-1\ge
2e-b\ge 0$ if $j<b+c$. This explains the term $(2b-2e-c)$ in the
($e<b-c/2$)-case of (3.62).

In order to explain the term $(b-c)$ in (3.62), we claim that for
$s=0,1,\dots,b-c-1$ we have
$$\align
&\sum _{i=0} ^{2c-b+s}
{{\left( -1 \right) }^{c + i}}
{{ \left(  b - c - s -1\right) !\,
     ({ \textstyle 1 + c - 2 e + s}) _{ b - i - s-1} \,
     ({ \textstyle 1 - b + 2 c - i + s}) _{ b - c - s-1} }\over 
   {2\,\left(  2 b - 2 c - 2 s -2\right) !\,
     \left(  b - i - s -1\right) !}}\\
&\hskip4cm
\cdot \big(\text {row $i$ of $\De_3(-e;b,c,e)$}\big)\\
&\sum _{i=c} ^{b-s-1}
{{\left(  b - c - s -1\right) !\,
     ({ \textstyle 1 - c + i}) _{ b - c - s-1} \,
     ({ \textstyle 1 - 2 b + c + 2 e + s}) _{ b - i - s-1} }\over 
   {\left(  2 b - 2 c - 2 s -2\right) !\,\left(  b - i - s -1\right) !}
   }\\
&\hskip4cm
\cdot \big(\text {row $i$ of $\De_3(-e;b,c,e)$}\big)
\\
&\hskip4cm=0\tag3.63
\endalign$$
if $s\ge 2b-2e-c$, and
$$\align
&\sum _{i=0} ^{2c-b+s}
{{\left( -1 \right) }^{i + s+1}}
{{     \left(  b - c - s -1\right) !\,
 \left(  b + c - 2 e - i -1\right) !\,
\left(  2e- c - s -1\right) !}\over 
   {\left(  2 b - 2 c - 2 s -2\right) !}   }\\
&\hskip3cm
\cdot\frac {     ({ \textstyle 1 - b + 2 c - i + s}) _{ b - c - s-1} }
 {\left(  b - i - s -1\right) !}
\cdot \big(\text {row $i$ of $\De_3(-e;b,c,e)$}\big)\\
&\sum _{i=c} ^{2e+c-b-1}
2\,{{\left( -1 \right) }^{c + s}}
{{   \left(  b - c - s -1\right) !\,
 \left(  2e+c- b - i -1\right) !\,
       \left(  2 b - c - 2 e - s -1\right) !}\over 
   {\left(  2 b - 2 c - 2 s -2\right) !}
   }\\
&\hskip3cm
\cdot\frac {     ({ \textstyle 1 - c + i}) _{ b - c - s-1} }
 {\left(  b - i - s -1\right) !}
\cdot \big(\text {row $i$ of $\De_3(-e;b,c,e)$}\big)\\
&\sum _{i=2e+c-b} ^{b-s-1}
{{\left(  b - c - s -1\right) !\,
     ({ \textstyle 1 - c + i}) _{ b - c - s-1}\,
    ({ \textstyle 1 - 2 b + c + 2 e + s}) _{ b - i - s-1}  }\over 
   {\left(  2 b - 2 c - 2 s -2\right) !\,\left(  b - i - s -1\right) !}
   }\\
&\hskip4cm
\cdot \big(\text {row $i$ of $\De_3(-e;b,c,e)$}\big)
\\
&\hskip4cm=0\tag3.64
\endalign$$
if $s\le 2b-2e-c$. In (3.63) and (3.64) we make a similar convention as
in Case~1 of how to understand $\De_3(x;b,c,e)$ in the case that $e\ge
b-c/2$.

It should be noted that in both cases these are indeed $b-c$ 
linear combinations of the rows, which are linearly independent. 

Let us first consider (3.63), i.e., in the following paragraphs we
assume $s\ge 2b-2e-c$. 
In the same way as for (3.54) it is seen that
the rows which are involved in the first sum in (3.63) are
from rows $0,1,\dots,c-1$, which form the top block in (3.38), and
that the bounds for the sum are proper bounds. 
Also in the same way, it is seen
that the rows which are involved in the second sum in (3.63) are
from rows $c,c+1,\dots,b-1$, which form the ``middle" block in (3.38) if
$s\ge 2b-2e-c$ (recall: the bottom block is empty in this case), and
that the bounds for the sum are proper.

Hence, in order to verify (3.63), we have to check
$$\multline
\sum _{i=0} ^{2c-b+s}
{{\left( -1 \right) }^{c + i}}
{{ \left(  b - c - s -1\right) !\,
     ({ \textstyle 1 + c - 2 e + s}) _{ b - i - s-1} \,
     ({ \textstyle 1 - b + 2 c - i + s}) _{ b - c - s-1} }\over 
   {2\,\left(  2 b - 2 c - 2 s -2\right) !\,
     \left(  b - i - s -1\right) !}}\\
\cdot\binom {c-e}{i-j+c}\,=\,0,\hskip2cm
\endmultline\tag3.65$$
which is (3.63) restricted to the $j$-th column, $j=c,c+1,\dots,b-1$ 
(note that this is indeed the restriction of (3.63) to the $j$-th
column, $c\le j<b$, since, due to $0\le i-e\le i-b/2\le i-b+c< i-j+c$, 
the entries $\binom {i-e} {i-j+c}$ in rows
$c,c+1,\dots,b-s-1$ of $\De_3(-e;b,c,e)$ vanish in such a column),
and
$$\align
&\sum _{i=0} ^{2c-b+s}
{{\left( -1 \right) }^{c + i}}
{{ \left(  b - c - s -1\right) !\,
     ({ \textstyle 1 + c - 2 e + s}) _{ b - i - s-1} \,
     ({ \textstyle 1 - b + 2 c - i + s}) _{ b - c - s-1} }\over 
   {2\,\left(  2 b - 2 c - 2 s -2\right) !\,
     \left(  b - i - s -1\right) !}}\\
&\hskip4cm
\cdot 2\binom {b-2e}{i-j+b}\\
&\sum _{i=c} ^{b-s-1}
{{\left(  b - c - s -1\right) !\,
     ({ \textstyle 1 - c + i}) _{ b - c - s-1} \,
     ({ \textstyle 1 - 2 b + c + 2 e + s}) _{ b - i - s-1} }\over 
   {\left(  2 b - 2 c - 2 s -2\right) !\,\left(  b - i - s -1\right) !}
   }\\
&\hskip4cm
\cdot \binom {i+b-c-2e}{i-j+b}\\
&\hskip4cm=0,\tag3.66
\endalign$$
which is (3.63) restricted to the $j$-th column,
$j=b,b+1,\dots,b+c-1$.

Identity (3.65) is now easily verified. In fact, the left-hand side
of (3.65) can be rewritten as the expression (3.58). It vanishes
since the crucial inequalities $1-b+c+s\le0$ and $e-c\le0$
are also valid here.

For verifying (3.66) we have to do little more work. 
We remind the reader that here $j$ is
restricted to $b\le j<b+c$.
We consider first the second term in (3.66).
We replace the
binomial $\binom {i+b-c-2e}{i-j+b}$ by the expansion
$\sum _{\ell=c}
^{i}\binom {b-2e}{\ell-j+b} \binom {i-c}{i-\ell}$, the 
equality of binomial and
expansion being again due to Chu--Vandermonde summation. Then we
interchange the summations over $i$ and $\ell$, and write the sum
over $i$ in hypergeometric notation. This gives
$$\multline
\sum _{\ell=c} ^{b-s-1}
{{     (b - c - s-1)!  \,
     ({ \textstyle 1 - 2 e + j - \ell}) _{\ell-j+b}  \,
     ({ \textstyle 1 - c + \ell}) _{ b - c - s-1}  }\over 
   {( b- j + \ell)!  \,
     (2 b - 2 c - 2 s-2)!  }}\\
\cdot
\frac {     ({ \textstyle 1 - 2 b + c + 2 e + s}) _{ b - \ell - s-1} }
 {     (b - \ell - s-1)! }
{} _{2} F _{1} \!\left [ \matrix { b - 2 c + \ell - s, 1 - b + \ell + s}\\ { 1 +
      b - c - 2 e + \ell}\endmatrix ; {\displaystyle 1}\right ]
\endmultline$$
as an equivalent expression for the second term in (3.66). 
The $_2F_1$-series can be
evaluated by the hypergeometric form (3.20) of the Chu--Vandermonde
summation. Thus we obtain the expression
$$\multline
\sum _{\ell=c} ^{b-s-1}
{{\left( -1 \right) }^{c + \ell}}  
{{     \left(  b - c - s -1\right) !\, 
     ({ \textstyle 1 + c - 2 e + s}) _{ b - \ell - s-1}  \,
     ({ \textstyle 1 - b + 2 c - \ell + s}) _{ b - c - s-1} }\over 
   {\left(  2 b - 2 c - 2 s -2\right) !\, \left(  b - \ell - s -1\right) !\,}
   }\\
\cdot {{b - 2 e}\choose {\ell-j+b}}.\hskip3cm
\endmultline\tag3.67$$
Now it is straight-forward to see that the 
first term in (3.66) and the above expression for the second term in
(3.66) can be combined into the single expression (3.59). Then the
same arguments as before apply to show that this expression vanishes
as well in the current case. For, the crucial inequalities
$2b-j-s-1\ge0$, $1-b-c+j\le0$, and $b-c-s-1\ge0$ are also valid here.
This establishes (3.66).

Now let us consider (3.64), i.e., in the following paragraphs we
assume $s\le 2b-2e-c$. 
In the same way as for (3.54), it is checked that the the rows which
are involved in the first sum in (3.64) are from rows
$0,1,\dots,c-1$, and that the bounds for the sum are proper bounds.
Clearly, the rows which are involved in the second sum in (3.64)
are from rows $c,c+1,\dots,2e+c-b-1$, which form the middle block in
(3.38). The assumption $e\ge b/2$ guarantees that the bounds for the
sum are proper (including the
possibility that $e=b/2$, in which case the sum is the empty sum).
Finally, since $s\ge0$, the rows which are involved in the third sum in (3.64)
are from rows $2e+c-b,2e+c-b+1,\dots,b-1$, which form the bottom block in
(3.38). That the bounds for the sum are proper follows from the 
condition $s\le 2b-2e-c$ (including the
possibility that $s=2b-2e-c$, in which case the sum is the
empty sum).

Hence, in order to verify (3.64), we have to check
$$\align
&\sum _{i=0} ^{2c-b+s}
{{\left( -1 \right) }^{i + s+1}}
{{     \left(  b - c - s -1\right) !\,
 \left(  b + c - 2 e - i -1\right) !\,
\left(  2e- c - s -1\right) !}\over 
   {\left(  2 b - 2 c - 2 s -2\right) !}   }\\
&\hskip3cm
\cdot\frac {     ({ \textstyle 1 - b + 2 c - i + s}) _{ b - c - s-1} }
 {\left(  b - i - s -1\right) !}
\binom {c-e}{i-j+c}\\
&\sum _{i=2e+c-b} ^{b-s-1}(-1)^{c+e+j+1}
{{\left(  b - c - s -1\right) !\,
     ({ \textstyle 1 - c + i}) _{ b - c - s-1}\,
    ({ \textstyle 1 - 2 b + c + 2 e + s}) _{ b - i - s-1}  }\over 
   {\left(  2 b - 2 c - 2 s -2\right) !\,\left(  b - i - s -1\right) !}
   }\\
&\hskip4cm
\cdot \frac {(i-e)!\,(e+c-j-1)!} {(i-j+c)!}
\\
&\hskip4cm=0,\tag3.68
\endalign$$
which is (3.64) restricted to the $j$-th column, $j=c,c+1,\dots,b-1$,
(recall that all the entries in rows $c,c+1,\dots,2e+c-b-1$ of 
$\De_3(-e;b,c,e)$ vanish in such a column),
and 
$$\align
&\sum _{i=0} ^{2c-b+s}
{{\left( -1 \right) }^{i + s+1}}
{{     \left(  b - c - s -1\right) !\,
 \left(  b + c - 2 e - i -1\right) !\,
\left(  2e- c - s -1\right) !}\over 
   {\left(  2 b - 2 c - 2 s -2\right) !}   }\\
&\hskip3cm
\cdot\frac {     ({ \textstyle 1 - b + 2 c - i + s}) _{ b - c - s-1} }
 {\left(  b - i - s -1\right) !}
\,2\binom {b-2e}{i-j+b}\\
&\sum _{i=c} ^{2e+c-b-1}
2\,{{\left( -1 \right) }^{c + s}}
{{   \left(  b - c - s -1\right) !\,
 \left(  2e+c- b - i -1\right) !\,
       \left(  2 b - c - 2 e - s -1\right) !}\over 
   {\left(  2 b - 2 c - 2 s -2\right) !}
   }\\
&\hskip3cm
\cdot\frac {     ({ \textstyle 1 - c + i}) _{ b - c - s-1} }
 {\left(  b - i - s -1\right) !}
\binom {i+b-c-2e}{i-j+b}\\
&\sum _{i=2e+c-b} ^{b-s-1}2\,(-1)^{c+j+1}
{{\left(  b - c - s -1\right) !\,
     ({ \textstyle 1 - c + i}) _{ b - c - s-1}\,
    ({ \textstyle 1 - 2 b + c + 2 e + s}) _{ b - i - s-1}  }\over 
   {\left(  2 b - 2 c - 2 s -2\right) !\,\left(  b - i - s -1\right) !}
   }\\
&\hskip4cm
\cdot \frac {(i+b-c-2e)!\,(2e+c-j-1)!} {(i-j+b)!}
\\
&\hskip4cm=0,\tag3.69
\endalign$$
which is (3.64) restricted to the $j$-th column, $j=b,b+1,\dots,b+c-1$.

We start with the proof of (3.68). We remind the reader that here $j$ is
restricted to $c\le j<b$.
The strategy is analogous to the
one used in the proof of (3.66) just before. We recast the second
term in (3.68) by replacing the subterm $(i-e)!\,(e+c-j-1)!/(i-j+c)!$ by
the expression
$$\lim_{\de\to0}\Bigg(\frac {1} {\de}\sum _{\ell=c} ^{i}\frac
{(c-e)!} {(\ell-j+c)!\,(1+\de)_{j-\ell-e}}\binom {i-c}{i-\ell}\Bigg),$$
the equality of second term and this expression following again from
Chu--Van\-der\-monde summation. Then we interchange sums over $i$ and
$\ell$, and evaluate the now inner sum by Chu--Vandermonde summation
(3.20). The computation is essentially the same as before in the
proof of (3.66). Eventually, we obtain the expression (3.67), with
the binomial $\binom {b-2e}{\ell-j+b}$ replaced by
$(c-e)!/\big(\de\,(\ell-j+c)!\, (1+\de)_{j-\ell-e}\big)$. Therefore, 
this expression and the first term in (3.68) can be combined into the
single expression $\lim_{\de\to0}2E_1/\de$, 
where $E_1$ is the expression in big parentheses in (3.58). Then we
may follow the arguments which proved that (3.58) vanishes, since the
crucial inequalities $1-b+c+s\le0$ and $e-c\le0$ are also valid here.
This establishes (3.68).

Now we turn to (3.69). We remind the reader that here $j$ is
restricted to $b\le j<b+c$.
We proceed again in the same way as in the
proof of (3.66). Once more using Vandermonde summation, we 
replace the
binomial $\binom {i+b-c-2e}{i-j+b}$ in the second term in (3.69) 
by the expansion $\sum _{\ell=c}
^{i}\binom {b-2e}{\ell-j+b} \binom {i-c}{i-\ell}$, and we replace the
subterm $(i+b-c-2e)!\,(2e+c-j-1)!/(i-j+b)!$ in the third term in
(3.69) by the expression
$$\lim_{\de\to0}\Bigg(\frac {1} {\de}\sum _{\ell=c} ^{i}
\frac {(1+j-\ell-2e+\de)_{\ell-j+b}} {(\ell-j+b)!}
\binom {i-c}{i-\ell}\Bigg).$$
Then we interchange sums over $i$ and $\ell$, and evaluate the now inner
sums over $i$ by the same instance of the Chu--Vandermonde summation
(3.20). Eventually, it is seen that the three terms on the left-hand
side of (3.69) can be combined into the single expression 
$\lim_{\de\to0}(2 E_2/\de)$,
where $E_2$ is the expression in big parentheses in (3.59). The same
arguments as in the proof of (3.57) can now be used since the crucial
inequalities $2b-j-s-1\ge0$, $1-b-c+j\le0$, and $b-c-s-1\ge0$
are also valid here. This establishes (3.69) and 
completes the proof that $(x+e)$
divides $\De'(x;b,c)$ with multiplicity $m(e)$ as given in (3.62).

\smallskip
{\it Case 4: $c\le e\le (b+c)/2$}.
By inspection of the expression (3.46),
we see that we have to prove that $(x+e)^{m(e)}$
divides $\De'(x;b,c)$, where
$$m(e)=\cases (b+c-2e)+(2b-2e-c)&c\le e< b-c/2\\
(b+c-2e)&b-c/2\le e\le (b+c)/2.\endcases
\tag3.70$$
Note that the first case in (3.70) could be empty, but not the second,
because of $b\le 2c$.

As in Case~3,
in order to explain the term $2b-2e-c$ in the ($e<b-c/2$)-case of
(3.70), 
we start with the determinant (3.13) with $e=b$. 
We take again $(x+e)$ out of rows $2e+c-b,2e+c-b+1,\dots,b-1$ (such
rows only exist under the assumption $e<b-c/2$), and obtain the
determinant in (3.38), which we denoted by $\De_3(x;b,c,e)$. 
Obviously, we have taken out $(x+e)^{2b-2e-c}$. 
As before, to see that 
this determinant has still entries which are
polynomial in $x$, it suffices to check that the entries in rows
$i=2e+c-b,2e+c-b+1,\dots,b-1$ are polynomials in $x$.
This follows in almost the same way as in Case~3: We have
$i-e\ge e+c-b\ge 2c-b\ge0$ by our assumptions, 
 $e+c-j-1\ge e+c-b\ge 0$ if
$j<b$, $i+b-c-2e\ge0$, and $2e+c-j-1\ge
2e-b\ge 2c-b\ge0$ if $j<b+c$. This explains the term $(2b-2e-c)$ in
(3.70).

In order to explain the term $(b+c-2e)$ in (3.70), we claim that for
$s=0,1,\dots,b+c-2e-1$ we have
$$\align
&\sum _{i=0} ^{2e-b+s}\Bigg(
\sum_{k = 0}^{c - i-1}{{\left( -1 \right) }^{c + e + i + k + s+1}}
       {{c - i-1}\choose { b + c - 2 e + k - s-1}} 
\frac {  \left(  b + c - 2 e - s -1\right) !}
 {k!}\\
&\cdot
{{            \left(  b - e - s -1\right) !\,\left(  c + k - s -1\right) !\,
       \left(  b - e + k - s -1\right) !}\over 
     {\left(  2 b - 2 e - 2 s -2\right) !\,
       \left(  2 b - 2 e + k - 2 s -1\right) !}}\Bigg)
\cdot \big(\text {row $i$ of $\De_3(-e;b,c,e)$}\big)\\
&+
\sum _{i=e} ^{\min\{2e+c-b-1,b-s-1\}}
2\,{{\left( -1 \right) }^{c + s}}
{{ \left(  2e+c- b - i -1\right) !\,
     \left(  b + c - 2 e - s -1\right) !}\over 
   {\left(  b - i - s -1\right) !}}\\
&\hskip2cm\cdot
\frac {     ({ \textstyle 1 - e + i}) _{ b - i - s-1} }
 {     ({ \textstyle b - 2 e + i - s}) _{ b - i - s-1} }
\cdot \big(\text {row $i$ of $\De_3(-e;b,c,e)$}\big)\\
&+\chi(s\ge 2b-2e-c)
\sum _{i=2e+c-b} ^{b-s-1}
{{\left( -1 \right) }^{b + i + s+1}} 
{{     ({ \textstyle 1 + b - c - 2 e + i}) _{ 2 c - i - s-1} }\over 
   {\left(  b - i - s -1\right) !}}\\
&\hskip4cm
\cdot
\frac {     ({ \textstyle 1 - e + i}) _{ b - i - s-1} } 
{     ({ \textstyle b - 2 e + i - s}) _{ b - i - s-1} }
\cdot \big(\text {row $i$ of $\De_3(-e;b,c,e)$}\big)
\\
&\hskip4cm=0.\tag3.71
\endalign$$
We make a similar convention as
in Case~1 of how to understand $\De_3(x;b,c,e)$ in the case that $e\ge
b-c/2$, as we already did in Case~3.

Note that these are indeed $b+c-2e$ linear combinations of the rows,
which are linearly independent. 
The latter fact comes from the observation that for fixed $s$ 
the last nonzero coefficient in the linear combination (3.71) 
is the one for row $b-s-1$, regardless whether $s\ge 2b-2e-c$ or not. 

Because of the condition $s\le b+c-2e-1$, we have $2e-b+s\le c-1$, and
therefore the rows which are involved in the first sum in (3.71) are
from rows $0,1,\dots,c-1$, which form the top block in (3.38). 
The assumptions $e\ge c$ and $b\le 2c$ imply $2e-b+s\ge 0$, and so the
bounds for the sum are proper bounds. 
Because of $e\ge c$, the rows which are involved in the second sum in 
(3.71) are
from rows $c,c+1,\dots,2e+c-b-1$, which form the middle block in (3.38).
The bounds for the sum are proper, since by our assumptions we have
$$s\le b+c-2e-1\le b-e-1\le b-c-1\le c-1\le e-1,\tag3.72$$
and therefore $e\le b-s-1$ and $e\le 2e+c-b$ (including the
possibility that $b/2=c=e$, in which case the second sum in (3.71) is the
empty sum).
Finally, because of the condition $s\ge 0$, we have $b-s-1\le b-1$, and
therefore the rows which are involved in the third sum in (3.71) (if
existent) are
from rows $2e+c-b,2e+c-b+1,\dots,b-1$, which form the bottom block in (3.38).

Hence, in order to verify (3.71), we have to check
$$\align
&\sum _{i=0} ^{2e-b+s}\Bigg(
\sum_{k = 0}^{c - i-1}{{\left( -1 \right) }^{c + e + i + k + s+1}}
       {{c - i-1}\choose { b + c - 2 e + k - s-1}} 
\frac {  \left(  b + c - 2 e - s -1\right) !}
 {k!}\\
&\hskip1cm\cdot
{{            \left(  b - e - s -1\right) !\,\left(  c + k - s -1\right) !\,
       \left(  b - e + k - s -1\right) !}\over 
     {\left(  2 b - 2 e - 2 s -2\right) !\,
       \left(  2 b - 2 e + k - 2 s -1\right) !}}\Bigg)
\binom {c-e}{i-j+c}\\
&+\chi(s\ge 2b-2e-c)
\sum _{i=2e+c-b} ^{b-s-1}
{{\left( -1 \right) }^{b + c+e+i+j + s}} 
{{     ({ \textstyle 1 + b - c - 2 e + i}) _{ 2 c - i - s-1} }\over 
   {\left(  b - i - s -1\right) !}}\\
&\hskip4cm
\cdot
\frac {     ({ \textstyle 1 - e + i}) _{ b - i - s-1} } 
{     ({ \textstyle b - 2 e + i - s}) _{ b - i - s-1} }
\cdot \frac {(i-e)!\,(e+c-j-1)!} {(i-j+c)!}
\\
&\hskip4cm=0,\tag3.73
\endalign$$
which is (3.71) restricted to the $j$-th column, $j=c,c+1,\dots,b-1$,
(note that this is indeed the restriction of (3.71) to the $j$-th
column, $c\le j<b$, since, due to $0\le i-e\le i-c
\le i-b+c<i-j+c$, the entries $\binom {i-e} {i-j+c}$ in rows
$e,e+1,\dots,2e+c-b-1$ of $\De_3(-e;b,c,e)$ vanish in such a column),
and
$$\align
&\sum _{i=0} ^{2e-b+s}\Bigg(
\sum_{k = 0}^{c - i-1}{{\left( -1 \right) }^{c + e + i + k + s+1}}
       {{c - i-1}\choose { b + c - 2 e + k - s-1}} 
\frac {  \left(  b + c - 2 e - s -1\right) !}
 {k!}\\
&\hskip1cm\cdot
{{            \left(  b - e - s -1\right) !\,\left(  c + k - s -1\right) !\,
       \left(  b - e + k - s -1\right) !}\over 
     {\left(  2 b - 2 e - 2 s -2\right) !\,
       \left(  2 b - 2 e + k - 2 s -1\right) !}}\Bigg)
2\binom {b-2e}{i-j+b}\\
&+
\sum _{i=e} ^{\min\{2e+c-b-1,b-s-1\}}
2\,{{\left( -1 \right) }^{c + s}}
{{ \left(  2e+c- b - i -1\right) !\,
     \left(  b + c - 2 e - s -1\right) !}\over 
   {\left(  b - i - s -1\right) !}}\\
&\hskip2cm\cdot
\frac {     ({ \textstyle 1 - e + i}) _{ b - i - s-1} }
 {     ({ \textstyle b - 2 e + i - s}) _{ b - i - s-1} }
\binom {i+b-c-2e}{i-j+b}\\
&+\chi(s\ge 2b-2e-c)
\sum _{i=2e+c-b} ^{b-s-1}
2\,{{\left( -1 \right) }^{b + c+i+j + s}} 
{{     ({ \textstyle 1 + b - c - 2 e + i}) _{ 2 c - i - s-1} }\over 
   {\left(  b - i - s -1\right) !}}\\
&\hskip3cm
\cdot
\frac {     ({ \textstyle 1 - e + i}) _{ b - i - s-1} } 
{     ({ \textstyle b - 2 e + i - s}) _{ b - i - s-1} }
\cdot \frac {(i+b-c-2e)!\,(2e+c-j-1)!} {(i-j+b)!}
\\
&\hskip4cm=0,\tag3.74
\endalign$$
which is (3.71) restricted to the $j$-th column, $j=b,b+1,\dots,b+c-1$.

For the proof of (3.73) and (3.74) we follow the strategy of the
proofs of (3.40) and (3.41). That is, first the first terms in (3.73)
and (3.74) are recast, by replacing the binomials by the expansions
that were described in the proofs of (3.40) and (3.41), then interchanging
sums, evaluating the inner sums, etc. Eventually, it turns out that
the two terms on the left-hand side of (3.73) can be combined into
a single expression, namely into (3.42) with
$(1+b-c-2e+\de+i)_{b-i-s-1}$ replaced by
$(1+b-c-2e+\de+i)_{2c-i-s-1}$. The same arguments as in Case~4 of the
proof of
Lemma~1 then prove that this expression vanishes, since the crucial
inequalities $b+c-j-s-1\ge0$, $1-c+s\le0$, $b-j-1\ge0$ are also valid
here, thanks to (3.72). Similarly, it turns out, eventually, that the
three terms on the left-hand side of (3.74) can be combined into one
expression, namely into (3.45), again with 
$(1+b-c-2e+\de+i)_{b-i-s-1}$ replaced by
$(1+b-c-2e+\de+i)_{2c-i-s-1}$. Now the same arguments as in Case~4 of
the proof of 
Lemma~1 apply to prove that this expression vanishes, since the
crucial inequalities $2b-j-s-1\ge0$, $1-b+e+s\le0$, and $b+e-j-1\ge0$
are also valid here, thanks to (3.72) again.

This proves (3.73) and (3.74), and thus completes the proof that $(x+e)$
divides $\De'(x;b,c)$ with multiplicity $m(e)$ as given in (3.70).

\medskip
This finishes the proof of Lemma~2.\quad \quad \qed
\enddemo

\proclaim{Lemma 3}For any integer $c$, any nonnegative integer $n$,
and any number $X$, there holds
$$\det_{1\le i,j\le n}\(\binom {X}{i-j+c}\)=
\prod _{i=1} ^{n}\frac {(X+i-c)_c} {(i)_c}.
\tag3.75$$
\endproclaim

\demo{Proof}This is an ubiquitous determinant, and there are numerous
proofs of its evaluation, see e.g\. \cite{\GoJaAJ, Lemma~3.1;
\KratAM, computation on p.~189 with $\la_s=c$ and $a=X-\al+b$;
\ZeilBL} for some conceptual ones that also include generalizations.
\quad \quad \qed
\enddemo

\proclaim{Lemma 4}Let $b$ and $c$ be even integers, $b>c$. Then
$$\det_{c\le i,j<b}\pmatrix \dsize\binom {c-b/2-1/2}{2c-1-j}&i=c\\
\dsize\binom {c-b/2+1/2}{i-j+c}&i>c
\endpmatrix=0.
\tag3.76$$
\endproclaim
\demo{Proof}  Let us denote the determinant in (3.76) by $D$.
We claim that the rows of $D$ are linearly dependent. To be precise,
we claim that
$$\sum _{j=c} ^{b-1}(-1)^j\frac {(1-b+c)_{j-c}\, 
\(\frac {1} {2}-\frac {b} {2}\)_{j-c}} {(j-c)!\,(1-b)_{j-c}}\cdot
\big(\text {column $j$ of $D$}\big)=0.\tag3.77$$
To see this, we have to check
$$\sum _{j=c} ^{b-1}(-1)^j\frac {(1-b+c)_{j-c}\, 
\(\frac {1} {2}-\frac {b} {2}\)_{j-c}} {(j-c)!\,(1-b)_{j-c}}
\binom {c-b/2-1/2}{2c-1-j}=0,$$
which is (3.77) restricted to row $c$, and
$$\sum _{j=c} ^{b-1}(-1)^j\frac {(1-b+c)_{j-c}\, 
\(\frac {1} {2}-\frac {b} {2}\)_{j-c}} {(j-c)!\,(1-b)_{j-c}}
\binom {c-b/2+1/2}{i-j+c}=0,$$
which is (3.77) restricted to row $i$, $c<i<b$.
Equivalently, in terms of hypergeometric series, this means to check
$${{({ \textstyle {3\over 2} - {b\over 2}}) _{c-1} }\over 
   {\left( c-1 \right) !}}\,
{} _{3} F _{2} \!\left [ \matrix { 1 - c, 1 - b + c, {1\over 2} - {b\over
      2}}\\ { {3\over 2} - {b\over 2}, 1 - b}\endmatrix ; {\displaystyle
      1}\right ]=0
\tag3.78$$
and
$${{({ \textstyle {3\over 2} - {b\over 2} + c - i}) _{i} }\over
     {i!}}\,
{} _{3} F _{2} \!\left [ \matrix { {1\over 2} - {b\over 2}, 1 - b + c, -i}\\
      { 1 - b, {3\over 2} - {b\over 2} + c - i}\endmatrix ; {\displaystyle
      1}\right ]=0.
\tag3.79$$

Equation (3.78) follows from Watson's $_3F_2$-summation (cf\. 
\cite{\SlatAC, (2.3.3.13); Appendix (III.23)}),
$$
{} _{3} F _{2} \!\left [ \matrix { A, B, C}\\ { {{1 + A + B}\over 2}, 2
   C}\endmatrix ; {\displaystyle 1}\right ]  =
  \frac {\Gamma \( {1\over 2}\) \Gamma \( {1\over 2} + C\) \Gamma \( {1\over 2} + {A\over 2} +
   {B\over 2}\) \Gamma \( {1\over 2} - {A\over 2} - {B\over 2} + C\)} 
  {\Gamma \( {1\over 2} + {A\over
   2}\) \Gamma \( {1\over 2} + {B\over 2}\) \Gamma \( {1\over 2} - {A\over 2} + C\) \Gamma \( {1\over 2} -
   {B\over 2} + C\)}.
\tag3.80$$
For, the term $\Ga(1/2+A/2)$ in the denominator of the right-hand
side of (3.80) implies that the $_3F_2$-series on the left-hand side will
vanish whenever $A$ is an odd negative integer. This is exactly the
case for the $_3F_2$-series in (3.78), where $A=1-c$ with $c$ being
even by assumption.

Equation (3.79) follows from the Pfaff--Saalsch\"utz summation (3.8).
For, a straight-forward
application of formula (3.8) gives for the $_3F_2$-series in (3.79) the
expression
$${{({ \textstyle {1\over 2} - {b\over 2}}) _{i} \,({ \textstyle -c}) _{i} \,
     ({ \textstyle {3\over 2} - {b\over 2} + c - i}) _{i} }\over 
   {i!\,({ \textstyle 1 - b}) _{i} \,
     ({ \textstyle -{1\over 2} + {b\over 2} - c}) _{i} }}.
$$
In the numerator of this expression there appears the term $(-c)_i$,
which vanishes because $i>c$. This completes the proof of the Lemma.\quad \quad
\qed

\enddemo

\enddocument